\documentclass[12pt]{amsart}

\usepackage{amsmath}
\usepackage{amsthm}

\usepackage{amsfonts}
\usepackage{amssymb}
\usepackage{graphicx}
\newcommand\cyr{%
\renewcommand\rmdefault{wncyr}%
\renewcommand\sfdefault{wncyss}%
\renewcommand\encodingdefault{OT2}%
\normalfont
\selectfont}
\DeclareTextFontCommand{\textcyr}{\cyr}

\usepackage{tikz-cd}
\usepackage{ amssymb }
\usepackage{tikz}
\usepackage{stix}
\usepackage{enumitem}
\usetikzlibrary{matrix,arrows}\usepackage[OT2,OT1]{fontenc}

\newcommand{\cB}{ \textcyr{B}} 
\newcommand{\ocB}{ \overline{\textcyr{B}}} 
\newcommand{\G}{\mathcal{G}}
\newcommand{\cT}{\mathcal{T}}
\newcommand{\cG}{\mathcal{G}}
\newcommand{\cN}{\mathcal{N}}
\newcommand{\cH}{\mathcal{H}}
\newcommand{\cJ}{\mathcal{J}}
\newcommand{\cK}{\mathcal{K}}
\newcommand{\cM}{\mathcal{M}}
\newcommand{\N}{\mathbb{N}}
\newcommand{\Z}{\mathbb{Z}}

\newcommand{\Gal}{\text{Gal}}

\newcommand{\noi}{\noindent}

\newcommand{\bs}{\bigskip}

\newcommand{\im}{\text{im}}
\newcommand{\Aut}{\text{Aut}}

\newcommand{\Br}{\text{Br}}

\newcommand{\res}{\text{res}}
\newcommand{\cor}{\text{cor}}
\newcommand{\Ind}{\text{Ind}}

\newcommand{\tF}{\widetilde{F}}
\newcommand{\tE}{\widetilde{E}}

\input xy
\xyoption{all}

\usepackage[top =.6in, bottom=.6in, left=.5in, right=.7in]{geometry}

\newtheorem{theorem}{Theorem}[section]

\newtheorem{corollary}[theorem]{Corollary}
\newtheorem{definition}[theorem]{Definition}

\newtheorem{lemma}[theorem]{Lemma}

\newtheorem{remark}[theorem]{Remark}

\input xy
\xyoption{all}
 
 \LARGE

\title{Cohomological Kernels for Cyclic by Cyclic Semi-Direct Product Extensions}

\author{ Nathan Schley}

\begin{document}

\maketitle
\vspace{-25pt}
\begin{center}
Department of Mathematics, University of California,\\
Santa Barbara, California, USA 93106\\  email: schley@math.ucsb.edu\\ 
\end{center}

\begin{abstract}
\noi Let $F$ be a field and $E$ an extension of $F$ 
with $[E:F]=d$ where the characteristic of $F$ is 
zero or prime to $d$. We assume $\mu_{d^2}\subset F$ where $\mu_{d^2}$ are the $d^2$th roots of unity. This paper studies the problem of determining 
the cohomological kernel $H^n(E/F):=\ker(H^n(F,\mu_d) \rightarrow
H^n(E,\mu_d))$ (Galois cohomology with
coefficients in the $d$th roots of unity)
when the Galois closure of $E$ is a semi-direct product of cyclic groups. The main result is a six-term exact sequence determining the kernel as the middle map and is based
on tools of Positelski \cite{Positselski}. When
$n=2$ this kernel is the relative 
Brauer group ${\rm Br}(E/F)$, the classes of 
central simple algebras in the Brauer group of $F$
split in the field $E$.  
The work of Aravire and Jacob (2008, 2018) \cite{AJ08} \cite{AJO} which calculated the groups 
$H^n_{p^m}(E/F)$ in
the case of semidirect products of cyclic groups
 in characteristic $p$, provides motivation
for this work. 
\end{abstract}  
 \bs

 \Large
 \noi {\bf Introduction}
 \normalsize
 
 \bs
 
This paper studies cohomological kernels of field extensions.
Historically, such kernels have played a key role in the computation of
relative Brauer groups, although the results presented here
apply more generally to  higher cohomology.  
Similarly, the computation of analogous
cohomological kernels have played important roles in the development of
the algebraic theory of quadratic forms.

In 1980, the work of Merkurjev and Suslin on the  conjecture
of Albert showed that in the presence of roots of unity, the Brauer group is generated
by classes of cyclic algebras.  
In terms of cohomology, this means that the
cup product map $H^1(F,\mu_d)\times H^1(F,\mu_d)\rightarrow H^2(F,\mu_d)\cong
{\rm Br}_d(F)$ is surjective. 
Their work was dependent upon detailed analyses of the K-theory
of Severi-Brauer varieties and the relationship between the Milnor K-theory of a field
and its Galois cohomology.   

In 2005, Positselki studied cohomological kernels of biquadratic extensions and certain degree $8$ extensions \cite{Positselski} using a four-term exact sequence of Galois group modules 
\[
0 \longrightarrow M_1 \longrightarrow M_2 \longrightarrow M_3 \longrightarrow M_4 \longrightarrow 0
\]
with homotopy maps and some other properties to produce a six-term exact sequence of cohomology. 
Prior to his work, it was known by work analyzing the Witt ring
that if $E=F(\sqrt{a},\sqrt{b})$ is biquadratic then the kernel
of the map $H^2(F,\Z/2\Z)\rightarrow H^2(E,\Z/2Z)$ is 
generated by the images of ``expected elements''
$(a)\smile(x)$ and $(b)\smile(y)$ for $x,y\in F$ 
\cite{ElmanLam} and that
the analogue of this expected result 
for triquadratic extensions was false
\cite{AmitsurRowenTignol}. 
The question of determining the kernel of
$H^n(F,\Z/2\Z)\rightarrow H^n(E,\Z/2\Z)$ for $n\geq3$
in the separable biquadratic case
was considered by a number of researchers (Merkurjev, Tignol, Kahn), and it is this problem that  Positselski solved
with his tools. Positselski's tools also applied to
dihedral extensions of degree $8$, indicating the applicability
of these techniques to the non-Galois case.  Characteristic
$p$ versions of Positselski's machinery  have been 
constructed by Aravire and Jacob (2012), for the separable
biquadratic case and the dihedral and quaternion cases
in characteristic 2 (2016), and more generally for the 
cyclic by cyclic semi-direct product cases in characteristic
$p>2$ by Aravire-Jacob-O'Ryan (2018).  It is these
latter constructions that this paper generalizes to the
case where the characteristic is prime to the field degree $d$
and the $d^2$th roots of unity are present in the field.  

The key to this work is determining the appropriate
modules $M_3$ and $M_4$ (see below for the set-up) and establishing the
requisite homotopies necessary to apply Positselski's
tools.  
This is spelled out in Section \ref{General Setup Section}. Sections \ref{Arason Section} and  \ref{Positselski Section}
develop the background as well as  provide details 
necessary for the application of Positselski's results that
are not clearly spelled out in his paper. In particular,
the ``connecting map'' $\eta: H^n(
\G
,M_4)
\rightarrow H^{n+1}(\G,M_1)$ needs to be
carefully computed. 
Section 5 covers the case of a cyclic extension and uses this machinery to prove a well-known result from Hilbert 90 as a way to get a sense of how the machinery works. 
Section 5 covers a dihedral extension.
Section 6 covers the more general case of an extension whose Galois group is a semi-direct product with certain conditions. 
The final section covers some cohomological interpretations of these results. 

\subsection{Notation and Further Background.} 
Let $F$ be a field, $d \in \N$ with $d> 1$, 
we will assume that 
${\rm char}(F)=0$ or (${\rm char}(F),d) = 1$, 
and that $\mu_{d^2} \subseteq F$, where $\mu_{d^2}$ are the ${d^2}$ distinct ${d^2}$th roots of unity. 
Let $F_{\text{sep}}$ denote the separable closure of $F$, 
$\cG = 
\Gal(F_{\text{sep}}/F)$, 
and $H^n(\cG, M)$, the $n$th cohomology groups 
for any $\Z[G]$-module $M$ \cite{GalCoh}. 
Let $E/F$ be an extension of degree $d$, $\cH \subseteq \cG$ be $\Gal(F_{\text{sep}}/E)$. We will also use the  notation  $H^n(F, M) = H^n(\cG,M)$, so $H^n(E,M) = H^n(\cH,M)$.
We also denote by $H^n(E/F, M):=
\ker(H^m(F,M)\rightarrow H^n(E,M))$.

The groups $H^0(F, \mu_d)$ and $H^1(F, \mu_d)$ have an interpretation from Kummer theory. Consider the following short exact sequence of $\Z[\cG]$-modules 
\[
0 \longrightarrow \mu_d \overset{\subseteq}{\longrightarrow} F_{\text{sep}}^{\times}  \overset{\cdot d}{\longrightarrow} F_{\text{sep}}^{\times}  \longrightarrow 0
\]
where the second map is multiplication by $d$ over the $\Z[\cG]$-modules. It is surjective because $F_{\text{sep}}$ is separably closed. This short exact sequence of $\Z[\cG]$-modules yields a long exact sequence of cohomology \cite{HomAlg}. 
\[
\begin{tikzpicture}[descr/.style={fill=white,inner sep=1.5pt}]
        \matrix (m) [
            matrix of math nodes,
            row sep=3em,
            column sep=3em,
            text height=1.5ex, text depth=0.25ex
        ]
        { 0 & H^0(\cG, \mu_d) & H^0(\cG, F_{\text{sep}}^{\times} ) & H^0(\cG, F_{\text{sep}}^{\times} ) &  \\ 
          \ & H^1(\cG, \mu_d) & H^1(\cG, F_{\text{sep}}^{\times} ) & H^1(\cG, F_{\text{sep}}^{\times} ) &  \\ 
          \ & H^2(\cG, \mu_d) & H^2(\cG, F_{\text{sep}}^{\times} ) & H^2(\cG, F_{\text{sep}}^{\times} ) &  \\ 
        \ & H^3(\cG, \mu_d) & \cdots \ \ \ \ \ \ \ \ \ \  \\ 
        }; 

        \path[overlay, font=\scriptsize,>=latex]
        (m-1-1) edge [->] (m-1-2)
        (m-1-2) edge [->]  node[descr,yshift=1.5ex] {$\subseteq $} (m-1-3) 
        (m-1-3) edge [->]  node[descr,yshift=1.5ex] {$\cdot d$} (m-1-4) 
        (m-2-2) edge [->]  node[descr,yshift=1.5ex] {$\subseteq $} (m-2-3) 
        (m-2-3) edge [->]  node[descr,yshift=1.5ex] {$\cdot d$} (m-2-4) 
        (m-3-2) edge [->]  node[descr,yshift=1.5ex] {$\subseteq$} (m-3-3) 
        (m-3-3) edge [->]  node[descr,yshift=1.5ex] {$\cdot d$} (m-3-4) 
        (m-4-2) edge [->]   (m-4-3) 
        (m-1-4) edge[out=355,in=175,->] node[descr,yshift=0.3ex] {$\partial$} (m-2-2)
        (m-2-4) edge[out=355,in=175,->] node[descr,yshift=0.3ex] {$\partial$} (m-3-2)
        (m-3-4) edge[out=355,in=175,->] node[descr,yshift=0.3ex] {$\partial$} (m-4-2)
;
\end{tikzpicture} 
\] 
Note that Fix$_\cG(\mu_d) = \mu_d$ and Fix$_\cG(F_{\text{sep}}^{\times} ) = F^{\times} $ by Galois theory. 
Furthermore, $H^1(\cG, F_{\text{sep}}^{\times} )$ is trivial by the cohomological version of Hilbert's Theorem 90. This information gives the long exact sequence.  
\[
\begin{tikzpicture}[descr/.style={fill=white,inner sep=1.5pt}]
        \matrix (m) [
            matrix of math nodes,
            row sep=3em,
            column sep=3em,
            text height=1.5ex, text depth=0.25ex
        ]
        { 0 & \mu_d & F^{\times}  & F^{\times}  &  \\ 
          \ & H^1(\cG, \mu_d) & 0 & 0 &  \\ 
          \ & H^2(\cG, \mu_d) & H^2(\cG, F_{\text{sep}}^{\times} ) & H^2(\cG, F_{\text{sep}}^{\times} ) &  \\ 
        \ & H^3(\cG, \mu_d) & \cdots \ \ \ \ \ \ \ \ \ \  \\ 
        };

        \path[overlay, font=\scriptsize,>=latex]
        (m-1-1) edge [->] (m-1-2)
        (m-1-2) edge [->]  node[descr,yshift=1.5ex] {$\subseteq $} (m-1-3) 
        (m-1-3) edge [->]  node[descr,yshift=1.5ex] {$\cdot d$} (m-1-4) 
        (m-2-2) edge [->]  node[descr,yshift=1.5ex] {$\subseteq $} (m-2-3) 
        (m-2-3) edge [->]  node[descr,yshift=1.5ex] {$\cdot d$} (m-2-4) 
        (m-3-2) edge [->]  node[descr,yshift=1.5ex] {$\subseteq$} (m-3-3) 
        (m-3-3) edge [->]  node[descr,yshift=1.5ex] {$\cdot d$} (m-3-4) 
        (m-4-2) edge [->]   (m-4-3) 
        (m-1-4) edge[out=355,in=175,->] node[descr,yshift=0.3ex] {$\partial$} (m-2-2)
        (m-2-4) edge[out=355,in=175,->] node[descr,yshift=0.3ex] {$\partial$} (m-3-2)
        (m-3-4) edge[out=355,in=175,->] node[descr,yshift=0.3ex] {$\partial$} (m-4-2)
;
\end{tikzpicture} 
\] 
In particular we have the following three results: 
\begin{enumerate}
\item $H^0(F, \mu_d) \cong \Z/dZ$. 
\item $H^1(F, \mu_d) \cong F^{\times} /F ^{\times d}$ and 
\item $H^2(F, \mu_d)$ is the $d$-torsion of $H^2(F, F_{\text{sep}}^{\times} ).$
\end{enumerate} 
For $a\in F^{\times}$ we use 
$(a)\in H^1(F, \mu_d)$ to denote
the class that  $aF^{\times d}\in F^{\times} /F^{\times ^d}$
corresponds to in the second identification.
Since $H^2(F, F_{\text{sep}}^{\times} ) \cong \Br(F)$ is the Brauer group (the cohomological Brauer group and the Brauer
group agree for fields), the third result will be of particular importance  because it means  
$H^2(F,  \mu_d)$ picks out the $d$-torsion in  $\Br(F)$.
We  use $\smile$ to denote the cup product:
$\smile:H^r(F,  \mu_d)\times H^s(F,  \mu_d)\rightarrow
H^{r+s}(F,  \mu_d)$, which makes sense in our
context because $\mu_d\subset F$ and therefore
has trivial $\G$-action so $\mu_d^{\otimes2}\cong\mu_d$ as $\G$-modules.


\subsection{The Problem Studied} 
The problem studied in this paper is that of determining the kernels of scalar extension
(restriction in group cohomology),
\[
\res_{E/F}:H^n(F,\mu_d) \ 
,{\longrightarrow}\, 
H^n(E, \mu_d)
\]
for various extension fields $E/F$ of degree $d$. 
The case where $E/F$ is cyclic Galois is basic. 
In the cyclic case when $n=2$, if the Brauer class of an $F$-division algebra 
$D$ of index $d$ lies in $H^2(E/F,\mu_d)$, then  this $D$ is a cyclic algebra (with maximal subfield $E$.)  
More specifically, we know that when $E=F(\sqrt[d]{a})$ (recall $\mu_d\subset F$) we  have $H^2(E/F,\mu_d)=(a)\smile H^1(F,\mu_d)$. 
Theorem \ref{Connecting Map Cyclic 6-term Sequence Theorem} shows that in the cyclic case this is valid for all $n$, 
namely $H^{n+1}(E/F,\mu_d)=(a)\smile H^n(F,\mu_d)$, a result already known \cite{Arason}.
The next cases generalize this situation, where either the Galois closure of  $E$ is dihedral or $E$
is an extension of degree $d$ that becomes a cyclic extension
when $F$ is extended by a cyclic extension of degree
prime to $d$. In this latter case the Galois
group of the Galois closure of $E$ is a cyclic by cyclic semi-direct product. 
In these latter cases one cannot describe the cohomological kernel as a
cup product by a class $(a)$ (indeed,
$H^1(E/F,\mu_d)=0$), but one does have the
connecting map $\eta$ from Positselski's theory to capture the kernel. 

In order to compute these kinds of kernels, Positselski used 4-term exact sequences of Galois group modules with homotopy maps
\[
0 \longrightarrow M_1 \longrightarrow M_2 \longrightarrow M_3 \longrightarrow M_4 \longrightarrow 0
\]
to produce a six-term exact sequence of cohomology. Here
with $M_1\cong\mu_d$ and $M_2$  an appropriately selected 
induced module with $H^n(\G,M_2)\cong H^n(E,\mu_d)$, so that the six-term sequence can be used to compute the cohomological
kernel $H^n(E/F,\mu_d)$. 
Aravire and Jacob \cite{AJ08} and Aravire, Jacob and O'Ryan \cite{AJO} have developed a variant of this machinery to compute cohomologial kernels in characteristic $p$ for $E/F$ of prime degree $d = p$ with Galois closure having Galois group a semidirect product of two cyclic groups of order $p$ and $s$, where $s | (p-1)$. 
This paper gives an analogous result when $E/F$ is degree $d$, $F$ has characteristic prime to $d$, and the Galois closure of $E/F$ has Galois group a semidirect product of cyclic groups of order $d$ and $s$, with $s | \phi(d)$ (the Euler $\phi$-function) and $\Z/s\Z$ acting faithfully on $\Aut(\Z/d\Z)$. 

We recall from applications that $H^2(\cG, \mu_d) = Br_d(F)$, and standard notation from this subject will be used.

\section{Arason's Theorem} \label{Arason Section} 

In his paper, Arason \cite{Arason} proved that the
third cohomological invariant, $e_3$ of quadratic forms is
well-defined.  To accomplish this, he determined the
cohomological kernel of a quadratic extension away from 
characteristic two (an equivalent result in group cohomology
was proved independently by D.L. Johnson \cite{Johnson} at the same time.)
We discuss Arason's results here because the approach he
took provides a model for understanding the work of Positselski,
and the computation of his connecting map lays a conceptual
framework for the computation of Positselksi's connecting map
$\eta$, see Example 3.7 of \cite{HW} for further discussion. 

Let $F$ be a field, char$(F) \neq 2$, $E = F(\sqrt{a})$ a quadratic extension, and $F_{\text{sep}}$ the separable closure of $F$. 
Let $\mu_2$ be the square roots of unity $\pm 1$; clearly $\mu_2 \subseteq F$. 
The result of Arason \cite{Arason} is the following, a result which has a critical role in the algebraic theory of quadratic forms. 
It is a cohomological analogue of an exact sequence for the Witt ring (see \cite{LamsBook} chap. 7 Sec. 3). 
We also sketch the proof. 

\begin{theorem} \label{Arason's Theorem} 
Let $F$ be a field, char$(F) \neq 2$, $E = F(\sqrt{a})$ a quadratic extension.
There is a long exact restriction/corestriction sequence
\[
\begin{tikzpicture}[descr/.style={fill=white,inner sep=1.5pt}]
        \matrix (m) [
            matrix of math nodes,
            row sep=3em,
            column sep=3em,
            text height=1.5ex, text depth=0.25ex
        ]
        { 0 & H^0(F, \mu_2) & H^0(E, \mu_2) & H^0(F, \mu_2) &  \\ 
          \ & H^1(F, \mu_2) & H^1(E, \mu_2) & H^1(F, \mu_2) &  \\ 
          \ & H^2(F, \mu_2) & H^2(E, \mu_2) & H^2(F, \mu_2) &  \\ 
        \ & \cdots \ \ \ \ \ \ \ \ \ \  \\ 
        };

        \path[overlay, font=\scriptsize,>=latex]
        (m-1-1) edge [->] (m-1-2)
        (m-1-2) edge [->]  node[descr,yshift=1.5ex] {$\res$} (m-1-3) 
        (m-1-3) edge [->]  node[descr,yshift=1.5ex] {$\cor$} (m-1-4) 
        (m-2-2) edge [->]  node[descr,yshift=1.5ex] {$\res$} (m-2-3) 
        (m-2-3) edge [->]  node[descr,yshift=1.5ex] {$\cor$} (m-2-4) 
        (m-3-2) edge [->]  node[descr,yshift=1.5ex] {$\res$} (m-3-3) 
        (m-3-3) edge [->]  node[descr,yshift=1.5ex] {$\cor$} (m-3-4) 
        (m-1-4) edge[out=355,in=175,->] node[descr,yshift=0.3ex] {$\partial$} (m-2-2)
        (m-2-4) edge[out=355,in=175,->] node[descr,yshift=0.3ex] {$\partial$} (m-3-2)
        (m-3-4) edge[out=355,in=175,->] node[descr,yshift=0.3ex] {$\partial$} (m-4-2)
;
\end{tikzpicture} .
\] 
where the connecting map $\partial$ is the cup product with the character function $\chi_K\in H^1(F,\mu_2)$, which corresponds
to the class $(a)\in F^{\times}/F^{\times2}$.
\end{theorem} 
\noi \underline{\bf Proof}: (Sketch) This long exact sequence is induced from a short exact sequence of $\cG$-modules 
\[
0 \longrightarrow \mu_2 \longrightarrow \Ind _{\cH}^{\G}(\mu_2) \longrightarrow \mu_2 \longrightarrow 0
\]
that will give the restriction and corestriction on cohomology once we replace $H^n(\cG, \Ind_\cH^\cG(\mu_2))$ with $H^n(\cH,\mu_2)$ using the Shapiro isomorphism.
We note that as $\G = \Gal(F_{\text{sep}}/F)$ and $\cH = \Gal(F_{\text{sep}}/E)$, 
we have $H^n(\cG, \mu_2) = H^n(F, \mu_2)$ and $H^n(\cH, \mu_2) = H^n(E, \mu_2)$. 

For our computations, we identify $Ind_\cH^\cG(\mu_2)$ with $\Z/2\Z \oplus \Z/2\Z$ and use the sequence 
\[
0 \longrightarrow \Z/2\Z \longrightarrow \Z/2\Z \oplus \Z/2\Z \longrightarrow \Z/2\Z \longrightarrow 0.
\]
Here $\Z/2\Z$ is a trivial $\G$-module since the field elements $(\pm 1)$ are fixed by the Galois group $\cG$, 
while the induced module has a $\G$-action of permuting the two entries for any $g \notin \cH$, and a trivial $\G$-action for any $h \in \cH$. 
Then the maps are the diagonal map $1 \mapsto 1 \oplus 1$ and the trace $x \oplus y \mapsto x+y$, respectively, 
which happen to be the only non-zero choices of $\Z[\cG]$-module homomorphisms.   

We need to show that the two $\G$-maps in the short exact sequence induce maps that commute with the restriction and corestriction through the Shapiro isomorphism, 
and we need to show that the snake-lemma connecting map $\partial$ is the cup product with the character function $\chi_K$. 
This latter fact is pulled out as Theorem \ref{Arason2}
below. The  Shapiro isomorphism gives the following,
\[
\begin{tikzpicture}[descr/.style={fill=white,inner sep=1.5pt}]
        \matrix (m) [
            matrix of math nodes,
            row sep=3em,
            column sep=3em,
            text height=1.5ex, text depth=0.25ex
        ]
        { \ & H^n(F, \Ind_\cH^\G(\Z/2\Z))   \\ 
        H^n(F, \Z/2\Z) \\ 
        \ & H^n(K, \Z/2\Z)  \\ 
        };

        \path[overlay, font=\scriptsize,>=latex]
        (m-2-1) edge [->] node[descr,yshift=1.5ex] {$\Delta$}(m-1-2)
        (m-3-2) edge [->] node[descr,xshift=1.7ex] {$\overset{\text{Shap.}}{\cong}$} (m-1-2) 
        (m-2-1) edge [->] node[descr,yshift=1.5ex] {$\res$} (m-3-2) 
;
\end{tikzpicture} 
\] 
The composition of the restriction with the Shapiro map is induced by the identity map on $\Z/2\Z$, followed by the diagonal map, which agrees with the induced map on the top. \\ 
\[
\begin{tikzpicture}[descr/.style={fill=white,inner sep=1.5pt}]
        \matrix (m) [
            matrix of math nodes,
            row sep=3em,
            column sep=3em,
            text height=1.5ex, text depth=0.25ex
        ]
        { H^n(F, \Ind_\cH^\G(\Z/2\Z))   \\
        \ & H^n(F, \Z/2\Z) \\ 
        H^n(K, \Z/2\Z)  \\ 
        };

        \path[overlay, font=\scriptsize,>=latex]
        (m-1-1) edge [->] node[descr,yshift=1.5ex] {$\text{Tr}$}(m-2-2)
        (m-3-1) edge [->] node[descr,xshift=1.7ex] {$\overset{\text{Shap.}}{\cong}$} (m-1-1) 
        (m-3-1) edge [->] node[descr,yshift=1.5ex] {$\cor$} (m-2-2) 
;
\end{tikzpicture} 
\] 
The Shapiro homomorphism on $\Z[G]$-modules is the diagonal map. When composed with the trace, 1 is sent to the sum of 1 over every $\cH$-coset. 
This is also the map that induces the corestriction, and therefore the diagrams commute on cohomology. 
This concludes the proof sketch of Theorem \ref{Arason's Theorem}. \hfill $_\square$ \vspace{10pt}

It remains to interpret the connecting map in the long exact sequence. Arason describes this as the cup product with the character function, and this will be shown next by direct computation.

\begin{theorem} \label{Arason2}
The connecting map $\partial$ in Arason's Theorem 
\ref{Arason's Theorem} sends a cocycle $\sigma$ to $\chi_E \smile \sigma$, the cup product with the character function $\chi_E$ defined to be the composite $\cG \longrightarrow \Gal(E/F) \longrightarrow \mu_2$, where the latter map is the unique non-trivial homomorphism from $\Gal(E/F) \longrightarrow \mu_2$. 
\end{theorem} 
\noi \underline{Proof}: We will compute the connecting map directly. 
Let $\sigma \in Z^{n-1}(F, \Z/2\Z).$ 
We pick the lifting $\ell$ of the trace map that sends $1$ to $1 \oplus 0$. 
So $\ell(\sigma)(g_1, \ldots, g_{n-1}) = (\sigma(g_1, \ldots, g_{n-1})) \oplus 0$. 
Note that though $\ell$ is not a $\G$-map (if it were, the connecting map would have to be zero), $\ell$ is an abelian group homomorphism. 
And we will use this fact in the next step of the computation. 

The next step after the lifting $\ell$ is the chain map $\delta$, which can be computed using the bar resolution. 
\[
    \begin{split}
        \delta(\ell(\sigma))(g_1, \ldots, g_n) &= g_1 \cdot \ell(\sigma(g_2, \ldots, g_n)) - \ell(\sigma(g_1 g_2, \ldots, g_n)) + \cdots + (-1)^n \ell(\sigma(g_1, \ldots, g_{n-1})) \\ 
        &= g_1 \cdot \ell(\sigma(g_2, \ldots, g_n)) - \ell(g_1 \cdot \sigma(g_2, \ldots, g_n)) + \ell(\delta(\sigma)(g_1, \ldots, g_n)) \\ 
        &= g_1 \cdot \ell(\sigma(g_2, \ldots, g_n)) - \ell(g_1 \cdot \sigma(g_2, \ldots, g_n))
    \end{split}
\] 
This simplification takes advantage of the fact that $\ell$ is an abelian group homomorphism as well as the fact that $\sigma$ is a cocycle, so $\delta$ applied to $\sigma$ is zero, which gets sent by $\ell$ to 0 as well. 
Now, because our modules are 2-torsion and $\Z/2\Z$ has a trivial $\G$-action, 
\[
g_1 \cdot \ell(\sigma(g_2, \ldots, g_n)) - \ell(g_1 \cdot \sigma(g_2, \ldots, g_n)) = (g_1 + 1) \ell(\sigma(g_2, \ldots, g_n)).
\]
Furthermore, $(h+1)$ annihilates all of $\Z/2\Z \oplus \Z/2\Z$ whenever $h \in \cH$, and for every $g_1\not\in \cH$, $(g_1 + 1)(1\oplus0) = 1\oplus1$. 
Thus 
\[
(g_1 + 1) \ell(\sigma(g_2, \ldots, g_n)) = \begin{cases} \sigma(g_2, \ldots, g_n) \oplus \sigma(g_2, \ldots, g_n) \ &\text{if} \ g_1 \notin \cH \\ 0 \oplus 0 \ &\text{if} \ g_1 \in \cH \\
\end{cases}\]
This is the same as the diagonal image of the function 
\[
\chi_\cH(g_1) \cdot \sigma(g_2, \ldots, g_n)) = \begin{cases} \sigma(g_2, \ldots, g_n) \ &\text{if} \ g_1 \notin \cH \\ 
0 \ &\text{if} \ g_1 \in \cH \\
\end{cases},
\]
which completes the computation of the connecting map $\partial$. Therefore $\partial(\sigma) = \chi_\cH \smile \sigma$, the cup product with the character function $\chi_\cH$. This completes the proof of Theorem \ref{Arason2}. \hfill $_\square$ \vspace{10pt}

These ideas in the proof of Theorem \ref{Arason2} will be generalized in the sections that follow when we compute  Positselski's connecting map $\eta$.

\section{Positselski's 6-Term Cohomological Sequence } \label{Positselski Section} 
If we attempt to construct a short exact sequence in the fashion of Arason's theorem for a cyclic extension of degree $d > 2$, 
then the dimensions of the modules over $\Z/d\Z$ would be 1, d, 1, which makes exactness (and hence this approach) impossible. 
However, with the right machinery we can obtain 6-term exact sequences with a connecting map whose image is the cohomological kernel. 
Positselski's theorem offers exactly this machinery. The following definition collects the hypotheses needed for this machinery \cite{Positselski}.

\begin{definition} \label{Positselski's Hypotheses} {\bf Positselski's Hypotheses}: Let $\G$ be a pro-finite group, let $d, n \in \Z$ with $d \geq 2$ and $n \geq 0$, and let 
\[
0 \longrightarrow A_2 \overset{d_1}{\longrightarrow} B_2 \overset{d_2}{\longrightarrow} C_2 \overset{d_3}{\longrightarrow} D_2 \longrightarrow 0
\]
be a 4-term exact sequence of free $\Z/d^2\Z$-modules with a discrete action of $\G$. Let $h_1, h_2, h_3$ be homotopy maps 
\[
A_2 \overset{h_1}{\longleftarrow} B_2 \overset{h_2}{\longleftarrow} C_2 \overset{h_3}{\longleftarrow} D_2.
\]
Furthermore, let $A_1, B_1, C_1, D_1$ be the $d$-torsion of $A_2, B_2, C_2, D_2$ (respectively), let $\overline{A_2}:=A_2/A_1$, noting that $A_1\cong \overline{A_2}$ as ${\cG}$-modules.

If the homotopy maps satisfy the ``prism'' condition, that $d_i   h_i + h_{i+1}   d_{i+1} = d \cdot \text{id}$ for all $i \in \{0, 1, 2, 3 \}$, and if the Bockstein maps are 0 for all 4 modules, i.e. the kernel/cokernel short exact sequences 
\[
0 \longrightarrow A_1 \stackrel{\subseteq}{\longrightarrow} A_2 \overset{\cdot d}{\longrightarrow} \overline{A_2} \longrightarrow 0 
\]
have 0 connecting maps $\partial: H^n(\cG,\overline{A_2}) \longrightarrow H^{n+1}(\cG, A_1)$ and the same thing applies for $B, C$ and $D$ alike, then the four modules and their associated maps satisfy the Positselski Hypotheses.  
\end{definition} 

Given the definition we can now give Positselksi's main result
\cite{Positselski} Theorem 6.

\begin{theorem} \label{Positselski's Theorem} (Positselski) Given a short exact sequence of $\G$-modules that satisfies the Positselski hypotheses, there is a 6-term exact sequence of $\G$-cohomology 
\[
\begin{tikzpicture}[descr/.style={fill=white,inner sep=1.5pt}]
        \matrix (m) [
            matrix of math nodes,
            row sep=3em,
            column sep=3em,
            text height=1.5ex, text depth=0.25ex
        ]
        { H^n(B_1 \oplus D_1) & H^n(C_1) & H^n(D_1) &  \\ 
        H^{n+1}(A_1) & H^{n+1}(B_1) & H^{n+1}(A_1 \oplus C_1) \\ 
        };

        \path[overlay, font=\scriptsize,>=latex]
        (m-1-1) edge [->]  node[descr,yshift=1.5ex] {$d_2 + h_3$} (m-1-2) 
        (m-1-2) edge [->]  node[descr,yshift=1.5ex] {$d_3$} (m-1-3) 
        (m-2-1) edge [->]  node[descr,yshift=1.5ex] {$d_1$} (m-2-2) 
        (m-2-2) edge [->]  node[descr,yshift=1.5ex] {$h_1 \oplus d_2$} (m-2-3) 
        (m-1-3) edge[out=355,in=175,->] node[descr,yshift=0.3ex] {$\eta$} (m-2-1)
        
;
\end{tikzpicture} 
\] 
with the connecting map $\eta$ defined in Lemma \ref{Positselski Eta} below. 
\end{theorem} 
We will spend the remainder of this section setting up the framework for this 6-term exact sequence as an exposition 
(and a few extra details) of how Positselski builds the framework in his paper \cite{Positselski} as well as providing the tools we
need for computing the map $\eta$ in our applications. 
This process begins with the following two technical lemmas and ends with a proof of Theorem \ref{Positselski's Theorem}. 
This first lemma is based on Positselski's Lemma 5 \cite{Positselski}.

\begin{lemma} \label{Positselski Prism} (Posetselski)
Let $X$ be a $\G$-module that is also a free $\Z/d^2\Z$-module. 
Then we have a natural short exact sequence: 
$0 \longrightarrow X_1 \longrightarrow X_2 \longrightarrow \overline{X_2} \longrightarrow 0$  of $\G$-modules, 
where $X_2$ is $d^2$-torsion, $X_1$ is the $d$-torsion of $X_2$, and $\overline{X_2}$ is the quotient group, 
which must also be $d$-torsion and isomorphic to $X_1$. Suppose $Y, Z$ are short exact sequences defined in the same way, and that 
\[
0 \longrightarrow X \overset{d_1}{\longrightarrow} Y \overset{d_2}{\longrightarrow} Z \longrightarrow 0
\]
is a short exact sequence of these short exact sequences with homotopy maps 
\[
X \overset{h_1}{\longleftarrow} Y \overset{h_2}{\longleftarrow} Z
\]
that satisfy the prism condition $d_i   h_i + h_{i+1}   d_{i+1} = d \cdot \text{id}$ \ for all \ $i \in \{0, 1, 2 \}$. Let $\Phi_X: \overline{X_2} \longrightarrow X_1$ be the isomorphism defined as follows: For any $\overline{x} \in \overline{X_2}$, find a (non-unique) $x \in X_2$ such that $\pi(x) = \overline{x}$. Then multiply it by $d$. $d\cdot x$ is now unique, since any two such $x$'s differ by a multiple of $d$. Furthermore, $d \cdot x$ is a $d$-torsion element of $X_2$ and is therefore an element of $X_1$. Let $\mathcal{B}_X, \mathcal{B}_Y, \mathcal{B}_Z$ denote the Bockstein homomorphisms, which are the snake-lemma connecting maps: $\mathcal{B}_X: H^n(\overline{X_2}) \longrightarrow H^{n+1}(X_1)$, etc. Finally, let $\overline{\partial}: H^n(\overline{Z_2}) \longrightarrow H^{n+1}(\overline{X_2})$. Then
\begin{enumerate} 
\item There are well-defined homomorphisms $\widetilde{h_2}: Z_1 \longrightarrow X_1$ and $\widetilde{\overline{h_2}}: \overline{Z_2} \longrightarrow \overline{X_2}$ defined as $\widetilde{h_2} = d_1^{-1} h_2$, $\widetilde{\overline{h_2}} = d_1^{-1} h_2$, both of which are equal to $-h_1   f$, for any lifting $f$ of $d_2$. 
\item $\Phi_X^* \overline{\partial} = \widehat{h_2} \mathcal{B}_Z - \mathcal{B}_X \widehat{\overline{h_2}},$ where $\widehat{h_2}$ and $\widehat{\overline{h_2}}$ are the induced maps on cohomology. 
\end{enumerate} 
\end{lemma}
\noi {\bf Proof}: 
\begin{enumerate} 
\item Note that because $d_1$ is injective, $d_1^{-1}$ is a map from $\im(d_1)$ to either $X_1$ or $\overline{X_2}$. 
We will show that $\im(h_2) \subseteq \im(d_1)$ in $Z_1$ and $\overline{Z_2}$, which will prove that the definition of $\widetilde{h_2} = d_1^{-1}h_2$ makes sense. 
Let $f$ be a lifting of $d_2$, such a lifting exists because $d_2$ is surjective, and let $z$ be an element of either $Z_1$ or $\overline{Z_2}$. Then from the prism condition, 
\[ 
\begin{split} 
h_2(z) &= h_2(d_2(f((z)))) \\ 
&= (h_2   d_2)f(z) \\ 
&= d \cdot f(z) - (d_1   h_1)(f(z)) \\ 
&= -d_1(h_1(f(z))) \\
&= d_1(h_1(-f(z)))
\end{split} 
\] 
which is in the image of $d_1$. 
Furthermore, ker$(d_2) = \text{im} (d_1) \subseteq$ ker$(h_1)$ for both $Y_1$ and $\overline{Y_2}$, 
since both modules are $d$-torsion, making the composition $h_1   d_1$ the 0 map from the prism condition at $X$.  
Thus, any two liftings of $d_2$ differ by an element in the kernel of $h_1$, which shows that $-h_1f$ is well-defined. 

To show the equality $d_1^{-1}h_2 = -h_1 f$, we can apply $d_1^{-1}$ to both sides of the previous equality: 
\[ 
\begin{split} 
d_1^{-1}h_2 &= d_1^{-1} d_1 h_1 (-f) \\ 
&= h_1 (-f) \\ 
&= -h_1 f 
\end{split} 
\] 

\item We assume without loss of generality that $X \subseteq Y$ and $d_1$ is the inclusion map. 
This allows the the snake-lemma connecting maps to be computed by applying the coboundary map $\delta$ to any pre-image of a given cocycle. 
Computation will be done this way for both Bockstein maps $\mathcal{B}_X$, $\mathcal{B}_Z$ and for $\overline{\partial}$. 
The inclusions also make $h_1 = -h_2   d_2$ in $Y_1$ and $\overline{Y_2}$ (but not $Y_2$), and furthermore $h_1$ restricted to $X \subseteq Y$ is multiplication by $d$. 
Below is a diagram of the exact square. \\ 
\[
\begin{tikzpicture}[descr/.style={fill=white,inner sep=1.5pt}]
        \matrix (m) [
            matrix of math nodes,
            row sep=3em,
            column sep=3em,
            text height=1.5ex, text depth=0.25ex
        ]
        { X_1 & Y_1 & Z_1  \\ 
          X_2 & Y_2 & Z_2  \\ 
          \overline{X_2} & \overline{Y_2} & \overline{Z_2}  \\ 
        };

        \path[overlay, font=\scriptsize,>=latex]
        (m-1-1) edge [->]  node[descr,yshift=1.5ex] {$\subseteq$} (m-1-2) 
        (m-1-2) edge [->]  node[descr,yshift=1.5ex] {$d_2$} (m-1-3) 
        (m-1-1) edge [->]  node[descr,xshift=-1.5ex] {$\subseteq$} (m-2-1) 
        (m-1-2) edge [->]  node[descr,xshift=-1.5ex] {$\subseteq$} (m-2-2) 
        (m-1-3) edge [->]  node[descr,xshift=-1.5ex] {$\subseteq$} (m-2-3) 
        (m-2-1) edge [->]  node[descr,yshift=1.5ex] {$\subseteq$} (m-2-2) 
        (m-2-2) edge [->]  node[descr,yshift=1.5ex] {$d_2$} (m-2-3) 
        (m-2-1) edge [->]  node[descr,xshift=-1.5ex] {$\pi$} (m-3-1) 
        (m-2-2) edge [->]  node[descr,xshift=-1.5ex] {$\pi$} (m-3-2) 
        (m-2-3) edge [->]  node[descr,xshift=-1.5ex] {$\pi$} (m-3-3) 
        (m-3-1) edge [->]  node[descr,yshift=1.5ex] {$\subseteq$} (m-3-2) 
        (m-3-2) edge [->]  node[descr,yshift=1.5ex] {$d_2$} (m-3-3)

;
\end{tikzpicture} 
\] 
Let $\ell$ be any lifting of $d_2: \overline{Y_2} \longrightarrow \overline{Z_2}$, 
and let $\ell'$ be a lifting of $\pi: Y_2 \longrightarrow \overline{Y_2}$ that maps $\overline{X_2}$ into $X_2$. 
We define a lifting for $\pi: \overline{Z_2} \longrightarrow Z_2$ by $\ell_1:= d_2   \ell'   \ell$, 
and a lifting for $d_2: Y_2 \longrightarrow Z_2$ by $\ell_2:= \ell'   \ell   \pi$. 
With these choices for $\ell_1$ and $\ell_2$, the bottom right square commutes for the liftings, as 
\[
\ell_2 \ell_1 = (\ell' \ell \pi)(d_2 \ell' \ell) = \ell' \ell (\pi d_2) \ell' \ell = \ell' \ell (d_2 \pi \ell' \ell) = \ell' \ell.
\]
Let $\sigma_z \in Z^n(\overline{Z_2})$, and let $\sigma_y = \ell'(\ell(\sigma_z)) = \ell_2(\ell_1(\sigma_z)) \in C^n(Y_2)$. 
The liftings can be seen in the following diagram.  
\[
\begin{tikzpicture}[descr/.style={fill=white,inner sep=1.5pt}]
        \matrix (m) [
            matrix of math nodes,
            row sep=3em,
            column sep=3em,
            text height=1.5ex, text depth=0.25ex
        ]
        { X_1 & Y_1 & Z_1  \\ 
          X_2 & Y_2 & Z_2  \\ 
          \overline{X_2} & \overline{Y_2} & \overline{Z_2}  \\ 
        };

        \path[overlay, font=\scriptsize,>=latex]
        (m-2-3) edge [->,out=155,in=25]  node[descr,yshift=1.5ex] {$\ell_2$} (m-2-2) 
        (m-1-3) edge [->,out=155,in=25]  node[descr,yshift=1.5ex] {$\ell_2$} (m-1-2) 
        (m-3-3) edge [->,out=155,in=25]  node[descr,yshift=1.5ex] {$\ell$} (m-3-2) 
        (m-3-3) edge [->,out=115,in=245]  node[descr,xshift=-1.5ex] {$\ell_1$} (m-2-3) 
        (m-3-2) edge [->,out=115,in=245]  node[descr,xshift=-1.5ex] {$\ell'$} (m-2-2) 
        (m-3-1) edge [->,out=115,in=245]  node[descr,xshift=-1.5ex] {$\ell'$} (m-2-1) 
        (m-1-1) edge [->, dotted]   (m-1-2) 
        (m-1-2) edge [->, dotted]   (m-1-3) 
        (m-1-1) edge [->, dotted]   (m-2-1) 
        (m-1-2) edge [->, dotted]   (m-2-2) 
        (m-1-3) edge [->, dotted]   (m-2-3) 
        (m-2-1) edge [->, dotted]   (m-2-2) 
        (m-2-2) edge [->, dotted]   (m-2-3) 
        (m-2-1) edge [->, dotted]   (m-3-1) 
        (m-2-2) edge [->, dotted]   (m-3-2) 
        (m-2-3) edge [->, dotted]   (m-3-3) 
        (m-3-1) edge [->, dotted]   (m-3-2) 
        (m-3-2) edge [->, dotted]   (m-3-3)    
;
\end{tikzpicture} 
\] 
We use $[\sigma]$ to denote the cohomology class of $\sigma$ whenever $\sigma$ is a cocycle. We recall that $d_1^{-1}$ is the identity map. Direct computation and an application of the prism condition yield the desired result as follows. \\ 
\[
    \begin{split} 
    -\mathcal{B}_X (\widehat{\overline{h}}_2([\sigma_z])) &= -\mathcal{B}_X ([\overline{h}_2\sigma_z)]) = \mathcal{B}_X ([-\overline{h}_2\sigma_z]) = \mathcal{B}_X ([h_1 \ell \sigma_z]) \\ 
    &= \mathcal{B}_X ([h_1 \pi \sigma_y]) = \mathcal{B}_X ([\pi h_1 \sigma_y]) = [\delta \ell' \pi h_1 \sigma_y] \\ 
    &= [\delta h_1 \sigma_y] = [h_1 (\delta \sigma_y)] \\ \\ \\ 
\end{split}
\] 
\[ 
\begin{split}
    \widehat{h_2}(\mathcal{B}_Z([\sigma_z])) &= \widehat{h_2}(\mathcal{B}_Z( [\pi d_2 \sigma_y])) = \widehat{h_2}([(\delta \ell') \pi d_2\sigma_y]) \\ 
    &= \widehat{h_2}([\delta d_2 \sigma_y]) = [h_2 \delta d_2 \sigma_y] = [(h_2   d_2)(\delta \sigma_y)] \\ \\ \\ 
    \Phi_X^*(\overline{\partial}([\sigma_z])) &= \Phi_X^*([\delta(\ell \sigma_z)]) = \Phi_X^*([\delta \pi \sigma_y]) = \Phi_X^*([\pi \delta \sigma_y]) \\ 
    &= [(\Phi_X  \pi) \delta \sigma_y] = [d \cdot \delta \sigma_y] \\ \\ 
    \end{split} 
\]
Using these computations and the prism condition at $Y_1$, we have the following equalities for any $\sigma_z \in Z^n(\overline{Z_2})$. 
\[
\Phi_X^*\overline{\partial}([\sigma_z]) = [d \cdot (\delta \sigma_y)] = [(h_2   d_2)(\delta\sigma_y) - h_1(\delta \sigma_y)] = \widehat{h_2}(\mathcal{B}_Z([\sigma_z])) - \mathcal{B}_X (\widehat{h_2}([\sigma_z]))
\]
and therefore 
\[
\Phi_X^* \overline{\partial} = \widehat{h_2}  \mathcal{B}_Z - \mathcal{B}_X \widehat{h_2} .
\]
\end{enumerate} 

This concludes the proof of Lemma 
\ref{Positselski Prism}. \hfill $_\square$ \vspace{10pt}

This next lemma gives the definition of the
map $\eta$ and uses a splitting of the four
term sequence into two three term exact sequences. 
It shows how $\eta$ is related to the connecting maps
of these short exact sequences via the homotopies provided by the Positeselski Hypotheses.

\begin{lemma} \label{Positselski Eta} 
Using the language of the previous lemma and viewing the 4-term exact sequence 

\[
0 \longrightarrow A  \longrightarrow B \longrightarrow C \longrightarrow D \longrightarrow 0
\]
as two short exact sequences 
\[
0 \longrightarrow A \longrightarrow B \longrightarrow \Delta \longrightarrow 0
\]
\[
0 \longrightarrow \Delta \longrightarrow C \longrightarrow D \longrightarrow 0
\]
with the second sequence starting with an inclusion $\Delta \subseteq C$, and given two liftings $f_2: \Delta \longrightarrow B$ and $f_3: D \longrightarrow C$ 
for the the terminal, surjective maps in each of these two respective short exact sequences, we have the following equivalent definitions of $\eta$: 
\[
-\widehat{h_2} \partial_{\Delta_1 D_1} = \eta = \partial_{A_1 \Delta_1} \widehat{h_3}
\]
where $\widetilde{h_2} := d_1^{-1}h_2 = -h_1 f_2 :C \longrightarrow A$ for the first short exact sequence as in the previous lemma, and 
$\widetilde{h_3} := h_3 = -h_2 f_3: D \longrightarrow \Delta$ is the analogue of $\widetilde{h_2}$ for the second short exact sequence. 
Also, $f_3$ is any lifting of $d_3$ (not necessarily a $\cG$-map) and $\delta: C^n(\G, C_1) \longrightarrow C^{n+1}(\G, C_1)$ is the group cohomology coboundary map. 

\end{lemma} 
Before proceeding with the proof, we note that though it is tempting to define $\widetilde{h_3}:D \longrightarrow \Delta$ as $f_2 h_3$, where $f_2$ is a lifting from $C$ to $B$ in the four-term exact sequence, however $\Delta$ is already a submodule of $C$, so this is not necessary. The lifting that proceeds $h_3$ needs to be a lifting of the inclusion, which is the identity map in this case. 

\noi {\bf Proof}: 
The first homotopy we need for the short exact sequence 
\[
0 \longrightarrow \Delta \overset{\subseteq}{\longrightarrow} C \overset{d_3}{\longrightarrow} D \longrightarrow 0
\]
is $d_2 h_2: C \longrightarrow \Delta$, and the fact that it satisfies the prism condition at $\Delta$ and $C$ follows from the prism condition being satisfied for the four-term exact sequence. We are ready to prove equality claimed between the two definitions of $\eta$. Let $\sigma \in Z^1(\cG, D_1)$. Then
\[ 
\begin{split} 
-\widehat{h_2} \partial_{\Delta_1 D_1}([\sigma]) &= [-(h_1 f_2)(\delta f_3 \sigma)] = [(-h_1 f_2)(\delta f_3 \sigma)] = [\widetilde{h_2}(\delta f_3\sigma)] = [(d_1^{-1} h_2)(\delta f_3\sigma)] \\ 
&= [d_1^{-1} h_2 \delta f_3 \sigma] = [d_1^{-1} \delta h_2 f_3 \sigma] = [d_1^{-1} \delta (f_2 d_2) h_2 f_3 \sigma] \\ 
&= [(d_1^{-1} \delta d_2^{-1})( d_2 h_2 f_3) \sigma] = \partial_{A_1 \Delta_1} \widehat{h_3} ([\sigma])
\end{split} 
\] 

The third equality above uses the prism condition to interchange the two definitions of $\widetilde{h_2}$ as in the previous lemma, 
the sixth equality follows from the commutativity of the chain map with the homotopies, 
the seventh equality is true because With the right choice of lifting (not necessarily $f_2$), applying $d_2$ and then the lifting is the same as the identity map,  
but this difference in lifting choice is just a coboundary, which means that any lifting of $d_2$ yields the desired equality. 
Finally, the last equality follows from the definitions of the connecting map $\partial_{A_1\Delta_1}$ and the homomorphism $\widetilde{h_3}$ from $D$ to $\Delta$. 

This concludes the proof of Lemma \ref{Positselski Eta}. \hfill $_\square$ \vspace{10pt}

The reader may note that the only difference between this definition of 
$\eta$ and the definition of $\partial$ from the snake lemma is the $h_2$ in the composition. 
$h_2$ can be viewed as the way to connects the two middle terms of the 4-term exact sequence, 
where the snake lemma has 1 middle term in a 3-term exact sequence and there is no need for this intermediate step. 
In this sense they are as close as can be considering the different number of modules in the exact sequence. 

With Lemmas \ref{Positselski Prism} and \ref{Positselski Eta} proved, we move on to the proof of Positselski's theorem. Up to this point we have not
used the assumption from the Positselski Hypotheses that the Bochstein maps are zero.
We will do so now. However as noted in Remark 2.5
below we will see that we
need not assume that the Bochstein map $\mathcal{B}_D$  is zero; we shall only use that the composite $\widehat{h_3}\circ \mathcal{B}_D$ is zero. \\ \\
\noi {\bf Proof of Theorem \ref{Positselski's Theorem}}:  
With $\eta$ defined, it is required to show exactness at the inner four terms of the 6-term sequence. 

\smallskip\noi
\underline{Exactness of $H^n(B_1 \oplus D_1) \overset{d_2^*+h_3^*}{\longrightarrow} H^n(C_1) \overset{d_3^*}{\longrightarrow} H^n(D_1)$}: The composition $d_3 d_2$ is the 0 map on the modules, which makes it the zero map on cohomology as well. 
The composition $d_3 h_3$ is also the 0 map because it is multiplication by $d$ and $D_1$ is $d$-torsion.   

Now we show that $d_2^*+h_3^*$ maps onto the kernel of $d_3^*$. 
Let $\sigma_c \in Z^n(\cG,C_1)$ such that $[\sigma_c]$ is in the kernel of $d_3^*$. Let $\Delta_1 = \ker(d_3) \subseteq C_1$.
Then, after adding a coboundary from $B^n(\cG, C_1)$ to $\sigma_c$ if necessary, we may assume that $\sigma_c \in Z^n(\cG, \Delta_1) \subseteq Z^n(\cG, C_1)$.
Let $\overline{\sigma}_c \in Z^n(\cG, \overline{\Delta}_2)$ be the corresponding element to $\sigma_c$ from the isomorphism $\Delta_1 \cong \overline{\Delta}_2$. 
Then $\mathcal{B}_C([\overline{\sigma}_c])=0$  in $H^{n+1}(\cG, C_1)$ by the assumption $\mathcal{B}_C=0$. 
Therefore $\mathcal{B}_\Delta([\overline{\sigma}_c])$ is in the kernel of $\iota^*$, where $\iota: \Delta \longrightarrow C$ is the inclusion. 
So $\mathcal{B}_\Delta([\overline{\sigma}_c]) = \partial_{\Delta D}([\sigma_d])$ for some $\sigma_d \in Z^n(\cG, D_1)$ 
from the exactness of the long exact sequence induced by the Snake Lemma from the short exact sequence 
\[
0 \longrightarrow \Delta_1 \longrightarrow C_1 \longrightarrow D_1 \longrightarrow 0.
\]
With $\mathcal{B}_\Delta([\overline{\sigma}_c])$ being in the image of $\partial_{\Delta_1 D_1}$, 
we can apply Lemma \ref{Positselski Prism} to the short exact sequence 
\[
0 \longrightarrow A \longrightarrow B \longrightarrow \Delta \longrightarrow 0
\]
with the following result: Let the $\widetilde{h_2}: \Delta_1 \longrightarrow A_1$ be as in the lemma. Then the following equalities hold in $H^{n+1}(A_1)$. 
\[
\partial_{A_1 \Delta_1}([\sigma_c]) = \mathcal{B}_A \widehat{h_2}([\overline{\sigma}_c]) - \widehat{h_2}\mathcal{B}_{\Delta}([\overline{\sigma}_c]) = - \widehat{h_2}\mathcal{B}_{\Delta}([\overline{\sigma}_c]) = - \widehat{h_2}\partial_{\Delta_1 D_1}([\sigma_d])  = \eta([\sigma_d])
\]
with the second equality holding because $\mathcal{B}_A=0$ and the last equality following from the definition of $\eta$. Furthermore, the equivalent definition of $\eta$ as $\partial_{A_1 \Delta_1} \widehat{h_3}$ yields the equality 
\[
\partial_{A_1 \Delta_1}([\sigma_c]) = \partial_{A_1 \Delta_1} \widehat{h_3} ([\sigma_d]),
\]
where $\widehat{h_3}$, as in Lemma \ref{Positselski Eta}, comes from the short exact sequence 
\[
0 \longrightarrow \Delta_1 \longrightarrow C_1 \longrightarrow D_1 \longrightarrow 0.
\]
The desired result follows from this: Since $\partial_{A_1 \Delta_1} \left( [\sigma_c] - \widehat{h_3}([\sigma_d])\right)=0 \in H^{n+1}(\cG, A_1)$, it follows that $[\sigma_c] - \widehat{h_3}([\sigma_d]) = d_2^*([\sigma_b])$ for $\sigma_b \in Z^n(B_1)$. 
Therefore $[\sigma_c] = \widehat{h_3}([\sigma_d]) + d_2^*([\sigma_b])$ for some $\sigma_D \in Z^n(\cG, D_1), \sigma_b \in H^n(\cG, B_1)$. 
This concludes exactness at $H^n(\cG, C_1)$.

\smallskip\noi
\underline{Exactness of $H^n(C_1) \overset{d_3^*}{\longrightarrow} H^n(D_1) \overset{\eta}{\longrightarrow} H^{n+1}(A_1)$}: 
We first show that $\eta d_3^* = 0$ by using the definition of $\eta = -\widehat{h}_2 \partial_{\Delta_1 D_1}$ in Lemma \ref{Positselski Eta} and the fact that $\partial_{\Delta_1 D_1} d_3^* = 0$ from exactness of the long exact sequence from the short exact sequence 
\[
0 \longrightarrow \Delta_1 \longrightarrow C_1 \longrightarrow D_1 \longrightarrow 0.
\]
In fact, this means the kernel of $\partial_{\Delta_1 D_1}$ is equal to the image of $d_3^*$, so the other containment amounts to showing that the kernel of $\partial_{\Delta_1 D_1}$ is no smaller than the kernel of $-\widehat{h}_2 \partial_{\Delta_1 D_1}$. We shall prove this fact. 
However, we will start with the equivalent definition of $\eta$ in Lemma \ref{Positselski Eta}, and then show this equality. 

Suppose $\sigma_d \in Z^n(\cG, D_1)$ such that $[\sigma_d]$ is in the kernel of $\eta$. 
From Lemma \ref{Positselski Eta}, $\eta([\sigma_d]) = \partial_{A_1 \Delta_1} \widehat{h_3}([\sigma_d])$, which means $\widehat{h}_3([\sigma_d])$ is in the kernel of $\partial_{A_1 \Delta_1}$, the connecting map for the short exact sequence 
\[ 
0 \longrightarrow A_1 \longrightarrow B_1 \longrightarrow \Delta_1 \longrightarrow 0, 
\]
and therefore $\widehat{h}_3([\sigma_d]) = d_2^*([\sigma_b])$ for some $\sigma_b \in Z^n(\cG, B_1)$. 
Let $\overline{\sigma}_b, \overline{\sigma}_d$ be the corresponding cocycles to $\sigma_b, \sigma_d$ from the isomorphisms $B_1 \cong \overline{B}_2, D_1 \cong \overline{D}_2$, respectively. 
Then from the commuting of the Bockstein maps with the cohomology maps $\widehat{h_3}$ induced from module homomorphisms $\widetilde{h_3}$ and the vanishing of $\mathcal{B}_B$, we have the following equalities. 
\[
0 = d_2^* (\mathcal{B}_B([\overline{\sigma}_b])) = \mathcal{B}_\Delta (d_2^*([\overline{\sigma}_b])) = \mathcal{B}_\Delta(\widehat{h_3}([\sigma_d])) \] 
From Lemma \ref{Positselski Prism} and the vanishing of the 
composite  $\widehat{h_3} \mathcal{B}_D$, we have
\[\mathcal{B}_\Delta \widehat{h_3}([\overline{\sigma}_d]) = \mathcal{B}_\Delta \widehat{h_3}([\overline{\sigma}_d]) - \widehat{h_3} \mathcal{B}_D([\overline{\sigma}_d]) = \Phi^* \partial_{\Delta_1 D_1}([\overline{\sigma}_d]) = \partial_{\Delta_1 D_1}([\sigma_d]).
\]
Therefore $\partial_{\Delta_1 D_1}([\sigma_d]) = 0$. This concludes exactness at $H^n(\cG, D_1)$. 

\smallskip\noi
\underline{Exactness of $H^n(D_1) \overset{\eta}{\longrightarrow} H^{n+1}(A_1)\overset{d_1}{\longrightarrow} H^{n+1}(B_1)$}: To show that $d_1 \eta = 0$, let $\sigma_d \in Z^n(D_1)$. Then by Lemma \ref{Positselski Eta}, 
\[
d_1^* \eta ([\sigma_d]) = d_1^* (-\partial_{A_1 \Delta_1} \widehat{h_3})([\sigma_d]) = -(d_1^* \partial_{A_1 \Delta_1})( \widehat{h_3}([\sigma_d])),
\]
with $d_1^* \partial_{A_1 \Delta_1}$ being zero because of exactness of the long exact sequence from the short exact sequence 
\[
0\longrightarrow A_1 \longrightarrow B_1 \longrightarrow \Delta_1 \longrightarrow 0.
\]
Now suppose $\sigma_a \in Z^{n+1}(\cG, A_1)$ and $[\sigma_a]$ is in the kernel of $d_1^*$. Then from exactness of the long exact sequence from the short exact sequence 
\[0 \longrightarrow A_1 \longrightarrow B_1 \longrightarrow \Delta_1 \longrightarrow 0,
\]
$[\sigma_a] = \partial_{A_1 \Delta_1} ([\sigma_x])$ for some $\sigma_x \in Z^n(\cG, \Delta_1)$. Let $\overline{\sigma}_x \in Z^n(\cG, \overline{\Delta}_2)$ be the element corresponding to $\sigma_x$ through $\Phi^*$. 
Then 
\[
[\sigma_a] = \partial_{A_1\Delta_1}\Phi_\Delta^*([\overline{\sigma}_x]) = \Phi_A^*\partial_{A_1\Delta_1}([\overline{\sigma}_x])
\]
From Lemma \ref{Positselski Prism} and the vanishing of the Bockstein map $\mathcal{B}_A$, 
\[
\Phi_A^*\partial_{A_1\Delta_1}([\overline{\sigma}_x]) = \widehat{h}_2 \mathcal{B}_\Delta([\overline{\sigma}_x])  - \mathcal{B}_A \widehat{\overline{h}}_2([\overline{\sigma}_x]) = \widehat{h}_2 \mathcal{B}_\Delta([\overline{\sigma}_x]).
\]
Furthermore, from the commuting of the Bockstein maps with the $d_i$'s and the vanishing of $\mathcal{B}_C$, we have
\[
d_3^* \mathcal{B}_\Delta([\overline{\sigma}_x]) = \mathcal{B}_C i^*([\overline{\sigma}_x]) = 0,
\]
where $i^*$ is induced by the inclusion $\subseteq: \Delta \longrightarrow C$. By exactness of the long exact sequence from the short exact sequence
\[
0 \longrightarrow \Delta_1 \longrightarrow C_1 \longrightarrow D_1 \longrightarrow 0, 
\]
$\mathcal{B}_\Delta([\overline{\sigma}_x]) = \partial_{\Delta_1 D_1}([\sigma_d])$ for some $\sigma_d \in Z^n(\cG, D_1)$. Therefore, using a definition of $\eta$ in Lemma \ref{Positselski Eta}, 
\[
[\sigma_a] = \widehat{h}_2 \mathcal{B}_\Delta([\overline{\sigma}_x]) = \widehat{h}_2 \partial_{\Delta_1 D_1}([\sigma_d]) = \eta([-\sigma_d]).
\]
This concludes exactness at $H^n(\cG, A_1)$.

\smallskip\noi
 \underline{Exactness of $H^{n+1}(A_1) \overset{d_1}{\longrightarrow} H^{n+1}(B_1) \overset{d_2 \oplus h_1}{\longrightarrow} H^n(A_1 \oplus C_1)$}:
Both of the compositions $d_2 d_1$ and $h_1 d_1$ are 0 maps. In the first case this comes from exactness of the $d_i$'s, and in the second case from $A_1$ being $d$-torsion. 

For the other containment, let $\widehat{d}_2: H^{n+1}(\cG, B_1) \longrightarrow H^{n+1}(\cG, \Delta_1)$ be the map induced by $d_2$ with codomain reduced to $\Delta_1$, so that $\iota^* \widehat{d}_2 = d_2^* : H^{n+1}(\cG, B_1) \longrightarrow H^{n+1}(\cG, C_1)$, where $\iota$ is induced from the inclusion $\iota: \Delta_1 \overset{\subseteq}{\longrightarrow} C_1$. 
Let $\sigma_b \in Z^{n+1}(\cG, B_1)$ such that $h_1^*([\sigma_b]) = 0$ and $d_2^*([\sigma_b]) = 0$. 
We will show that $[\sigma_b] = d_1^*([\sigma_a])$ for some $\sigma_a \in Z^{n+1}(\cG, A_1)$ by showing that $\widetilde{d}_2([\sigma_b]) = 0$. 

From the hypothesis, $0 = d_2^*([\sigma_b]) = \iota^* (\widehat{d}_2([\sigma_b]))$, which means 
$\widehat{d}_2([\sigma_b]) = \partial_{\Delta_1 D_1}([\sigma_d])$ for some $\sigma_d \in Z^n(\cG, D_1)$. 
This follows from the exactness of the long exact sequence from the short exact sequence
\[ 
0 \longrightarrow \Delta_1 \longrightarrow C_1 \longrightarrow D_1 \longrightarrow 0. 
\]
Using the definitions of $\eta$ from Lemma \ref{Positselski Prism}, we have 
\[
0 = h_1^*([\sigma_b]) = \widehat{h}_2 \widehat{d}_2 ([\sigma_b]) = \widehat{h}_2 \partial_{\Delta_1 D_1}([\sigma_b]) = - \eta ([\sigma_b]) = -\partial_{A_1 \Delta_1} \widehat{h}_3 ([\sigma_b]), 
\]
which means that $\widehat{h}_3 ([\sigma_b]) = \widehat{d}_2 ([\sigma_b'])$ for some $\sigma_b \in Z^{n+1}(\cG, B_1)$ from exactness of the long exact sequence induced by the short exact sequence 
\[
0 \longrightarrow A_1 \longrightarrow B_1 \longrightarrow \Delta_1 \longrightarrow 0.
\]

Let $\overline{\sigma}_d$ and $\overline{\sigma}_b'$ be the corresponding element to $\sigma_d$ and $\sigma_b'$ through the isomorphism $\Phi^*$ (respectively). Then from Lemma \ref{Positselski Prism} and the assumption  $\widehat{h_3} \mathcal{B}_D=0$,
\[
\begin{split}
\widetilde{d}_2([\sigma_b]) = \partial_{\Delta_1 D_1}([\sigma_d]) &= \widehat{h_3}\mathcal{B}_D([\overline{\sigma}_d]) - \mathcal{B}_\Delta \widehat{h}_3([\overline{\sigma}_d]) = - \mathcal{B}_\Delta \widehat{h}_3([\overline{\sigma}_d]) \\
&= - \mathcal{B}_\Delta \widehat{d}_2([\overline{\sigma}_b']) = - \widehat{d}_2 \mathcal{B}_B([\overline{\sigma}_b']) = 0.
\end{split}
\]
This yields the desired result of $\widehat{d}_2([\sigma_b]) = 0$. Therefore, from the exactness of the long exact sequence induced by the short exact sequence
\[
0 \longrightarrow A_1 \longrightarrow B_1 \longrightarrow \Delta_1 \longrightarrow 0
\]
we have $[\sigma_b] = d_1^*([\sigma_a])$ for some $\sigma_a \in Z^{n+1}(\cG, A_1)$. \hfill $_\square$ \vspace{10pt}

\begin{remark} The proof given above follows the proof given by Positselski 
\cite{Positselski}. In Positeselski's write-up all four of the Bockstein maps were assumed to be zero for the six-term sequence to be exact, 
however as noted above this hypothesis can be weakened to assuming that $\mathcal{B}_A$, $\mathcal{B}_B$ and $\mathcal{B}_C$  are zero along with the vanishing of the composite $\widehat{h_3} \mathcal{B}_D$. Tabulating the use of these hypotheses, the exactness at $H^n(C_1)$ only required $\mathcal{B}_A=0$ and $\mathcal{B}_C=0$, at $H^n(D_1)$ exactness requires 
$\widehat{h_3} \mathcal{B}_D=0$, at $H^{n+1}(A_1)$ exactness requires $\mathcal{B}_A=0$ and $\mathcal{B}_C=0$, and exactness at $H^{n+1}(B_1)$ requires that $\widehat{h_3} \mathcal{B}_D=0$ and $\mathcal{B}_B=0$. Exactness at $H^{n+1}(A_1)$ is of particular interest when $n = 1$ because it can be used to calculate cohomological (in particular Brauer) kernels.  These  observations will be important in a sequel to this paper where four term exact sequences with homotopies are constructed, but 
where the fourth Bochstein map $\mathcal{B}_D$ fails to be zero, while the composite
$\widehat{h_3} \mathcal{B}_D$ is shown to be zero.
\end{remark}

\section{The General Setup and the Bockstein Maps} \label{General Setup Section}

For the sections that follow we adopt the following notation:  

Let $G$ be a semi-direct product of $\langle\tau\rangle$ by $\langle\sigma\rangle$ with $|\tau| = d$, 
$|\sigma| = s$ and $\langle \sigma \rangle$ has a faithful action on $\langle \tau \rangle$. 
Let $\cG = \Gal(F^{\text{sep}}/F)$ for some field $F$ whose extension $E$ is the degree-$d$ extension we wish to study. 
Let $\widetilde{E}$ be the Galois closure of $E/F$, $\cN = \Gal(F^{\text{sep}}/\widetilde{E}) \triangleleft \cG$. 
Assume further that $G \cong {\rm Gal}(\tilde{E}/F)\cong \cG/\cN$, $E = F(\beta)$, $\widetilde{E} = F(\alpha, \beta)$, 
and let $\widetilde{F}$ be the cyclic  extension of $F$ given by Fix$(\langle \tau \rangle)$ in this setup. 
Let $\cH = \Gal(F^{\text{sep}}/E)$, $\cJ = \Gal(F^{\text{sep}}/\widetilde{F})$, 
$H = \langle \sigma \rangle$ and $J = \langle \tau \rangle$ so that 
$\sigma(\beta) = \beta$ and $\tau(\alpha) = \alpha$. 
We also assume that the $d^2$th roots of unity $\mu_{d^2} \subseteq F$, 
and char$(F)$ does not divide $d$ so that $\mu_{d^2}$ contains $d^2$ distinct roots of unity.  

The following three sections will cover the cases of $E/F$ in increasing generality. 
In the next section $E/F$ will be a cyclic extension so that $\widetilde{E} = E$, 
$\widetilde{F} = F$ and $s = |\sigma| = 1$ with no restrictions on $d$. 
In the section that follows we assume $s = |\sigma| = 2$ so that $G$ 
is a dihedral group with $d$ odd. 
Lastly we let $s = |\sigma|$ be any even positive integer with $s | (d-1)$, hence $d$ is still odd.

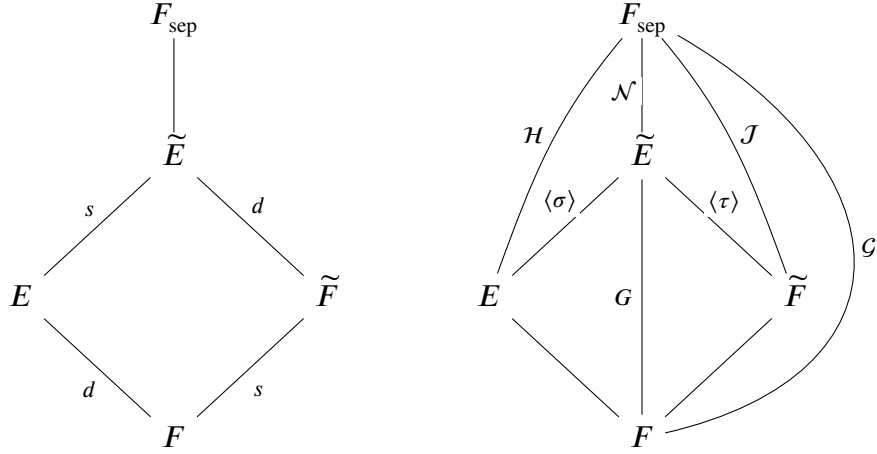
\[
\begin{tikzpicture}[descr/.style={fill=white,inner sep=1.5pt}]
        \matrix (m) [
            matrix of math nodes,
            row sep=3em,
            column sep=3em,
            text height=1.5ex, text depth=0.25ex
        ]
        {\ & F_{\text{sep}} & \ \\ 
         \ & \widetilde{E} & \                                        \\ 
            E & \ & \widetilde{F} \\ 
           \ & F & \                                        \\ \\ 
        };

        \path[overlay, font=\scriptsize,>=latex]
        (m-2-2) edge  node[descr,xshift=-1.75ex] {$\ $} (m-1-2) 
        (m-3-1) edge  node[descr,xshift=-3,yshift=1.1ex] {$s $}(m-2-2) 
        (m-3-3) edge node[descr,xshift=3,yshift=1.5ex] {$d$ }(m-2-2) 
        (m-4-2) edge node[descr,xshift=-3,yshift=-1.5ex] {$d$ } (m-3-1) 
        (m-4-2) edge node[descr,xshift=3,yshift=-1.6ex] {$s$ } (m-3-3) 
;
\end{tikzpicture} \hspace{.5in}
\begin{tikzpicture}[descr/.style={fill=white,inner sep=1.5pt}]
        \matrix (m) [
            matrix of math nodes,
            row sep=3em,
            column sep=3em,
            text height=1.5ex, text depth=0.25ex
        ]
        {\ & F_{\text{sep}} & \ \\ 
         \ & \widetilde{E} & \                                        \\ 
            E & \ & \widetilde{F} \\ 
           \ & F & \                                        \\ \\ 
        };

        \path[overlay, font=\scriptsize,>=latex]
        (m-2-2) edge  node[descr,xshift=-1.25ex,yshift=-2] {$\cN$} (m-1-2) 
        (m-3-1) edge[out=70,in=230] node[descr,xshift=-1ex,yshift=1ex] {$\cH$} (m-1-2) 
        (m-3-3) edge[out=110,in=310] node[descr,xshift=1ex,yshift=1ex] {$\cJ$} (m-1-2) 
        (m-3-1) edge  node[descr,xshift=-2.1,yshift=1.75ex] {$\langle \sigma \rangle$} (m-2-2) 
        (m-3-3) edge  node[descr,xshift=2.1,yshift=1.75ex] {$\langle \tau \rangle$} (m-2-2) 
        (m-4-2) edge[out=13,in=-30,looseness = 1.7] node[descr,xshift=1ex,yshift=1ex] {$\cG$} (m-1-2) 
        (m-4-2) edge node[descr,xshift=-7,yshift=0ex] {$G$} (m-2-2) 
        (m-4-2) edge (m-3-1) 
        (m-4-2) edge (m-3-3) 
;
\end{tikzpicture} 
\] 

We denote by $\theta:\{0,1,\ldots,d-1\}\rightarrow \{0,1,\ldots,d-1\}$ the conjugation in $\langle\tau\rangle$ by $\sigma$, that is, $\sigma\tau^i\sigma^{-1}=\tau^{\theta(i)}$.  
We define $\theta_j$ by $\sigma^j\tau\sigma^{-j}=\tau^{\theta_j}$ (in fact $\theta_j=\theta^j(1)$ where the 
latter is the $j$'th iterate of $\theta$, but the notation $\theta_j$ is less cumbersome.) 
We assume that $\theta$ has order $s$, that is, conjugation by $\sigma$ on 
$\langle\tau\rangle$ has order $s$. 
We assumed this from the faithful action in the semi-direct product setup of $G$ above at the beginning of this section. 
As $\tau$ has odd order,  $\sigma^{\frac{s}{2}}\tau^i\sigma^{-\frac{s}{2}}= \tau^{-i}$ for all $i$. 
From this, $\theta_{j+\frac{s}{2}}\equiv-\theta_j$ (mod $d$) and since $0<\theta_j<d$ 
we must have $\theta_j+\theta_{j+\frac{s}{2}}=d$.  

\bigskip

Part of the Positselski Hypotheses is the requirement that the Bockstein maps are zero for the four modules in the exact sequence, 
and we will show that this is indeed the case for the next three sections. 
To facilitate this in the later sections, 
we will prove a couple lemmas here. To set up these lemmas, let us first examine the long exact sequence over which the Bockstein map is defined for the module 
$\mu_{d^2}$ for a given field $F$ that contains those roots of unity. 
We will be using $M_1 = \mu_{d^2}$ as a trivial $\cG$-module in all three of the following sections, 
and $M_2$ will be an induced module with the same cohomology when taken over a slightly larger field.  

The long exact sequence associated with 
$0\rightarrow \mu_d\rightarrow \mu_{d^2}\rightarrow \mu_{d^2}/\mu_d\rightarrow0$
is the following 
\[
\cdots \longrightarrow H^n(F, \mu_d) \overset{i}{\longrightarrow} H^n(F, \mu_{d^2}) \overset{\pi}{\longrightarrow} H^n(F, \mu_{d^2}/\mu_d) \overset{\beta_{\mu_{d^2}}}{\longrightarrow} H^{n+1}(F, \mu_d) \longrightarrow \cdots 
\]
and the relevant Bockstein map is the connecting map
labelled $\beta_{\mu_{d^2}}$ in this sequence (see also \cite{SV}). In this case the
vanishing of the Bockstein map is given next.

\begin{lemma}\label{FirstBocksteinVanishing}
Suppose $F$ is a field with $\mu_{d^2}\subset F$. Then the Bockstein map $\beta_{\mu_{d^2}}$ associated with
the short exact sequence $0\rightarrow \mu_d\rightarrow \mu_{d^2}\rightarrow \mu_{d^2}/\mu_d\rightarrow0$
is zero.
\end{lemma}

\noi {\bf Proof.}
The Bloch-Kato Conjecture, proved by Veovodski in \cite{V} (see section 1.7 of \cite{HW} for discussion of history), states that the norm residue homomorphisms below are surjective since $(\text{char}(F),d)=1$. Furthermore, the identity map commutes with the canonical quotient map through the norm residue homomorphism. This means
we have a commutative diagram,
\[ 
\begin{tikzpicture}[descr/.style={fill=white,inner sep=1.5pt}]
        \matrix (m) [
            matrix of math nodes,
            row sep=3em,
            column sep=3em,
            text height=1.5ex, text depth=0.25ex
        ]
        { \ & K^M_n & K^M_n & \ \\ 
        H^n(F, \mu_d) & H^n(F, \mu_{d^2}) & H^n(F, \mu_{d^2}/\mu_d) & H^{n+1}(F, \mu_d) \\ 
            };

        \path[overlay, font=\scriptsize,>=latex]
        (m-1-2) edge [->] node[descr,yshift=1.5ex] {$ \text{id} $}(m-1-3)
        (m-2-1) edge [->] node[descr,yshift=1.5ex] {$\subseteq$} (m-2-2)
        (m-2-2) edge [->] node[descr,yshift=1.5ex] {$ \pi $} (m-2-3)
        (m-2-3) edge [->] node[descr,yshift=1.5ex] {$\beta_{\mu_{d^2}}$} (m-2-4)
        (m-1-2) edge [->>]  node[descr,xshift=1.7ex] {$s_{n,d^2}$} (m-2-2) 
        (m-1-3) edge [->>]  node[descr,xshift=1.5ex] {$s_{n,d}$} (m-2-3) 
;
\end{tikzpicture} 
\] 
Therefore the map $\pi$ is surjective as well, so that $\beta_{\mu_{d^2}} = 0$ by exactness. This proves the lemma.
\hfill $_\square$ \vspace{10pt}

Lemma \ref{FirstBocksteinVanishing}
checks the vanishing of the Bockstein map for $\mu_{d^2}$, which takes care of $M_1$ in every case
we will consider. Similarly, $M_2$ is either an induced module or a sum of such in every case, and Lemma 7 will apply. 
The field $E$ may be used instead of $F$ to apply Lemma \ref{FirstBocksteinVanishing} and obtain the result that $\beta_{M_2} = 0$. 
For the cyclic case, $M_3=M_2$ and $M_4=M_1$, so all the Bockstein maps have been shown to be zero for the cyclic case. The dihedral and semi-direct cases have different $M_3$'s and $M_4$'s, and these cases will be treated separately in the development of the verification of their respective Positselski Hypotheses. In doing so, the following lemma will be used. 

\begin{lemma} \label{Zero Bockstein} 
Let $\cJ \subseteq \cG$ be an index $s$ subgroup, $X_2$ a free $\Z/d^2\Z$-module with a discrete action of $\G$ with $(s,d) = 1$. If $\beta_{\cJ,X} = 0$ then $\beta_{\cG,X} = 0$ as well. 
\end{lemma} 
\noi {\bf Proof}: The restriction map $res:H^n(\cG,X_2) \longrightarrow H^n(\cJ,X_2)$ is injective because $cor\circ res$ = $\cdot s$, which is invertible because $(s,d) = 1$. Now, we have the following commutative diagram 
\[ 
\begin{tikzpicture}[descr/.style={fill=white,inner sep=1.5pt}]
        \matrix (m) [
            matrix of math nodes,
            row sep=3em,
            column sep=3em,
            text height=1.5ex, text depth=0.25ex
        ]
        {H^n(\cG, X_1) & H^n(\cG, X_2) & H^n(\cG, \overline{X_2})) & H^{n+1}(\cG, X_1) \\ 
        H^n(\cJ, X_1) & H^n(\cJ, X_2) & H^n(\cJ, \overline{X_2})) & H^{n+1}(\cJ, X_1) \\ 
            };

        \path[overlay, font=\scriptsize,>=latex]
        (m-1-1) edge [->] node[descr,yshift=1.5ex] {$i$} (m-1-2)
        (m-1-2) edge [->] node[descr,yshift=1.5ex] {$ \pi $} (m-1-3)
        (m-1-3) edge [->] node[descr,yshift=1.5ex] {$\beta_{\cG,X}$} (m-1-4)
        (m-2-1) edge [->] node[descr,yshift=1.5ex] {$i$} (m-2-2)
        (m-2-2) edge [->] node[descr,yshift=1.5ex] {$ \pi $} (m-2-3)
        (m-2-3) edge [->] node[descr,yshift=1.5ex] {$\beta_{\cJ,X}$} (m-2-4)
        (m-1-1) edge [->] node[descr,xshift=1.5ex] {$res$} (m-2-1) 
        
        (m-1-2) edge [->] node[descr,xshift=1.5ex] {$res$} (m-2-2) 
        
        (m-1-3) edge [->] node[descr,xshift=1.5ex] {$res$} (m-2-3) 
        
        (m-1-4) edge [->] node[descr,xshift=1.5ex] {$res$} (m-2-4) 
        ;
\end{tikzpicture} 
\] 
in which case 
\[
{res} \circ \beta_{\cG,X} = \beta_{\cJ,X} \circ {res} = 0 \circ {res} = 0
\]
this means $\beta_{\cG,X} = 0$ because $res$ is injective. This concludes the proof of Lemma \ref{Zero Bockstein}. \hfill $_\square$ \vspace{10pt}

With the previous discussion and the lemma proved, we have the framework necessary to show that all four modules have zero Bockstein maps for all three cases considered in this paper.

\newpage 
\section{The Cyclic Case} \label{Cyclic Section}

We begin our analysis with the cyclic case, namely when $E/F$ is cyclic Galois of degree $d$, and where we assume $\mu_{d^2}\subseteq F$.  
In this case the kernel $H^2(E/F,\mu_d)$ has been understood since
the early days of Class Field Theory. 
For suppose $E=F(\sqrt[d]{a})$. 
Then one has the well-known exact sequence
\[
H^1(E, \mu_d) \overset{{\rm N}_{E/F}}{\longrightarrow} H^1(F, \mu_{d}) \overset{(a)\smile}{\longrightarrow} H^2(F, \mu_{d}) \overset{i_{E/F}}{\longrightarrow} H^{2}(E, \mu_d) 
\]
which describes cohomology classes in $H^2(F, \mu_{d})$ that
vanish in $E$ as those corresponding to ``symbol algebras'' of the
form $(a,b)_F$ for some $b\in F$. 
This result also encodes
the classes $(a,b)_F$ which vanish in $H^2(F, \mu_{d})$
as those where $b\in {\rm N}_{E/F}(E^{\times})$, that is $b$
is a norm from $E$.  
This section generalizes this classical information to
all higher cohomology. 
The generalization of this result due to Voevodsky \cite{HW}, 
however the result there is the direct generalization of the four-term exact sequence in the presence of $d$th roots of unity, 
whereas the result obtained here extends this sequence one term to the right and left because the machinery of [P] gives a six-term sequence, but requiring $d^2$th roots of unity. 
Of course, the $H^2$ result just mentioned is a consequence of Hilbert's Theorem $90$ so this is not a surprise. 
Moreover, the generalization of 
Hilbert's Theorem 90 to higher K-theory is essential to Voevodsky's work, so this is also to be expected.

\subsection{The 4-Term Exact Sequence with Homotopies}

Let $E/F$ be a cyclic extension of degree $d$, 
with $d$th roots of unity 
$\mu_d \subseteq F$ and $\Gal(E/F) = G \cong \G/\cH = \langle \tau \rangle$. 
Let $F_\text{sep}$ denote the separable closure of $F$. 
Because $E/F$ is cyclic, the Galois closure $\widetilde{E}$ of $E/F$ is $E$, 
and $\widetilde{F} = F$. 
This simplifies the general setup.

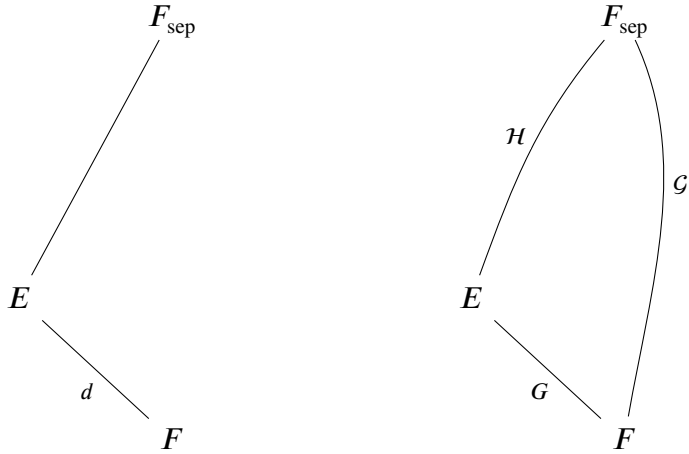
\[
\begin{tikzpicture}[descr/.style={fill=white,inner sep=1.5pt}]
        \matrix (m) [
            matrix of math nodes,
            row sep=3em,
            column sep=3em,
            text height=1.5ex, text depth=0.25ex
        ]
        {\ & F_{\text{sep}} & \ \\ 
         \ & \ & \                                        \\ 
            E & \ & \ \\ 
           \ & F & \                                        \\ \\ 
        };

        \path[overlay, font=\scriptsize,>=latex]
        (m-3-1) edge  node[descr,xshift=-1.75ex] {$\ $} (m-1-2) 
        (m-4-2) edge node[descr,xshift=-3,yshift=-1.5ex] {$d$ } (m-3-1) 
;
\end{tikzpicture} \hspace{.5in}
\begin{tikzpicture}[descr/.style={fill=white,inner sep=1.5pt}]
        \matrix (m) [
            matrix of math nodes,
            row sep=3em,
            column sep=3em,
            text height=1.5ex, text depth=0.25ex
        ]
        {\ & F_{\text{sep}} & \ \\ 
         \ & \ & \                                        \\ 
            E & \ & \ \\ 
           \ & F & \                                        \\ \\ 
        };

        \path[overlay, font=\scriptsize,>=latex]
        (m-4-2) edge node[descr,xshift=-3,yshift=-1.5ex] {$G$ } (m-3-1) 
        (m-3-1) edge[out=70,in=230] node[descr,xshift=-1ex,yshift=1ex] {$\cH$} (m-1-2) 
        (m-4-2) edge[out=80,in=-65,looseness = .9] node[descr,xshift=1.3ex,yshift=2.8ex] {$\cG$} (m-1-2) 
;
\end{tikzpicture} 
\]

It would be nice to use the restriction, corestriction sequence 
used in Arason's theorem, as was discussed in the introduction to Positselski's 6-term sequence, but the sequence 
\[
0 \longrightarrow \Z \overset{\Delta}{\longrightarrow} \Ind_{\cH}^\G(\Z) \overset{Tr}{\longrightarrow} \Z \longrightarrow 0 
\]
is not exact for every $d>2$. 
This is why we turn to Positselski's machinery, 
where two induced modules in the middle resolve this dimension problem. 
We use the 4-term exact sequence
\[
0 \longrightarrow \Z/d^2\Z \overset{\Delta}{\longrightarrow} \Ind_\cH^\G(\Z/d^2\Z) \overset{(1-\tau)\cdot}{\longrightarrow} \Ind_\cH^\G(\Z/d^2\Z) \overset{\text{Tr}}{\longrightarrow} \Z/d^2\Z \longrightarrow 0
\]
of $\G$-modules, where the maps $\Delta$ and Tr are defined in
Definition \ref{CyclicDefs} below.

In the notation of this sequence, since  $\mu_{d^2}\subseteq F$, 
the $\cG$-module $\mu_{d^2}$ will be identified with $\Z/d^2\Z$ and then
the module $\Z/d^2\Z$ acts as a vessel for the short exact sequence 
\[
0 \longrightarrow \mu_d \overset{\subseteq}{\longrightarrow} \mu_{d^2} \overset{\pi}{\longrightarrow} \mu_{d^2}/\mu_d \longrightarrow 0 
\]
making the $\G$-module homomorphisms chain maps. 
To do this, we identify the roots-of-unity short exact sequence with the additive short exact sequence 
\[
0 \longrightarrow d\Z/d^2\Z \overset{\subseteq}{\longrightarrow} \Z/d^2\Z \overset{\pi}{\longrightarrow} \Z/d\Z \longrightarrow 0.
\]

We also need a characterization of the induced module in order to facilitate computations in the four-term sequence given in Definition \ref{CyclicDefs} below.  
This is the subject of the following lemma, which allows us to do induced module computation in the group ring. 

\begin{lemma} \label{Induced Cyclic} 
We denote by $\tilde{\tau}$ be a lifting of $\tau$ to $F_\text{sep}$, and let $G = \langle \tau \rangle = \langle \tilde{\tau} \cH \rangle = \G/\cH$.  
Define $\phi:\Ind_\cH^\cG(\Z) \longrightarrow \Z[G]$ as follows: For any $f: \cG \longrightarrow \Z$ with the property $f(hg) = h\cdot f(g)$ for every $g \in \G, h \in \cH$, 
\[
\phi(f) = \sum_{i=0}^{d-1} f(\tilde{\tau}^i)\tau^{-i}
\]
Then $\phi$ is an isomorphism of $\cG$-modules, that is, 
$\Ind_\cH^\G(\Z) \cong \Z[G].$ 
\end{lemma} 
\noi {\bf Proof}: $\phi$ is a $\Z$-linear map that restricts to a bijection between $\Z$-bases for the two $\G$-modules, so it is a $\Z$-module isomorphism. So we need only check that the action is preserved. Every $g \in \G$ can be expressed as $g = h \tilde{\tau}^k$, where $h \in \cH$. Let $h_i = \tilde{\tau}^i h \tilde{\tau}^{-i} \in \cH$ for each $i \in \{0, \ldots, d-1 \}$ (these will be used to ``hop over'' the $\tilde{\tau}$ terms). Then 
\[
\begin{split} 
\phi(h \tilde{\tau}^k \cdot f) &= \sum_{i=0}^{d-1} (h\tilde{\tau}^k \cdot f)(\tilde{\tau}^i)\tau^{-i} = \sum_{i=0}^{d-1} f(\tilde{\tau}^i h \tilde{\tau}^k)\tau^{-i} = \sum_{i=0}^{d-1} f(h_i\tilde{\tau}^i \tilde{\tau}^k)\tau^{-i} \\ 
&= \sum_{i=0}^{d-1} h_i \cdot f(\tilde{\tau}^i \tilde{\tau}^k)\tau^{-i} = \sum_{i=0}^{d-1} f(\tilde{\tau}^i \tilde{\tau}^k)\tau^{-i} = \sum_{j=0}^{d-1} f(\tilde{\tau}^j )\tau^{-(j-k)} \\ 
&= \tau^k \sum_{j=0}^{d-1} f(\tilde{\tau}^j )\tau^{-j} = \tilde{\tau}^kh \cdot \phi(f) \\ 
\end{split} 
\]
\noi This concludes the proof of Lemma \ref{Induced Cyclic}. \hfill $_\square$ \vspace{10pt}

We next give the four term sequence in the cyclic case.
We will define the homotopies and verify the computational conditions for Positselski's 6-term sequence in the group ring $\Z[G]$
in Theorem \ref{Positselski Hypotheses Cyclic} below. 

\begin{definition}\label{CyclicDefs} Let $G$ be as above with
$\Z$ a trivial $G$-module and with $\Z[G]$ a $G$-module via 
multiplication on the left. We define $G$-module maps
\[ 
\begin{split} 
\Delta&: \Z \longrightarrow \Z[G], \ \ \ \ n \mapsto \mathop{\oplus}_{g \in G} ng, {\rm and}  \\ 
\text{Tr}&: \Z[G] \longrightarrow \Z, \ \ \mathop{\oplus}_{g \in G}c_g g \mapsto \sum_{g \in G} c_g.
\end{split} 
\] 
The Positselski modules
$M_1, M_2, M_3, M_4$ and maps $d_1, d_2, d_3$ for the
cyclic case are defined as follows: 
\[
0 \longrightarrow \Z \overset{\Delta}{\longrightarrow} \Z[G] \overset{\cdot(1-\tau)}{\longrightarrow} \Z[G] \overset{\text{Tr}}{\longrightarrow} \Z \longrightarrow 0.
\]
The homotopies
$h_1, h_2, h_3$ are defined as follows:  
\[
\Z \overset{\text{Tr}}{\longleftarrow} \Z[G] \overset{\cdot \sum-i\tau^i}{\longleftarrow} \Z[G] \overset{\Delta}{\longleftarrow} \Z.
\]
\end{definition} 

\begin{remark} 
The $\G$-module homomorphism 
$d_2 = (1-\tau): \Ind_\cH^\G(\Z) \longrightarrow \Ind_\cH^\G(\Z)$ 
is multiplication on the left by $(1-\tau)$ 
in the group ring $\Z[G]$ as in the Hilbert 90 sequence, 
which is also multiplication on the right by $(1-\tau)$, 
since $G=\langle \tau \rangle$ is abelian. 
We will use this last convention for the more general cases to follow. Another way of viewing this map is as the unique $\G$-module homomorphism that sends $1$ to $(1-\tau)$. 
The same remarks apply to
$h_2 = \cdot \sum -i \tau^i$, though the $h_2$ map will differ in the later cases. 
\end{remark} 

\begin{theorem} \label{Positselski Hypotheses Cyclic}
The $d_i$s are exact, and the $h_i$s satisfy the prism condition. 
\end{theorem} 
\noi \underline{\bf Proof}: Exactness follows from the usual Hilbert's Theorem 90 projective resolution argument for a cyclic extension. To show the prism condition,  
\[
\begin{split} 
h_2   d_2 = d_2   h_2 &= (1-\tau)\sum_{i=0}^{d-1} -i\tau^i \ = \ \sum_{i=0}^{d-1} -i\tau^i + \sum_{i=0}^{d-1} i\tau^{i+1} \\ 
&= -0\tau^0 + \sum_{i=1}^{d-1} -i\tau^i + \sum_{i=1}^{d-1} (i-1) \tau^i + (d-1)\tau^d \ = \ \sum_{i=1}^{d-1} -\tau^i + d
\end{split} 
\]
while $d_1   h_1 = h_3   d_3 = \sum_{i=1}^{d-1} \tau^i$. Therefore $d_1   h_1 + h_2   d_2 = h_2   d_2 + d_3   h_3 = d \cdot $. Furthermore, the fact that 
\[
\text{Tr}   \Delta= d \cdot
\]
follows from $d = [\G, \cH]$, and thus the prism condition for $h_2   d_1 = d_3   h_3$ is verified. \\ 
This concludes the proof of Theorem \ref{Positselski Hypotheses Cyclic}. \hfill $_\square$ \vspace{10pt}

The only remaining requirement to check for the Positselski hypotheses is the Bockstein maps being zero. 
In the cyclic case, all four modules are either $\Z$ or $\Ind_{\cH}^\cG(\Z)$, 
with each mod $d^2\Z$ to represent $\mu_{d^2}$ with a trivial action. 
Both of these moduoles have been shown to have a zero Bockstein map in Lemma \ref{Zero Bockstein}. 
Therefore the Positselski Hypotheses are satisfied by the 
four module exact sequence with homotopies defined in this section.

\subsection{The Connecting Map for the Cyclic Case}  
Now we compute the connecting map $\eta$. Let $\overline{M_i} := M_i/dM_{i}$ for $i \in \{1,2,3,4\}$ and let $c \in Z^{n-1}(F, \overline{M_4})$. We will begin the computation with a choice of lifting $\ell$ of $d_3$. For any $x \in \overline{M_4} = \Z/d\Z$, define $\ell(x) := x \cdot 1_G \in \overline{M_3} = \Z[G]/d\Z[G]$. 
Note that this is not the diagonal map, nor is it a $\G$-map. 
With $\ell$ chosen, we are ready to compute
the connecting map $\eta$ as the composition $d_1^{-1} h_2 \delta \ell$, 
where $\delta$ is the chain complex map from Galois cohomology 
for which the Bar Resolution is used. We used the customary choice of $\sigma$ for a cocycle in the Positselski framework, but we will reserve this symbol for the enlarged group $G$ in future sections. 
So let $c \in Z^{n-1}(\cG, \overline{M_4})$ denote our cocycle. We will compute the connecting map $\eta$ applied to $c$. 
What follows is a lemma that will reduce the complexity of this computation. 
\begin{lemma} \label{Cyclic Connecting Map} The following are true for $c$ and for any $g_1, \ldots, g_n \in \cG$. 
\begin{enumerate} 
\item $\delta(\ell(c))(g_1, \ldots, g_n) = -d_2\left( \sum_{i=0}^{k-1} \tau^i \cdot \ell(c(g_2, \ldots, g_n)) \right)$, \\ where $k \in \{ 0, \ldots, d-1 \}$ are such that $g_1 \cN = \tau^k$. 
\item $h_2   d_2 \equiv -d_1   h_1 $ (mod $dM_2)$. 
\end{enumerate} 
\end{lemma} 
\noi {\bf Proof}: (2) follows from the prism condition at $M_2$. 
To prove (1), we will start by using a similar argument to that used in the 
connecting map $\partial$ for Arason's theorem, 
using the fact that $\ell$ is an abelian group homomorphism to make the computation of $\delta(\ell(c))$ easier: 
\[
\begin{split}
\delta(\ell(c))(g_1, \ldots, g_n) &= g_1 \cdot \ell(c(g_2, \ldots, g_n)) - \ell(g_1 \cdot c(g_2, \ldots, g_n)) \\ 
&= \tau^k \cdot \ell(c(g_2, \ldots, g_n)) - \ell(c(g_2, \ldots, g_n)) \\ 
&= (\tau^k - 1) \cdot \ell(c(g_2, \ldots, g_n)) \\ 
&= -(1-\tau)\left( \sum_{i=0}^{k-1} \tau^i \cdot \ell(c(g_2, \ldots, g_n)) \right) \\ 
&= -d_2\left( \sum_{i=0}^{k-1} \tau^i \cdot \ell(c(g_2, \ldots, g_n)) \right). 
\end{split}
\]
\hfill $_\square$ \vspace{10pt}

This next Theorem describes the connecting map in the cyclic case.
\begin{theorem} \label{Connecting Map Cyclic Theorem} Let $\chi_\cH: \G \longrightarrow \Z/d\Z$ denote the character function that factors through the isomorphism $\G/\cH \overset{\cong}{\longrightarrow} \Z/d\Z$, namely $\chi_\cH(\tau^k) = k$. Let $\smile$ denote the cup product. Then
\[
\eta(c) = -\chi_\cH \smile c.
\]
\end{theorem} 
\noi {\bf Proof}: Since $\ell(c(g_2, \ldots, g_n) = c(g_2, \ldots, g_n)$, the equality 
\[
\tau^i \cdot c(g_2, \ldots, g_n) = c(g_2, \ldots, g_n) \tau^i
\]
follows. 
We next compute using parts (1) and (2) of Lemma \ref{Cyclic Connecting Map} as indicated. 
\[
\begin{split} 
-\eta(c)(g_1, \ldots, g_n) &= d_1^{-1} h_2 \delta \ell(c)(g_1, \ldots, g_n) \ \overset{(1)}{=} \ -d_1^{-1} h_2 d_2\left( \sum_{i=0}^{k-1} \tau^i \cdot \ell(c(g_2, \ldots, g_n)) \right) \\ 
&= -d_1^{-1} h_2 d_2\left( \sum_{i=0}^{k-1} c(g_2, \ldots, g_n) \tau^i \right) \ \overset{(2)}{=} \ d_1^{-1} d_1 h_1\left( \sum_{i=0}^{k-1} c(g_2, \ldots, g_n) \tau^i \right) \\ 
&= h_1\left( \sum_{i=0}^{k-1} c(g_2, \ldots, g_n) \tau^i \right) \ = \ Tr\left( \sum_{i=0}^{k-1} c(g_2, \ldots, g_n) \tau^i\cH \right) \\ 
&= \sum_{i=0}^{k-1} c(g_2, \ldots, g_n) \ = \ k c(g_2, \ldots, g_n) \ = \ (\chi_\cH \smile c)(g_1, \ldots, g_n). \\ 
\end{split} 
\]
This concludes the proof of Theorem \ref{Connecting Map Cyclic Theorem}. \hfill $_\square$ \vspace{10pt}

In view of Theorems \ref{Positselski Hypotheses Cyclic} 
and  \ref{Connecting Map Cyclic Theorem} the machinery in Theorem \ref{Positselski's Theorem} gives the following result.
\begin{theorem} \label{Connecting Map Cyclic 6-term Sequence Theorem} 
In the cyclic case we have the following 6-term exact sequence. 
\[
\begin{tikzpicture}[descr/.style={fill=white,inner sep=1.5pt}]
        \matrix (m) [
            matrix of math nodes,
            row sep=3em,
            column sep=3em,
            text height=1.5ex, text depth=0.25ex
        ]
        { H^n(E,\mu_d)\oplus H^n(F,\mu_d) & H^{n+1}(E,\mu_d) & H^n
        (F,\mu_d)   \\ 
        H^{n+1}(F,\mu_d) & H^{n+1}(E,\mu_d) & H^n(F,\mu_d)\oplus H^n(E,\mu_d) \\ 
        };

        \path[overlay, font=\scriptsize,>=latex]
        (m-1-1) edge [->]  node[descr,yshift=1.5ex] {$d_2 + h_3 $} (m-1-2) 
        (m-1-2) edge [->]  node[descr,yshift=1.5ex] {$d_3 $} (m-1-3) 
        (m-2-1) edge [->]  node[descr,yshift=1.5ex] {$d_1 $} (m-2-2) 
        (m-2-2) edge [->]  node[descr,yshift=1.5ex] {$h_1 \oplus d_2$} (m-2-3) 
        (m-1-3) edge[out=355,in=175,->] node[descr,yshift=0.3ex] {$\eta$} (m-2-1)
        
;
\end{tikzpicture} 
\] 
where $d_3$ is the norm, $\eta(c) = -\chi \smile c$ and $d_1$ is scalar extension.
\end{theorem}

\section{The Dihedral Case} \label{The Dihedral Case}
We noted in the introduction that the case where
$[E:F]=4$ and $\Gal(\tE/F)$ is the dihedral group of order $8$
was handled by Positselski in \cite{Positselski}. In fact, Positselski handled every case where $[E:F]$ is a multiple of 4.  
In this section we turn to the dihedral cases where  $d=[E:F]$ is odd. 
\subsection{The 4 Term Exact Sequence with Homotpies} 
For this section, we have the following notation. Let $G$ be a dihedral group, $G = \langle \sigma, \tau \rangle$ with $|\tau| = d$ for some odd integer $d$, $|\sigma| = 2$ and the relation $\sigma \tau = \tau^{-1} \sigma$.
We use the notation described earlier. In this case the diagrams of fields and groups are as follows. 
\[
\begin{tikzpicture}[descr/.style={fill=white,inner sep=1.5pt}]
        \matrix (m) [
            matrix of math nodes,
            row sep=3em,
            column sep=3em,
            text height=1.5ex, text depth=0.25ex
        ]
        {\ & F_{\text{sep}} & \ \\ 
         \ & \widetilde{E} & \                                        \\ 
            E & \ & \widetilde{F} \\ 
           \ & F & \                                        \\ \\ 
        };

        \path[overlay, font=\scriptsize,>=latex]
        (m-2-2) edge  node[descr,xshift=-1.75ex] {$\ $} (m-1-2) 
        (m-3-1) edge  node[descr,xshift=-3,yshift=1.1ex] {$2$}(m-2-2) 
        (m-3-3) edge node[descr,xshift=3,yshift=1.5ex] {$d$ }(m-2-2) 
        (m-4-2) edge node[descr,xshift=-3,yshift=-1.5ex] {$d$ } (m-3-1) 
        (m-4-2) edge node[descr,xshift=3,yshift=-1.6ex] {$2$ } (m-3-3) 
;
\end{tikzpicture} \hspace{.5in}
\begin{tikzpicture}[descr/.style={fill=white,inner sep=1.5pt}]
        \matrix (m) [
            matrix of math nodes,
            row sep=3em,
            column sep=3em,
            text height=1.5ex, text depth=0.25ex
        ]
        {\ & F_{\text{sep}} & \ \\ 
         \ & \widetilde{E} & \                                        \\ 
            E & \ & \widetilde{F} \\ 
           \ & F & \                                        \\ \\ 
        };

        \path[overlay, font=\scriptsize,>=latex]
        (m-2-2) edge  node[descr,xshift=-1.25ex,yshift=-2] {$\cN$} (m-1-2) 
        (m-3-1) edge[out=70,in=230] node[descr,xshift=-1ex,yshift=1ex] {$\cH$} (m-1-2) 
        (m-3-3) edge[out=110,in=310] node[descr,xshift=1ex,yshift=1ex] {$\cJ$} (m-1-2) 
        (m-3-1) edge  node[descr,xshift=-2.1,yshift=1.75ex] {$\langle \sigma \rangle$} (m-2-2) 
        (m-3-3) edge  node[descr,xshift=2.1,yshift=1.75ex] {$\langle \tau \rangle$} (m-2-2) 
        (m-4-2) edge[out=13,in=-30,looseness = 1.7] node[descr,xshift=1ex,yshift=1ex] {$\cG$} (m-1-2) 
        (m-4-2) edge node[descr,xshift=-7,yshift=0ex] {$G$} (m-2-2) 
        (m-4-2) edge (m-3-1) 
        (m-4-2) edge (m-3-3) 
;
\end{tikzpicture} 
\] 

In the previous section in the cyclic case with $G = \langle \tau \rangle$ the two short exact sequences that make up our 4-term sequence of modules come from adjacent terms in this projective resolution. 
\[
\cdots \longrightarrow \Z[G] \overset{\cdot(1-\tau)}{\longrightarrow} \Z[G] \overset{\cdot T_\tau}{\longrightarrow} \Z[G] \overset{\cdot(1-\tau)}{\longrightarrow} \Z[G] \overset{\cdot T_\tau}{\longrightarrow} \Z[G] T_\tau \longrightarrow 0
\]
In this resolution, $T_\tau=1 + \tau + \cdots + \tau^{d-1}$ is the $\tau$-trace, though it will act as a norm map on roots of unity in the field. 
This projective resolution is commonly used in a proof of Hilbert's Theorem 90 for cyclic Galois extensions. 
The alternating short exact sequences are
\[
0 \longrightarrow \Z[G]T_\tau \overset{\subseteq}{\longrightarrow} \Z[G] \overset{\cdot(1-\tau)}{\longrightarrow} \Z[G] \longrightarrow 0
\]
and 
\[
0 \longrightarrow \Z[G] (1-\tau)  \overset{\subseteq}{\longrightarrow} \Z[G] \overset{\cdot T_\tau}{\longrightarrow} \Z[G] T_\tau \longrightarrow 0.
\]
For the cyclic case where $d=2$ these two short exact sequences are the same because $T_\tau = 1-\tau$. Furthermore, both short exact sequences begin and end with $\Z/2\Z$ as a trivial $G$-module. This avoids the need for Positselski's machinery altogether, as both short exact sequences are the Arason sequence. 

The four term exact sequence of modules for the dihedral case is also similar to that of the cyclic case. The two short exact sequences that define it are
\[
0 \longrightarrow \Z[G]T_\sigma T_\tau \overset{\subseteq}{\longrightarrow} \Z[G]T_\sigma \overset{\cdot(1-\tau)}{\longrightarrow} \Z[G]T_\sigma(1-\tau) \longrightarrow 0
\]
and 
\[
0 \longrightarrow 
\Z[G] 
(1-\tau)
\cB 
\overset{\subseteq}{\longrightarrow} 
\Z[G]\cB \overset{\cdot T_\tau}{\longrightarrow} 
\Z[G] T_\tau\cB \longrightarrow 0
\]
where $T_\sigma=(1+\sigma)$ is the $\sigma$-trace, and $\cB = (1-\sigma\tau) \in \Z[G]$. Of course, we need $\Z[G]T_\sigma(1-\tau) \cong \Z[G](1-\tau)\cB$ 
for the two short exact sequences to build a 4-term exact sequence, and in fact $T_\sigma(1-\tau) = (1-\tau)\cB$. 
Thus, the 4-term exact sequence is
\[
0 \longrightarrow \Z[G]T_\sigma T_\tau \overset{\subseteq}{\longrightarrow} \Z[G]T_\sigma \overset{\cdot(1-\tau)}{\longrightarrow} \Z[G]\cB \overset{\cdot T_\tau}{\longrightarrow} \Z[G]\cB T_\tau \longrightarrow 0.
\]

Our first lemma describes the induced module $\Ind_\cH^\cG(\Z)$ needed for this case and relates it to $\Z[G]T_\sigma$ for computation.

\begin{lemma} \label{Induced Dihedral} Let $\tilde{\sigma}$ and $\tilde{\tau}$ be liftings of $\sigma$ and $\tau$ from $G$ to $\G$. Then for the trivial $\Z[\G]$-module $\Z$, the map $\phi: \Ind_\cH^\G(\Z) \longrightarrow \Z[G]T_\sigma$ given by 
\[
\phi(f) = \sum_{i=0}^{d-1} f(\tilde{\tau}^i)\tau^{-i}T_\sigma
\]
is an isomorphism of $\Z[\G]$-modules. 
\end{lemma} 
\noi {\bf Proof}: The proof is similar to that for the cyclic case,
but it is worth noting that this proof would not work for $\Z[\langle \tau \rangle]$ instead of $\Z[G]T_\sigma$ 
because $\cH$ is not a normal subgroup of $\G$. 
However, $\cN$ is normal in $\G$ and we will use this fact. 

The bijectivity of $\phi$ follows as in Lemma \ref{Induced Cyclic} and we proceed to check the compatibility of $\G$-actions. 
Every element of $\G$ has a unique expression as $n\tilde{\sigma}^m\tilde{\tau}^k$, with $n \in \cN$, $m \in \Z/2\Z, k \in \Z/d\Z$. 
Let $\theta_m: \Z/d\Z \longrightarrow \Z/d\Z$ be the $\sigma^m$-conjugation automorphisms so that $\sigma^m \tau^i = \tau^{\theta_m(i)} \sigma^m$ and $\tau^i \sigma^m = \sigma^m \tau^{\theta_{-m}(i)}$. Finally, let $n_{i,m} = \tilde{\tau}^i n \tilde{\sigma}^m \tilde{\tau}^{\sigma_{-m}(-i)} \tilde{\sigma}^{-m} \in \cN$ so that $ \tilde{\tau}^i n \tilde{\sigma}^m = n_{i,m}\tilde{\sigma}^m \tilde{\tau}^{\theta_{-m}(i)}$. \\
Then 
\[
\begin{split} 
\phi(n \tilde{\sigma}^m \tilde{\tau}^k \cdot f) &= \sum_{i=0}^{d-1} ( n\tilde{\sigma}^m\tilde{\tau}^k \cdot f)(\tau^i)\tau^{-i}T_\sigma \ = \ \sum_{i=0}^{d-1} f(\tilde{\tau}^i n\tilde{\sigma}^m\tau^k)\tau^{-i}T_\sigma \\ 
&= \sum_{i=0}^{d-1} f(n_{i,m}\tilde{\sigma}^m \tilde{\tau}^{\theta_{-m}(i)} \tau^k)\tau^{-i}T_\sigma \ = \ \sum_{i=0}^{d-1} (n_{i,m}\tilde{\sigma}^m) \cdot f(\tilde{\tau}^{\theta_{-m}(i)} \tau^k)\tau^{-i}T_\sigma \\ 
&= \sum_{i=0}^{d-1} f(\tilde{\tau}^{\theta_{-m}(i)+k})\tau^{-i}T_\sigma \ = \ \sum_{j=0}^{d-1} f(\tilde{\tau}^j)\tau^{-\theta_m(j-k)}T_\sigma \ = \ \sum_{j=0}^{d-1} f(\tilde{\tau}^j)\tau^{-\theta_m(j-k)}\sigma^m T_\sigma \\ 
&= \sum_{j=0}^{d-1} f(\tilde{\tau}^j) \sigma^m \tau^{-(j-k)}T_\sigma \ = \ \sum_{j=0}^{d-1} f(\tilde{\tau}^j) \sigma^m \tau^k \tau^{-j}T_\sigma \ = \ (\sigma^m \tau^k) \cdot \phi(f) \\ 
&= (n \sigma^m \tau^k) \cdot \phi(f). \\ 
\end{split} 
\]
With the compatability of $\Z[\G]$-action checked, this completes the proof of Lemma \ref{Induced Dihedral}. \hfill  $_\square$ 
\vspace{10pt}

\begin{remark}\label{sigma doesn't need to have order 2}
The fact that $|\sigma| = 2$ was not used in this proof. 
And indeed, this same proof can be used in the analogous claim for the more general semi-direct case later. 
We will therefore refer to this lemma for the semi-direct case as well. 
\end{remark}

With $\Z[G]T_\sigma \cong \Ind_\cH^\cG(\Z)$ established, we move on to defining the homomorphisms. 
\begin{definition} 
Let $d_1, d_2, d_3$ be the maps 
\begin{equation} \label{Dihedral d's} 0 \longrightarrow \Z[G]T_\sigma T_\tau \overset{\subseteq}{\longrightarrow} \Z[G]T_\sigma \overset{\cdot(1-\tau)}{\longrightarrow} \Z[G]\cB \overset{\cdot T_\tau}{\longrightarrow} \Z[G]\cB T_\tau \longrightarrow 0 
\end{equation} 
Let $h_1, h_2, h_3$ be the homotopy maps  
\begin{equation} \label{Dihedral h's}
\Z[G]T_\sigma T_\tau \overset{\cdot T_\tau}{\longleftarrow} \Z[G]T_\sigma \overset{h_2}{\longleftarrow} \Z[G]\cB \overset{\supseteq}{\longleftarrow} \Z[G]\cB T_\tau 
\end{equation} 
where 
\[h_2(\cB) = \sum_{i=0}^{d-1} \left( \frac{d-1}{2} - i \right) \tau^i T_\sigma. \] \\ 
\end{definition} 

We observe that $\Z[G]T_\sigma T_\tau = \Z[G] \left( \sum_{g \in G} g \right) \cong \Z$ is a trivial $\Z[G]$-module, while $\Z[G]\cB T_\tau$ has a trivial $\tau$-action but $\sigma$ and $\sigma\tau$ act as multiplication by $(-1)$. 

The main result needed for applying Positselski's machinery is given next.

\begin{theorem} \label{Dihedral d's and h's}
The sequence (\ref{Dihedral d's}) is an exact sequence, and the homotopy maps in sequence (\ref{Dihedral h's}) satisfy the prism condition. Furthermore, $d_1, d_2, d_3, h_1, h_2, h_3$ are $\Z[G]$-module homomorphisms. 
\end{theorem} 

\noi {\bf Proof}: The exactness of the $d_i$'s was discussed at the beginning of this section, so we move on to checking the prism condition. For the first and last modules, 
the composition is multiplication by $d$ because $T_\tau T_\tau = d T_\tau$. For the second module, $d_1   h_1(T_\sigma) = T_\tau T_\sigma$ and 
\[
\begin{split} 
h_2   d_2(T_\sigma) &= h_2(T_\sigma(1-\tau)) \ = \ h_2((1-\tau)\cB) \ = \ (1-\tau) \sum_{i=0}^{d-1} \left( \frac{d-1}{2} - i \right) \tau^i T_\sigma \\ 
&= (1-\tau) \sum_{i=0}^{d-1} \left( - i \right) \tau^i T_\sigma \ = \ (d - T_\tau)T_\sigma
\end{split} 
\]
Therefore $(d_1h_1 + h_2d_2)(T_\sigma) = T_\tau T_\sigma + (d - T_\tau)T_\sigma = dT_\sigma$. To verify the prism condition for the third module, $h_3   d_3(\cB) = h_3(\cB T_\tau) = \cB T_\tau = T_\tau \cB$; and using the identity $T_\sigma (1-\tau) = (1-\tau)\cB$, 
\[
\begin{split} 
d_2   h_2 (\cB) &= d_2 \left( \sum_{i=0}^{d-1} \left( \frac{d-1}{2} - i \right) \tau^i T_\sigma \right) \ = \ \sum_{i=0}^{d-1} \left( \frac{d-1}{2} - i \right) \tau^i T_\sigma (1-\tau) \\ 
&= \sum_{i=0}^{d-1} \left( \frac{d-1}{2} - i \right) \tau^i (1-\tau) \cB \ = \ \sum_{i=0}^{d-1} \left( - i \right) \tau^i (1-\tau) \cB \\ 
&= (d-T_\tau) \cB. \\ 
\end{split} 
\]
Therefore $(h_3d_3 + d_2h_2)(\cB) = T_\tau\cB+(d-T_\tau)\cB = d\cB$. 

Every one of these maps except for $h_2$ is defined by left multiplication in the group ring $\Z[G]$, so it remains to check that $h_2$ preserves action by $\G$. $h_2$ is also a map from a cyclic $\Z[G]$-module to a cyclic $\Z[G]$-module defined by sending one generator to another. So it need only be checked that the annihilator of $\cB$, which is $\Z[G](1 + \sigma\tau)$, also annihilates the $h_2$-image of $\cB$. We will check this by showing that $h_2$ preserves the action of $\sigma \tau$.  
\[
\begin{split}  
\sigma \tau \cdot h_2( \cB ) &= \sigma \tau \sum_{i=0}^{d-1} \left( \frac{d-1}{2} - i \right) \tau^i T_\sigma = \sigma \sum_{i=0}^{d-1} \left( \frac{d-1}{2} - i \right) \tau^{i+1} T_\sigma \\ 
&= \sum_{i=0}^{d-1} \left( \frac{d-1}{2} - i \right) \tau^{d-(i+1)} \sigma T_\sigma = \sum_{i=0}^{d-1} \left( \frac{d-1}{2} - i \right) \tau^{(d-1)-i} T_\sigma \\ 
&= \sum_{i=0}^{d-1} \left( - \frac{d-1}{2} + \left( (d-1) - i \right) \right) \tau^{(d-1)-i} T_\sigma = \sum_{j=0}^{d-1} \left( - \frac{d-1}{2} + j \right) \tau^i T_\sigma \\ 
&= - \sum_{j=0}^{d-1} \left( \frac{d-1}{2} - j \right) \tau^i T_\sigma = h_2 \left( - \cB \right) = h_2 \left( \sigma \tau \cdot \cB \right) 
\end{split} 
\]
This concludes the proof of Theorem \ref{Dihedral d's and h's}. \hfill $_\square$ \vspace{10pt}

The only remaining requirements to check for the Positselski hypotheses are that the Bockstein homomorphisms are zero. 
The first two modules are isomorphic to $\Z$ and the induced module $\Ind_{\cH}^\cG(\Z)$ respectively. 
Both of these modules were shown to have zero Bockstein maps in the general setup section. 
For $M_3$ and $M_4$, we will show that as $\cJ$-modules, $M_3 \cong \Z[J] \cong \Ind_\cN^\cJ(\Z)$ and that $M_4 \cong \Z$ as a trivial $\cJ$-module. 
Then Lemma \ref{Zero Bockstein} will imply that the Bockstein map vanishes for $M_3$ and $M_4$. 
The isomorphisms to be defined have the property that $\cB \mapsto 1 \in \Z[\cJ]$ and $\cB T_\tau \mapsto 1 \in \Z$ for $M_3$ and $M_4$ respectively.  

\begin{lemma} \label{Dihedral Zero Bockstein}  
As $\cJ$-modules, 
\begin{enumerate} 
\item $M_3 = \Z[G]\cB \cong \Ind_\cN^\cJ(\Z)$ 
\item $M_4 = \Z[G]\cB T_\tau \cong \Ind_\cJ^\cJ(\Z)$ 
\end{enumerate} 
\end{lemma}

\noi {\bf Proof}: There are two steps.
\begin{enumerate} 
\item We begin with the observation that 
\[
M_3 = \Z[G] \cdot (1-\sigma \tau) = \Z[\langle \tau \rangle ] \cdot (1 - \sigma \tau).
\] 
The second equality above follows from 
\[
\tau^i\sigma(1-\sigma\tau) = \tau^{i+1} \cdot \sigma \tau(1-\sigma\tau) = -\tau^{i+1}(1-\sigma\tau)
\]
and hence 
\[
\left( \sum_{i=0}^{d-1}c_{i,0} \tau^i + \sum_{i=0}^{d-1}c_{i,1}\tau^i\sigma \right)(1-\sigma\tau) = \left( \sum_{i=0}^{d-1}(c_{i,0} - c_{i,1}\tau)\tau^i \right)(1-\sigma\tau) \in\Z[\langle \tau \rangle](1-\sigma\tau).
\]
Furthermore, the set $\{ \tau^i(1-\sigma \tau \ | \ 0 \leq i \leq d-1 \}$ is a $\Z$-basis for $M_3$, with $\tau^j \cdot \tau^i(1-\sigma\tau) = \tau^{j+i}(1-\sigma\tau)$. 
\item Similarly 
\[
M_4 = \Z[G](1-\sigma \tau)T_\tau = \Z[\langle \tau \rangle] (1-\sigma\tau)T_\tau = \Z[\langle \tau \rangle] T_\tau (1-\sigma\tau) = \Z T_\tau (1-\sigma \tau) = \Z (1-\sigma\tau)T_\tau
\]
and $\Z(1-\sigma\tau)T_\tau \cong \Z$ via the $\cJ$-module isomorphism $(1-\sigma\tau)T_\tau \mapsto 1$. 
\end{enumerate}
This concludes the proof of Lemma \ref{Dihedral Zero Bockstein}. \hfill $_\square$ \vspace{10pt}

Now, with $M_3$ isomorphic to $\Z[J]$ and $M_4$ isomorphic to $\Z$ as $\cJ$-modules, the fact that $\widetilde{E}/\widetilde{F}$ is a cyclic extension allows the application of Lemma \ref{Zero Bockstein} to reduce the problem to the cyclic case. This makes the Bockstein maps zero for $M_3$ and $M_4$ as well as $M_1$ and $M_2$. 
Therefore the Positselski Hypotheses are satisfied by the four module exact sequence with homotopies defined in this section.

\subsection{The Connecting Map for the Dihedral Case} 
Now we compute the connecting map. Let $\overline{M}_i := M_i/dM_{i}$ for $i \in \{1,2,3,4\}$. 
Note that Lemmas 
\ref{Induced Dihedral} and 
\ref{Dihedral Zero Bockstein}
extend to the same modules mod $d^2$ and mod $d$. 
In this section we will use the exact sequence of modules with homotopies defined in the previous section to describe the connecting map $\eta: H^{n-1}(\cG, \overline{M}_4) \longrightarrow H^n(\cG, \overline{M}_1)$ given by Positselski's machinery.   
\begin{definition}: With this notation we define the following. \begin{enumerate} 
\item $\ell: \overline{M}_4 \longrightarrow \overline{M}_3$, the $d_3$-lifting defined as $\ell(z \cB T_\tau) = z \cB$ for every $z \in \Z/d\Z$. 
\item $ \delta: C^{n-1}(\cG, \overline{M}_3) \longrightarrow C^n(\cG, \overline{M}_3)$ the cochain map from the bar resolution.
\item $\widetilde{\eta}_{\ell} = \widetilde{\eta} : Z^{n-1}(\cG, \overline{M}_4) \longrightarrow Z^n(\cG, \overline{M}_1)$, \ $\widetilde{\eta}(c) := d_1^{-1} h_2 \delta \ell(c)$. 
\item $\eta: H^{n-1}(\cG, \overline{M}_4) \longrightarrow H^n(\cG, \overline{M}_1)$, \ $\eta([c]) := [\widetilde{\eta}(c)]$. 
\end{enumerate} 
\end{definition}

We observe that our choice of lifting $\ell$ is a $\Z$-module homomorphism, though it is not a $\Z[G]$-module homomorphism. Note, for $x \in \im(d_1)$, we let $d_1^{-1}(x)$ denote the unique preimage element. 

We collect some basic properties of these maps next.

\begin{lemma} \label{Connecting Map Dihedral Lemma}
Let $c \in Z^{n-1}(\cG, \overline{M}_4)$ be a cocycle, $g_1, \ldots, g_n \in \cG$, and let $c' \in \Z/d\Z$ such that $c(g_2, \ldots, g_n) = c' \cdot \cB T_\tau$. 
Express the coset of $g_1$ in $\cG/\cN$ as $(\sigma\tau)^j \tau^i$. 
Then  
\begin{enumerate} 
\item $\delta(\ell(c))(g_1, \ldots, g_n) = g_1 \cdot \ell(c(g_2, \ldots, g_n)) - \ell(g_1 \cdot c(g_2, \ldots, g_n))$ \\ 
$= (\sigma\tau)^j \left( -\sum_{k=0}^{i-1} \tau^k \right) \cdot c' (1-\tau)\cB $. 
\item For any $x \in C^n(\cG, \overline{M}_2)$ such that $\delta(\ell(c)) = d_2(x)$, $\widetilde{\eta}(c) = h_1(x).$
\end{enumerate} 
\end{lemma} 

\underline{\bf Proof}: For the first part, by definition, 
\[
\delta(\ell(c))(g_1, \ldots, g_n) = g_1 \cdot \ell(c(g_2, \ldots, g_n)) - \ell(c(g_1 g_2, \ldots, g_n)) + \cdots \pm \ell(c(g_1, \ldots, g_{n-1}))
\]
We will next subtract an expanded form of $\ell(\delta(c))(g_1, \ldots, g_n)$ from the right side of this equation. 
This term is 0 because $\ell$ is a $\Z$-homomorphism and $c$ is a cocycle. Subtracting the expanded form will leave us with only two remaining terms. The expansion is as follows.  \\ 
\[
\begin{split}
&\text{ \ } \ \ \ell(\delta(c))(g_1, \ldots, g_n) \\ 
&= \ell(\delta(c)(g_1, \ldots, g_n)) \\ 
&= \ell \left( g_1 \cdot c(g_2, \ldots, g_n) - c(g_1g_2, \ldots, g_n) + \cdots \pm c(g_1, \ldots, g_{n-1})\right) \\  
&= \ell(g_1 \cdot c(g_2, \ldots, g_n)) - \ell(c(g_1g_2, \ldots, g_n)) + \cdots \pm \ell(c(g_1, \ldots, g_{n-1}))) \\ 
\end{split} 
\]
We can now subtract and simplify 
\[
\begin{split} 
&\text{ \ } \ \ \delta(\ell(c))(g_1, \ldots, g_n) \\ 
&\text{ \ } \ \ - 0 \\ 
&= \delta(\ell(c))(g_1, \ldots, g_n)  \\ 
&\text{ \ } \ \ - \ell(\delta(c))(g_1, \ldots, g_n) \\ 
&= g_1 \cdot \ell(c(g_2, \ldots, g_n)) - \ell(c(g_1 g_2, \ldots, g_n)) + \cdots \pm \ell(c(g_1, \ldots, g_{n-1})) \\ 
&\text{ \ } \ \ - \ell(g_1 \cdot c(g_2, \ldots, g_n)) + \ell(c(g_1 g_2, \ldots, g_n)) - \cdots \mp \ell(c(g_1, \ldots, g_{n-1})) \\ 
&= g_1 \cdot \ell(c(g_2, \ldots, g_n))  \\ 
&\text{ \ } \ \ - \ell(g_1 \cdot c(g_2, \ldots, g_n)) 
\end{split}
\]
This shows the first equality in part $(1)$ of the lemma. \\ \\ 
For the next equality, we replace $g_1$ with $(\sigma\tau)^j \tau^i$, $g_1$'s coset representative in $\cG/\cN$, and use the fact that $\ell$ preserves the action by $\sigma \tau$, which is also multiplication by $(-1)$ for both $\cB \in \overline{M}_3$ and $\cB T_\tau \in \overline{M}_4$.  
\[
\begin{split} 
\delta(\ell(c))(g_1, \ldots, g_n) 
&= g_1 \cdot \ell(c(g_2, \ldots, g_n)) - \ell(g_1 \cdot c(g_2, \ldots, g_n)) \\ 
&= (\sigma\tau)^j \tau^i \cdot \ell(c(g_2, \ldots, g_n)) - \ell((\sigma\tau)^j \tau^i \cdot c(g_2, \ldots, g_n)) \\ 
&= (\sigma\tau)^j \tau^i \cdot \ell(c(g_2, \ldots, g_n)) - \ell((\sigma\tau)^j c(g_2, \ldots, g_n)) \\ 
&= (\sigma\tau)^j \tau^i \cdot \ell(c(g_2, \ldots, g_n)) - \ell((-1)^j c(g_2, \ldots, g_n)) \\ 
&= (\sigma\tau)^j \tau^i \cdot \ell(c(g_2, \ldots, g_n)) - (-1)^j \ell( c(g_2, \ldots, g_n)) \\ 
&= (\sigma\tau)^j \tau^i \cdot \ell(c(g_2, \ldots, g_n)) - (\sigma\tau)^j \ell(c(g_2, \ldots, g_n)) \\ 
&= (\sigma\tau)^j \left( \tau^i \cdot \ell(c(g_2, \ldots, g_n)) - \ell(c(g_2, \ldots, g_n)) \right) \\ 
&= (\sigma\tau)^j (\tau^i - 1) \cdot \ell(c(g_2, \ldots, g_n)) \\ 
&= (\sigma\tau)^j (\tau^i - 1) \cdot c' \cB \\ 
&= (\sigma\tau)^j \left( -\sum_{k=0}^{i-1} \tau^k \right) (1 - \tau) \cdot c' \cB \\ 
&= (\sigma\tau)^j \left( -\sum_{k=0}^{i-1} \tau^k \right) \cdot c' (1-\tau)\cB \\ 
\end{split} 
\]
This proves the first statement of the lemma. 

The second part of the lemma follows from the prism condition. Modulo $d$,  
\[
d_1   h_1 + h_2   d_2 = \cdot d = 0
\]
This means 
\[
h_2   d_2 = -d_1   h_1
\]
And hence 
\[
\begin{split}
\eta(c) &= -(d_1^{-1}) (h_2(\delta(\ell(c))) \\
&= -(d_1^{-1}) (h_2(d_2(x))) \\ 
&= -(d_1^{-1}) (-d_1 (h_1(x))) \\ 
&= h_1(x) \\ 
\end{split} 
\]

This concludes the proof of Lemma \ref{Connecting Map Dihedral Lemma}. \hfill $_\square$ \vspace{10pt}

As an application we obtain a description of the connecting map $\eta$
in this case.  Although it is not a cup product, its description is
almost one.
 
\begin{corollary}  \label{Connecting Map Dihedral Corollary} Let $c \in Z^{n-1}(\cG, \overline{M}_4)$, $g_1, g_2, \ldots, g_n \in \cG$, $c' \in \Z/d\Z$ such that $c(g_2, \ldots, g_n) = c' \cB T_\tau$. 
Let $\sigma^i \tau^j$, an element of $\cG/\cN$, be the coset of $g_1$. 
Then \\ 
\[
\widetilde{\eta}(c)(g_1, \ldots, g_n) = -i c' \in \overline{M}_1.
\]
\end{corollary} 
\noi \underline{Proof}: From part (1) of Lemma \ref{Connecting Map Dihedral Lemma}, we know that 
\[
\delta(\ell(c))(g_1, \ldots, g_n) = (\sigma\tau)^j \left( -\sum_{k=0}^{i-1} \tau^k \right) \cdot c' (1-\tau)\cB.
\]
We will first find an $x \in \overline{M}_2$ with this $d_2$-image, and then use part (2) of Lemma \ref{Connecting Map Dihedral Lemma} to compute $\widetilde{\eta}(c)$ by finding $-h_1(x)$. 
This process begins by showing that $(\sigma\tau)^j \left( -\sum_{k=0}^{i-1} \tau^k \right) \cdot c' (1+\sigma)$ is a suitable choice for $x$. 
\[
\begin{split} 
\delta(\ell(c))(g_1, \ldots, g_n) &= (\sigma\tau)^j \left( -\sum_{k=0}^{i-1} \tau^k \right) \cdot c' (1-\tau)\cB \\ 
&= (\sigma\tau)^j \left( -\sum_{k=0}^{i-1} \tau^k \right) \cdot c' (1+\sigma)(1-\tau) \\
&= (\sigma\tau)^j \left( -\sum_{k=0}^{i-1} \tau^k \right) \cdot c' (1+\sigma)(1-\tau) \\ 
&= d_2\left( (\sigma\tau)^j \left( -\sum_{k=0}^{i-1} \tau^k \right) \cdot c' (1+\sigma) \right) \\ 
\end{split} 
\]
Applying part $(b)$ allows the computation of $\eta(c)(g_1, \ldots, g_n)$ as follows. 
\[
\begin{split} 
\widetilde{\eta}(c)(g_1, \ldots, g_n) &= h_1\left((\sigma\tau)^j \left( -\sum_{k=0}^{i-1} \tau^k \right) \cdot c' (1+\sigma) \right) \\
&= (\sigma\tau)^j \left( -\sum_{k=0}^{i-1} \tau^k\right) c' \cdot h_1\left(1+\sigma \right) \\ 
&= (\sigma\tau)^j \left( -\sum_{k=0}^{i-1} \tau^k\right) c' \cdot 1 \\ 
&= -i c' \cdot 1 \\ 
&= -i c' \\ 
\end{split} 
\]
This concludes the proof of Corollary \ref{Connecting Map Dihedral Corollary} . \hfill $_\square$  \vspace{10pt}

Let $\chi: \cG \longrightarrow \Z/d\Z$
be defined by $\chi((\sigma \tau)^j \tau^i \cN) = i$.  It should be noted that in this case $\chi$ is no longer a character because it is not a homomorphism. However, the Corollary shows that the map $\eta$ can be understood through a cup-product like structure which we
denote by $\smile'$ whereby $(\chi \smile' c)(g_1, \ldots, g_n)=
\chi(g_1)\cdot c(g_2, \ldots, g_n)$.

In view of Theorem \ref{Dihedral d's and h's} 
and Corollary \ref{Connecting Map Dihedral Corollary}
the machinery in Theorem \ref{Positselski's Theorem} gives the following result, due to Positselski (Prop. 17 in \cite{Positselski}). 
\begin{theorem} \label{Connecting Map Dihedral Theorem} 
In the dihedral case we have the following 6-term exact sequence. 
\[
\begin{tikzpicture}[descr/.style={fill=white,inner sep=1.5pt}]
        \matrix (m) [
            matrix of math nodes,
            row sep=3em,
            column sep=3em,
            text height=1.5ex, text depth=0.25ex
        ]
        { H^n(E,\mu_d)\oplus H^n(\cG, \overline{M}_4) & H^n(\cG,\overline{M}_3) & H^n(\cG, \overline{M}_4)   \\ 
        H^{n+1}(F,\mu_d) & H^{n+1}(E,\mu_d) & H^{n+1}(F,\mu_d)\oplus H^{n+1}(\cG,\overline{M}_3) \\ 
        };

        \path[overlay, font=\scriptsize,>=latex]
        (m-1-1) edge [->]  node[descr,yshift=1.5ex] {$d_2 + h_3 $} (m-1-2) 
        (m-1-2) edge [->]  node[descr,yshift=1.5ex] {$d_3 $} (m-1-3) 
        (m-2-1) edge [->]  node[descr,yshift=1.5ex] {$d_1 $} (m-2-2) 
        (m-2-2) edge [->]  node[descr,yshift=1.5ex] {$h_1 \oplus d_2$} (m-2-3) 
        (m-1-3) edge[out=355,in=175,->] node[descr,yshift=0.3ex] {$\eta$} (m-2-1)
        
;
\end{tikzpicture} 
\] 
where $\eta(c) = \chi \smile' c$ and $d_1$ is scalar extension.
\end{theorem}

\begin{corollary} \label{Dihedral Cohomological Kernel Result}
$\eta$ induces an isomorphism 
\[
\frac{H^n(\cG,\overline{M}_4)}{d_3 H^n(\cG,\overline{M}_3)} \overset{\cong}{\longrightarrow} H^n(E/F).
\]
\end{corollary} 
In Theorem \ref{Dihedral Corestriction}, it will be shown that $H^n(\widetilde{F},\mu_d)$ maps onto $H^n(\cG,\overline{M}_4)$ in such a way that the image of the corestriction from $H^n(\widetilde{E},\mu_d)$, cor$_{\widetilde{E}/\widetilde{F}}$, maps onto the image of $d_3$ from $\overline{M}_3$. This will be used to characterize the cohomological kernel for the dihedral setup as follows: 
\[
\frac{H^n(\widetilde{F},\mu_d)}{\text{cor}_{\widetilde{E}/\widetilde{F}} H^n(\widetilde{E},\mu_d)} \cong H^n(E/F).
\]

\section{The Semi-Direct Case} \label{Semi-direct Section}
This section expands the ideas of the previous section.  
Technically, the previous two sections could be interpreted as applications of the results contained in this section, 
but for this paper we decided it would be prudent to present them separately illustrate the development of these tools.
\subsection{The 4 Term Exact Sequence with Homotopies} 
For this section we adopt the following notation:  $G$ is a semi-direct product of $\langle\tau\rangle$ by $\langle\sigma\rangle$.  
We assume the order of $\tau$ is $d$, the order of $\sigma$ is $s$ and we will assume that $s$ is even and divides $d-1$, making $d$ odd. 
We denote by $\theta:\{0,1,\ldots,d-1\}\rightarrow \{0,1,\ldots,d-1\}$ the permutation defined by $\sigma$-conjugation on $\langle\tau\rangle$ that makes
$\sigma\tau^i\sigma^{-1}=\tau^{\theta(i)}$ for every $i$. 
We define $\theta_j$ in a similar fashion to be defined by conjugation by $\sigma^j$, so that $\sigma^j\tau\sigma^{-j}=\tau^{\theta_j}$ for every $j$. In fact $\theta_j=\theta^j(1)$, where the latter is the $j$'th iterate of $\theta$, but the notation $\theta_j$ is less cumbersome. 
We assume that $\theta$ has order $s$, that is, conjugation by $\sigma$ on $\langle\tau\rangle$ does not have a smaller order than $s$. 
As $\tau$ has odd order, this means   $\sigma^{\frac{s}{2}}\tau^i\sigma^{-\frac{s}{2}}= \tau^{-i}$ for all $i$. 
From this, $\theta_{j+\frac{s}{2}}\equiv-\theta_j$ (mod $d$) and since $0<\theta_j<d$ we must have $\theta_j+\theta_{j+\frac{s}{2}}=d$. 
Here are the diagrams of the fields and Galois groups. 
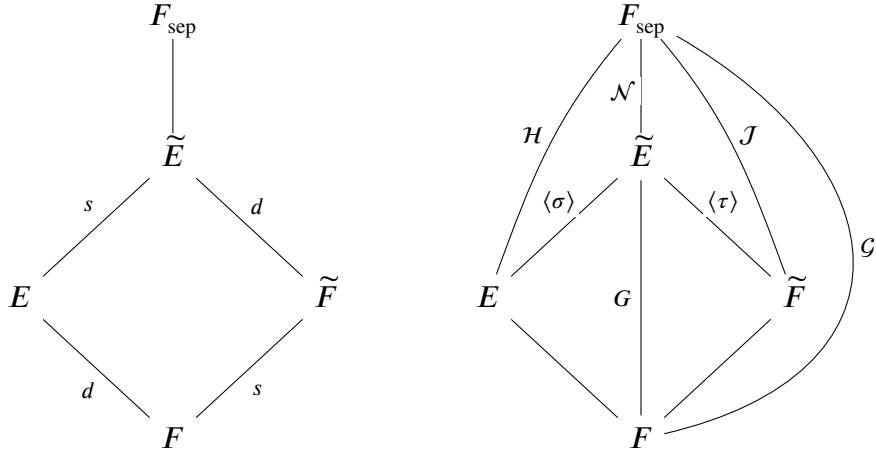
\[
\begin{tikzpicture}[descr/.style={fill=white,inner sep=1.5pt}]
        \matrix (m) [
            matrix of math nodes,
            row sep=3em,
            column sep=3em,
            text height=1.5ex, text depth=0.25ex
        ]
        {\ & F_{\text{sep}} & \ \\ 
         \ & \widetilde{E} & \                                        \\ 
            E & \ & \widetilde{F} \\ 
           \ & F & \                                        \\ \\ 
        };

        \path[overlay, font=\scriptsize,>=latex]
        (m-2-2) edge  node[descr,xshift=-1.75ex] {$\ $} (m-1-2) 
        (m-3-1) edge  node[descr,xshift=-3,yshift=1.5ex] {$s$}(m-2-2) 
        (m-3-3) edge node[descr,xshift=3,yshift=1.5ex] {$d$}(m-2-2) 
        (m-4-2) edge node[descr,xshift=-3,yshift=-1.5ex] {$d$} (m-3-1) 
        (m-4-2) edge node[descr,xshift=3,yshift=-1.5ex] {$s$} (m-3-3) 
;
\end{tikzpicture} \hspace{.5in}
\begin{tikzpicture}[descr/.style={fill=white,inner sep=1.5pt}]
        \matrix (m) [
            matrix of math nodes,
            row sep=3em,
            column sep=3em,
            text height=1.5ex, text depth=0.25ex
        ]
        {\ & F_{\text{sep}} & \ \\ 
         \ & \widetilde{E} & \                                        \\ 
            E & \ & \widetilde{F} \\ 
           \ & F & \                                        \\ \\ 
        };

        \path[overlay, font=\scriptsize,>=latex]
        (m-2-2) edge  node[descr,xshift=-1.25ex,yshift=-2] {$\cN$} (m-1-2) 
        (m-3-1) edge[out=70,in=230] node[descr,xshift=-1ex,yshift=1ex] {$\cH$} (m-1-2) 
        (m-3-3) edge[out=110,in=310] node[descr,xshift=1ex,yshift=1ex] {$\cJ$} (m-1-2) 
        (m-3-1) edge  node[descr,xshift=-2.1,yshift=1.75ex] {$\langle \sigma \rangle$} (m-2-2) 
        (m-3-3) edge  node[descr,xshift=2.1,yshift=1.75ex] {$\langle \tau \rangle$} (m-2-2) 
        (m-4-2) edge[out=13,in=-30,looseness = 1.7] node[descr,xshift=1ex,yshift=1ex] {$\cG$} (m-1-2) 
        (m-4-2) edge node[descr,xshift=-7,yshift=0ex] {$G$} (m-2-2) 
        (m-4-2) edge (m-3-1) 
        (m-4-2) edge (m-3-3) 
;
\end{tikzpicture} 
\] 

In order to generalize from the dihedral case a special element
$\cB_{d,s}=\cB\in\Z[G]$ is the essential tool. It is described next.

\noi \begin{definition} For $\tau$ and $\sigma$ as above,
\[
\cB_{d,s}=\cB=(1-\sigma^{\frac{s}{2}})\tau^{\frac{d+1}{2}}\sum_{j=0}^{\frac{s}{2}-1}\left(\sum_{i=0}^{\theta_j-1}\tau^i \right) \sigma^j.
\]
We set $T_{\sigma,s} = T_{\sigma}=
1+\sigma+\sigma^2+\cdots+\sigma^{s-1}$ and define
\[
C_{d,s,i}=C_i=\tau^iT_\sigma(1-\tau).
\]
\end{definition} 

Note that this $\cB$ is not quite the same it was in the previous section, even when restricting to  the dihedral case (making $s=2$). 
This $\cB$ is a generalization of what was used in the original calculations for the dihedral case. 
The dihedral $\cB$ was modified along with the maps to make the computations and proofs smoother. 
We begin with some basic properties of $\cB$. 

\noi \begin{lemma} \label{b lemma} Given the above assumptions and notation we have,\\
(i) $\sigma^{\frac{s}{2}}\cdot\cB=-\cB$.\\
(ii) $(1-\tau)\cB=C_{\frac{d+1}{2}}$.
\end{lemma}

\noi {\bf Proof}: Part (i)  is clear by the definition of $\cB$
since $\sigma^{\frac{s}{2}}(1-\sigma^{\frac{s}{2}})=-(1-\sigma^{\frac{s}{2}})$. \\ 
For (ii) We have $\tau\sigma^{\frac{s}{2}} = \sigma^{\frac{s}{2}}\tau^{-1}$, and hence 
\[
(1-\tau)(1-\sigma^{\frac{s}{2}}) = 1 - \tau + \sigma^{\frac{s}{2}}\tau^{-1} - \sigma^{\frac{s}{2}} = (1 + \sigma^{\frac{s}{2}}\tau^{-1})(1-\tau). 
\]
This allows us to make the substitution 
\[
(1-\tau)(1-\sigma^{\frac{s}{2}})\tau^{\frac{d+1}{2}} = \tau^{\frac{d+1}{2}}(1 + \sigma^{\frac{s}{2}})(1-\tau)
\]
Now, if we multiply the inner sum in the defintion of $\cB$ by $(1-\tau)$, we get 
\[
(1-\tau)\left(\sum_{i=0}^{\theta_j-1}\tau^i \right) = 1 - \tau^{\sigma_j}
\]
These facts allow us to establish $(ii)$.
\begin{eqnarray*}
(1-\tau) \cB&=&(1-\tau)(1-\sigma^{\frac{s}{2}})\tau^{\frac{d+1}{2}}\sum_{j=0}^{\frac{s}{2}-1}\left(\sum_{i=0}^{\theta_j-1}\tau^i \right) \sigma^j\\
&=&\tau^{\frac{d+1}{2}}(1 + \sigma^{\frac{s}{2}})(1-\tau)\sum_{j=0}^{\frac{s}{2}-1}\left(\sum_{i=0}^{\theta_j-1}\tau^i \right) \sigma^j\\
&\overset{(1)}{=}&\tau^{\frac{d+1}{2}}(1 + \sigma^{\frac{s}{2}})\sum_{j=0}^{\frac{s}{2}-1}\left( 1-\tau^{\theta_j} \right)\sigma^j\\
&\overset{(2)}{=}&\tau^{\frac{d+1}{2}}(1 + \sigma^{\frac{s}{2}})\sum_{j=0}^{\frac{s}{2}-1}\sigma^j\left( 1-\tau \right)\\
&=&\tau^{\frac{d+1}{2}}T_\sigma \left( 1-\tau \right)\\
&=&C_{\frac{d+1}{2}}
\end{eqnarray*}
with the equality labeled $(1)$ using the geometric series identity $(1-\tau)\sum_{i=0}^{\theta_j-1} \tau^i = 1 - \tau^{\theta_j}$ 
and the equality labeled $(2)$ following from 
$(1+\sigma^{\frac{s}{2}})\sum_{j=0}^{\frac{s}{2}-1} \sigma^j$ 
being the sum of every power of $\sigma$, i.e. $T_\sigma$. 
This shows directly that $C_{\frac{d+1}{2}}=(1-\tau)\cB$ giving (ii). \hfill $_\Box$ \\

We next define the modules we need.

\noi \begin{definition} Define all four modules to be submodules of $\Z[G]$, defined as follows. 
\begin{eqnarray*}
M_1&=&\Z[G]T_\sigma T_\tau \cong \Z \text{ as a trivial $\cG$-module} \\ 
M_2&=&\Z[G]T_\sigma = \sum_{i=0}^{d-1}{\Z}\tau^iT_\sigma \\ 
M_3&=&\Z[G]\cB = \sum_{j=0}^{\frac{s}{2}-1} \Z[\langle \tau \rangle] \sigma^j\cB \\ 
M_4&=&\Z[G]\cB T_\tau = \sum_{j=0}^{\frac{s}{2}-1}{\Z}\cdot\sigma^{j}\cB T_\tau
\end{eqnarray*}
\end{definition}

This next lemma gives key properties of the modules.

\begin{lemma} \label{M_3 Decomposition Semidirect}
\begin{eqnarray*}
\Z[G]\cB &=&\Z[G](1-\tau)\cB \oplus \Z[\langle \sigma \rangle] \cB \\
&=&M_3'\oplus M_{\cB}\\
&&\text{where} \ M_3'=\sum_{i=0}^{d-1}{\Z}\cdot\tau^iT_\sigma(1-\tau)  \ \text{ and } M_{\cB}=\sum_{j=0}^{\frac{s}{2}-1}{\Z}\cdot\sigma^{j}\cB_{d,s}\text{ and the direct }\\
&&  \text{sum between $M_3'$ and $M_\cB$} \text{is that of }\Z\text{-modules, not of }\cG\text{-modules}. \\
\end{eqnarray*}
\end{lemma}

\noi {\bf Proof}: For the (additive) direct summands of $M_3$, the first identification 
\[
\Z[G](1-\tau)\cB = \sum_{i=0}^{d-1}\Z \cdot \tau^i T_\sigma(1-\tau)
\]
follows from Lemma \ref{b lemma}$(i)$ while the second identification 
\[
\Z[\langle \sigma \rangle ] \cB = \sum_{j=0}^{\frac{s}{2}-1}\Z \cdot \sigma^j \cB
\]
follows from Lemma \ref{b lemma} (ii). \hfill $_\square$ \vspace{10pt}
 
We also note that  the ${\Z}$-ranks of $M_1$, $M_2$, $M_3$, and $M_4$ are, 
respectively $1$, $d$, $d-1+\frac{s}{2}$, and $\frac{s}{2}$ 
(although for the latter one has to check the linear independence of the $\sigma^j\cB$ from 
$M'_3$.) 
We next define the $d_i$ maps, which are similar to those from the dihedral case.

 \begin{definition}
The $d_i:M_i\rightarrow M_{i+1}$  are  as follows:

\smallskip

$d_1:M_1 \longrightarrow M_2$ is the inclusion $\Z[G]T_\sigma T_\tau \overset{\subseteq}{\longrightarrow} \Z[G]T_\sigma$. 

\smallskip 

If we view $M_1$ as $\Z$, then 
$d_1(n)=\sum_{i=0}^{d-1}n\tau^iT_\sigma$. 

\smallskip

$d_2:M_2 \rightarrow M_3$ is given by
$\cdot (1-\tau):\Z[G]T_\sigma \longrightarrow \Z[G]T_\sigma(1-\tau) \subseteq \Z[G]\cB$. 

\smallskip 

with ``$\subseteq$'' coming from the identity $\Z[G]T_\sigma(1-\tau) = \Z[G](1-\tau)\cB$ in Lemma \ref{b lemma}. 

\smallskip

$d_3:M_3\rightarrow M_4$ is given by $\cdot T_\tau:\Z[G]\cB \longrightarrow \Z[G]\cB T_\tau$. 
\end{definition} 

We note that each map  is a $\cG$-module homomorphism. 
The map $d_1$, which is an inclusion, can be thought of as the diagonal embedding if $M_1$ is viewed as $\Z$, with an image that has a trivial $\cG$-action. 
The map $d_2$ is right mulplitplication by
$(1-\tau)$ and hence is a $\cG$-map. 
By construction $M_3'$ is the image of $d_2$.
The map $d_3$ is right multiplication by $T_\tau$ and can be viewed as a trace map on $M_\cB$, the right summand of $M_3$. 
The trace is also trivial on $M_3'$ because $(\tau-1)T_\tau = 0$, so we need only consider $d_3$ applied to $M_{\cB}$. 
The homotopy maps are given next.

\noi \begin{definition} 
The $h_i:M_{i+1} \rightarrow M_i$  are  as follows:

\smallskip

 $h_1:\Z[G]T_\sigma \longrightarrow \Z[G]T_\sigma T_\tau :=
\cdot T_\tau$ is given by $xT_\sigma \mapsto xT_\sigma T_\tau$. 

\smallskip

$h_2:\Z[G]\cB \longrightarrow \Z[G]T_\sigma$ is given by $h_2(x \cB) = x \sum_{i=0}^{d-1}(\frac{d-1}{2} - i)\tau^i \tau^{\frac{d+1}{2}} T_\sigma$ for every $x \in \Z[G]$. 

\smallskip

$h_3:\Z[G] \cB T_\tau \overset{\subset}{\longrightarrow} \Z[G] \cB$ is the inclusion. 
\end{definition}

The next result verifies that the maps just defined satisfy the Positselski hypotheses.

\begin{theorem} \label{Semi-direct Checks}  
In the semi-direct case, given the above 
definitions we have the following.
\begin{enumerate} 
\item The $d_i$'s are exact 
\item $h_2$ is well-defined 
\item The prism condition is satisfied at all 4 modules. 
\end{enumerate}
\end{theorem} 

\noi {\bf Proof}: 
For part (1), exactness follows from extending the Hilbert 90 sequence discussed in the dihedral case, 
\[
0 \longrightarrow \Z[\langle \tau \rangle]\cdot T_\tau \overset{\subseteq}{\longrightarrow} \Z[\langle \tau \rangle] \overset{\cdot(1-\tau)}{\longrightarrow} \Z[\langle \tau \rangle](1-\tau) \longrightarrow 0.
\]
Here we replace $\Z[\langle \tau \rangle]$ with $\Z[G]T_\sigma$, which is isomorphic as a $\Z[\langle \tau \rangle]$-module to $\Z[\langle \tau \rangle]$. 
The second short exact sequence is immediate from the direct sum decomposition, 
although it should be stated that it is not split exact, since the direct sum is only that of $\Z$-modules, not $\cG$-modules. 

For part (2), as in the previous section we will show that the left annihilator of $\cB$ is in the left annihilator of $h_2(\cB)$. 
Let $x \in \Z[G]$ such that $x \cB = 0$. We will use the direct sum decomposition of $\Z[G]$ to express $x$ as follows: 
\[
x = (x_1(1-\tau) , x_2),
\]
where $x_1 \in \Z[G]$ and $x_2 \in \Z[\langle \sigma \rangle]$. 
Now we use the direct sum decomposition of $M_3 = \Z[G]\cB$: 
\[
x \cB = (x_1(1-\tau)\cB , x_2 \cB) = (0 , 0). \]
The fact that $x_1 \cB = 0$ and $x_2 \cB = 0$ will be used after we apply $h_2$.
But before applying $h_2$ to the elements of each direct summand we first make two observations, $i)$ and $ii)$, about $x_1$ and $x_2$ respectively that will be important in understanding where $h_2$ sends both of these elements. 
\begin{enumerate}[label=\roman*)] 
\item $x_1 \tau^{\frac{d+1}{2}} T_\sigma = k T_\tau$ for some $k \in \Z[\langle \sigma \rangle]$. \\ 
This follows because once $x_1 \tau^{\frac{d+1}{2}} T_\sigma$ is multiplied by $(1-\tau)$ we get 0, 
which is shown directly as follows, 
\[
x_1 \tau^{\frac{d+1}{2}} T_\sigma(1-\tau) = x_1 (1-\tau) \cB = 0 \]
with the first equality following from Lemma \ref{b lemma} $(ii)$. 
The left-annihilator of $(1-\tau)$ in $\Z[G]$ is $\Z[G] T_\tau$. 

\item $x_2$ is a left-multiple of $(1 + \sigma^{\frac{s}{2}})$. 
This follows from the direct sum decomposition
\[
\Z[\langle \sigma \rangle] \cB = \oplus_{j=0}^{\frac{s}{2}-1} \Z \sigma^j \cB.
\]
\end{enumerate} 
With the above observations, we are ready to compute, starting with $h_2(x_1(1-\tau)) \cB$: 
\[
\begin{split} 
h_2(x_1(1-\tau) \cB) &= x_1(1-\tau) \sum_{i=0}^{d-1} \left( \frac{d-1}{2} - i \right) \tau^i \tau^{\frac{d+1}{2}} T_\sigma = x_1 \tau^{\frac{d+1}{2}} \left( (1-\tau) \sum_{i=0}^{d-1} \left( \frac{d-1}{2} - i \right) \tau^i \right) T_\sigma \\ 
&= x_1 \tau^{\frac{d+1}{2}} \left( d - T_\tau \right) T_\sigma = x_1 \tau^{\frac{d+1}{2}} T_\sigma \left( d - T_\tau \right) \overset{i)}{=} kT_\tau (d - T_\tau) = k(dT_\tau - T_\tau^2) \\ 
&= k(dT_\tau - dT_\tau) = 0
\end{split} 
\]
Now we compute $h_2(x_2 \cB)$: 
\[
\begin{split} 
h_2(x_2 \cB) &\overset{ii)}{=} h_2(x_2'(1+\sigma^{\frac{s}{2}}) \cB) = x_2'(1+\sigma^{\frac{s}{2}}) \sum_{i=0}^{d-1} \left( \frac{d-1}{2} - i \right) \tau^i \tau^{\frac{d+1}{2}} T_\sigma \\ 
&= x_2' (d-1)T_\sigma + x_2'(1+\sigma^{\frac{s}{2}}) \sum_{i=0}^{d-1} \left( -i \right) \tau^i \tau^{\frac{d+1}{2}} T_\sigma \\ 
&= x_2' (d-1)T_\sigma + x_2' \sum_{i=0}^{d-1} \left( -i \right) \tau^i \tau^{\frac{d+1}{2}} T_\sigma + x_2'\sum_{i=0}^{d-1} \left( -i \right) \tau^{-i} \tau^{\frac{d-1}{2}} \sigma^{\frac{s}{2}}T_\sigma \\
&\overset{*}{=} x_2' (d-1)T_\sigma + x_2' (- (d-1)) T_\sigma = 0
\end{split} 
\]
For the second to last equality above (*), note that the coefficients in the two summands add to $-(d-1)$ for each power of $\tau$. 
With $x_1(1-\tau)$ and $h_2(x_2)$ both in left-annihilator of $h_2(\cB)$, 
we have shown that $h_2$ is well-defined. 

For part (3), the prism condition at $M_1$, which is $h_1   d_1 = \cdot d$, holds because $T_\tau \cdot T_\tau = d \cdot T_\tau$, as $h_1   d_1 = \cdot T_\tau$. The same goes for the prism condition at $M_4$, since $d_3   h_3 = \cdot T_\tau$ as well. 

Furthermore, the prism condition at $M_2$ holds from the following calculations: $(d_1   h_1) (T_\sigma) = T_\sigma \cdot T_\tau$, while 
\[
\begin{split} 
(h_2   d_2)(T_\sigma) &= h_2(T_\sigma(1-\tau)) = h_2((1-\tau)\tau^{\frac{d-1}{2}}\cB) =(1-\tau)\tau^{\frac{d-1}{2}}\sum_{i=0}^{d-1}\left(\frac{d-1}{2}-i\right)\tau^{\frac{d+1}{2}+i}T_\sigma \\ 
&= (1-\tau)\sum_{i=0}^{d-1}\left(-i\right)\tau^iT_\sigma = -T_\tau T_\sigma + d\cdot T_\sigma.  
\end{split} 
\]
Therefore 
\[
(d_1   h_1 + h_2   d_2)(T_\sigma) = T_\sigma T_\tau - T_\sigma T_\tau + d \cdot T_\sigma = d \cdot T_\sigma.
\]
The prism condition holds at $M_3$ from the following calculations: $(h_3   d_3)(\cB) = \cB \cdot T_\tau = T_\tau \cdot \cB$, while 
\[
\begin{split}
(d_2   h_2)(\cB) &= d_2 \left( \sum_{i=0}^{d-1}\left(\frac{d-1}{2}-i\right)\tau^{\frac{d+1}{2}+i}T_\sigma \right) = \left( \sum_{i=0}^{d-1}\left(\frac{d-1}{2}-i\right)\tau^{\frac{d+1}{2}+i}T_\sigma \right)(1-\tau) \\
&= \sum_{i=0}^{d-1}\left(\frac{d-1}{2}-i\right)\tau^{i}(1-\tau)\cB =  \sum_{i=0}^{d-1}-i\tau^{i}(1-\tau)\cB = \left(-T_\tau + d \right) \cB. 
\end{split}
\]
Therefore 
\[
(d_2   h_2 + h_3   d_3)(\cB) = T_\tau \cdot \cB + \left( -T_\tau + d \right)\cB = d \cdot \cB.
\]
And finally, $d_3h_3(\cB T_\tau) = d_3(\cB T_\tau) = \cB T_\tau T_\tau = d \cB T_\tau$. 
This concludes the proof of Theorem \ref{Semi-direct Checks}. \hfill $\square$ \vspace{10pt}

The only remaining part to check for the Positselski hypotheses is that the Bockstein homomorphisms are zero. 
$M_1$ and $M_2$ are $\Z$ (as a trivial $\cG$-module) and the induced module $\Ind_{\cH}^\cG(\Z)$ respectively, each are viewed modulo $d^2\Z$ to represent $\mu_{d^2}$ with a trivial $\cG$-action on cohomology with their respective Galois groups. 
Both of these modules were shown to have zero Bockstein maps in the general setup section. 
For $M_3$ and $M_4$, we will show that as $\cJ$-modules, $M_3/dM_3 \cong \Z[J]/d^2\Z[J] \cong \oplus_{j=0}^{\frac{s}{2}-1} \Ind_\cN^\cJ(\Z/d^2\Z)$ 
and that $M_4/dM_4 \cong \oplus_{j=0}^{\frac{s}{2}-1} \Z/d^2\Z$ as a trivial $\cJ$-module. 
The isomorphisms are defined by 
$\cB \mapsto 1 \in \Z[\cJ]$ and $\cB T_\tau \mapsto 1 \in \Z$ respectively.

\begin{lemma} \label{Semi-Direct Bockstein} As $\cJ$-modules,  
\begin{enumerate} 
\item $M_3 = \Z[G]\cB \cong \oplus_{j=0}^{\frac{s}{2}-1} 
\Ind_\cN^\cJ(\Z)$ 
\item $M_4 = \Z[G]\cB T_\tau \cong \oplus_{j=0}^{\frac{s}{2}-1} \Ind_\cJ^\cJ(\Z)$ 
\end{enumerate} 
\end{lemma} 
\noi {\bf Proof}: 
Both parts of this lemma follow from the fact that 
\[
\Z[G]\cB = \oplus_{j=0}^{\frac{s}{2}-1}\sigma^j \Z[\langle \tau \rangle]\cB
\]
followed by the arguments in Lemma \ref{Dihedral Zero Bockstein} applied to each direct summand. Note that in the dihedral case, $s=2$ and therefore $\sigma = \sigma^{\frac{s}{2}}$. So $\sigma^{\frac{s}{2}}$ must be used to apply the same arguments. \hfill $_\square$ \vspace{10pt}

Thus, the Bockstein maps have component maps for each direct summand that are identical to those in the dihedral case. And the dihedral Bockstein maps were shown to be zero. Therefore the Bockstein maps are zero on $M_3$ and $M_4$, completing the verification of the Positselski hypotheses for the 4-term exact sequence of $\cG$-modules.

\subsection{The Connecting Map for the Semi-Direct Case} 
Now we compute the connecting map $\eta$. Let $\overline{M_i} := M_i/dM_{i}$ for $i \in \{1,2,3,4\}$. 
In this section we will use the exact sequence of modules with homotopies defined in the previous section to describe connecting map $\eta: H^{n-1}(\cG, \overline{M_4}) \longrightarrow H^n(\cG, \overline{M_1})$. 

\begin{definition} Given the above notation we define the following.
\begin{enumerate} 
\item $\ell: \overline{M_4} \longrightarrow \overline{M_3}$, the $d_3$-lifting, where 
$\ell \left( \sum_{j=0}^{\frac{s}{2}-1} x_j \sigma^j \cB T_\tau \right) = \sum_{j=0}^{\frac{s}{2}-1} x_j \sigma^j \cB $ 
for every $x_j \in \Z/d\Z.$ 
\item $ \delta: C^{n-1}(\cG, \overline{M_3}) \longrightarrow C^n(\cG, \overline{M_3})$ the cochain map from the bar resolution.
\item $\widetilde{\eta}_{\ell} = \widetilde{\eta} : Z^{n-1}(\cG, \overline{M_4}) \longrightarrow Z^n(\cG, \overline{M_1})$, \ $\widetilde{\eta}(c) := -[d_1^{-1} h_2 \delta \ell(c)]$. 
\item $\eta: H^{n-1}(\cG, \overline{M_4}) \longrightarrow H^n(\cG, \overline{M_1})$, \ $\eta([c]) := [\widetilde{\eta}(c)]$. 
\end{enumerate} 
\end{definition} 

Our choice of lifting $\ell$ is a $\Z$-module homomorphism, though it is not a $\Z[G]$-module homomorphism. 
We will see in a moment that though $\ell$ is not a $\Z[G]$-module homomorphism, it is a module homomorphism for the subring $\Z[\langle \sigma \rangle]$ of $\Z[G]$ than properly contains $\Z = \Z\cdot 1_G$. 
Also, since $d_1$ is injective we let $d_1^{-1}(x)$ denote the unique preimage element for every $x \in \im(d_1)$.

This next lemma shows that $\ell$ preserves action by $\sigma$, a fact that was used at this point in the dihedral case in the previous section. 
However, in the dihedral case it was plain that $\sigma$ had the same action as $-1$ on $\cB$, since $\sigma$ was just $\sigma^{\frac{s}{2}}$,
so no lemma was necessary. 
Furthermore, the fact that $\sigma^{\frac{s}{2}}$ has this action of $-1$ in the semi-direct case will be central to the proof that follows. 

\begin{lemma} \label{Lifting Sigma Semidirect}
$\ell$ is a $\Z[\langle \sigma \rangle]$-module homomorphism. 
\end{lemma} 

\noi
{\bf Proof}: 
We will show that $\ell$ preserves action by $\sigma$. 
The last summand will be the only term of significance that needs to be addressed with respect to how it passes through the lifting after action by $\sigma$. 
\[
\begin{split} 
\ell\left( \sigma \cdot \sum_{j=0}^{\frac{s}{2} - 1} x_j \sigma^j \cB T_\tau \right) &= \ell\left( \sum_{j=0}^{\frac{s}{2} - 1} x_j \sigma^{j+1} \cB T_\tau \right) \\
&= \ell \left( x_{\frac{s}{2} - 1} \sigma^{\frac{s}{2}} \cB T_\tau + \sum_{j=0}^{\frac{s}{2} - 2} x_j \sigma^{j+1} \cB T_\tau \right) = \ell \left(  -x_{\frac{s}{2} - 1} \cB T_\tau + \sum_{j=0}^{\frac{s}{2} - 2} x_j \sigma^{j+1} \cB T_\tau \right) \\
&= - x_{\frac{s}{2} - 1} \cB + \sum_{j=0}^{\frac{s}{2} - 2} x_j \sigma^{j+1} \cB 
= x_{\frac{s}{2} - 1} \sigma^{\frac{s}{2}} \cB + \sum_{j=0}^{\frac{s}{2} - 2} x_j \sigma^{j+1} \cB \\ 
&= \sum_{j=0}^{\frac{s}{2} - 1} x_j \sigma^{j+1} \cB = \sigma \cdot \sum_{j=0}^{\frac{s}{2} - 1} x_j \sigma^j \cB = \sigma \cdot \ell \left( \sum_{j=0}^{\frac{s}{2} - 1} x_j \sigma^j \cB T_\tau \right).  
\end{split}
\]
This proves the lemma.
\hfill $_\square$ \vspace{10pt}

Now we are ready to compute the connecting map for the semi-direct case. 
The next lemma will allow us to find $\widetilde{\eta}(c)(g_1, \ldots, g_n)$ 
by first expressing $\delta(\ell(c))(g_1, \ldots, g_n)$ as 
$z\cdot \ell(c(g_2, \ldots, g_n))$ for an appropriate $z \in \Z[G]$ based on $g_1$, 
lifting this through $d_2$, 
and finally applying $-h_1$. 

\begin{lemma} \label{Connecting Map Semi-direct Lemma} 
Let $c \in Z^{n-1}(\cG, \overline{M_4})$ be a cocycle, $g_1, \ldots, g_n \in \cG$, and let $c' \in \Z/d\Z$ such that $c(g_2, \ldots, g_n) = c' \cdot \cB T_\tau$. 
Let $\sigma^i \tau^j$, an element of $\cG/\cN$, be the coset of $g_1$. 
Then \\ 
\begin{enumerate} 
\item $\delta(\ell(c))(g_1, \ldots, g_n) = g_1 \cdot \ell(c(g_2, \ldots, g_n)) - \ell(g_1 \cdot c(g_2, \ldots, g_n))$ \\ 
$= \sigma^j \left( \tau^i - 1 \right) \cdot \ell(c(g_2, \ldots, g_n)) $.  
\item For any $x \in C^n(\cG, \overline{M_2})$ such that $\delta(\ell(c)) = d_2(x)$, 
$\widetilde{\eta}(c) = -h_1(x).$ 
\end{enumerate} 
\end{lemma} 

\noi {\bf Proof}: 
 For (1), the first equality follows from the fact that $\ell$ is a $\Z$-module homomorphism, and hence the proof is identical to the analogous proof in the previous section. The second equality of (1) follows from $\ell$ being a $\Z[\langle \sigma \rangle]$-module homomorphism and $\tau$ acting trivially on $\overline{M_4}$. Note that $\tau$ acts trivially on $\overline{M_4}$ because $T_\tau$ is in the center of $\Z[G]$ because $\langle \tau \rangle$ is normal in $G$, and $\tau T_\tau = T_\tau$. We use these two reasons in tandem to show the second equality: 
\[
\begin{split} 
g_1 \cdot \ell(c(g_2, \ldots, g_n)) - \ell(g_1 \cdot c(g_2, \ldots, g_n)) &= \sigma^j \tau^i \ell(c(g_2, \ldots, g_n)) - \ell(\sigma^j \tau^i \cdot c(g_2, \ldots, g_n) \\ 
&= \sigma^j \tau^i \cdot \ell(c') - \sigma^j \cdot \ell(\tau^i c') \\ 
&= \sigma^j \tau^i \cdot \ell(c') - \sigma^j \cdot \ell(c') \\ 
&= \sigma^j(\tau^i - 1) \cdot \ell(c'). 
\end{split} 
\]
For (2), showing that $\widetilde{\eta}(c) = h_1(x)$ relies only on the prism condition, and these details are identical to those in part (2) of
Lemma \ref{Connecting Map Dihedral Lemma}, the analog to this lemma in the dihedral case. They will therefore be omitted here to prevent repetition. 
This concludes the proof of Lemma \ref{Connecting Map Semi-direct Lemma}. \hfill $_\square$ \vspace{10pt}

As a corollary we obtain a characterization of the connecting map. As in the case for the dihedral extensions, this connecting map behaves
like a cup product.

\begin{corollary} \label{Connecting Map Semi-direct Corollary}
Let $g_1, g_2, \ldots, g_n, c, c',i,j$ be defined as in Lemma \ref{Connecting Map Semi-direct Lemma}, and let $c_m \in \Z/d\Z$ such that $c(g_2, \ldots, g_n) = \sum_{m=0}^{\frac{s}{2}-1} c_m \sigma^m \cB T_\tau$. 
Then 
\[
\widetilde{\eta}(c)(g_1, \ldots, g_n) = -i\sum_{m=0}^{\frac{s}{2}-1} c_m (\theta_{s-m}-1) T_\sigma T_\tau
\]
\end{corollary} 

\noi {\bf Proof}: From part (1) of Lemma \ref{Connecting Map Semi-direct Lemma}, we know that 
\[
\delta(\ell(c))(g_1, \ldots, g_n) = \sigma^j \left( -\sum_{k=0}^{i-1} \tau^k \right) \cdot (1-\tau) c' \cB.
\]
Now we will find an $x \in \overline{M_2}$ with the above $d_2$-image so we can apply $-h_1$ to $x$. By part (2) of Lemma \ref{Connecting Map Semi-direct Lemma}, this yields $\widetilde{\eta}(c) = -h_1(x)$. 

We show now that for $x=\sigma^j \left( -\sum_{k=0}^{i-1} t^k \right) \sum_{m=0}^{\frac{s}{2}-1} c_m \sigma^m \left( \sum_{r=0}^{\theta_m^{-1}-1} \tau^r \right) \tau^{\frac{d+1}{2}} T_\sigma$, 
we have the desired equality: 
\[d_2(x) = \delta(\ell(c))(g_1, \ldots, g_n).
\]
We start with the right hand side.
\[
\begin{split} 
\delta(\ell(c))(g_1, \ldots, g_n) &= \sigma^j \left( \tau^i - 1 \right) \cdot \ell(c(g_2, \ldots, g_n)) \\ 
&= \sigma^j \left( -\sum_{k=0}^{i-1} t^k \right) (1-\tau) \cdot \ell(c') \\ 
&= \sigma^j \left( -\sum_{k=0}^{i-1} t^k \right) (1-\tau) \cdot \sum_{m=0}^{\frac{s}{2}-2} c_m \sigma^m \cB \\ 
&= \sigma^j \left( -\sum_{k=0}^{i-1} t^k \right) \sum_{m=0}^{\frac{s}{2}-2} c_m \sigma^m  (1-\tau^{\theta_m^{-1}})\cB \\ 
&= \sigma^j \left( -\sum_{k=0}^{i-1} t^k \right) \sum_{m=0}^{\frac{s}{2}-2} c_m \sigma^m  (1-\tau^{\theta_m^{-1}}) \cB \\ 
&= \sigma^j \left( -\sum_{k=0}^{i-1} t^k \right) \sum_{m=0}^{\frac{s}{2}-1} c_m \sigma^m \left( \sum_{r=0}^{\theta_m^{-1}-1} \tau^r \right) (1-\tau) \cB \\ 
&= \sigma^j \left( -\sum_{k=0}^{i-1} t^k \right) \sum_{m=0}^{\frac{s}{2}-1} c_m \sigma^m \left( \sum_{r=0}^{\theta_m^{-1}-1} \tau^r \right) \tau^{\frac{d+1}{2}} T_\sigma (1-\tau) \\ 
&= d_2 \left( \sigma^j \left( -\sum_{k=0}^{i-1} t^k \right) \sum_{m=0}^{\frac{s}{2}-1} c_m \sigma^m \left( \sum_{r=0}^{\theta_m^{-1}-1} \tau^r \right) \tau^{\frac{d+1}{2}} T_\sigma \right). \\ \\ 
\end{split} 
\]
With $x$ having the desired property, we can now apply $h_1$ to $x$.
\[
\begin{split} 
h_1(x) &= h_1 \left( \sigma^j \left( -\sum_{k=0}^{i-1} t^k \right) \sum_{m=0}^{\frac{s}{2}-1} c_m \sigma^m \left( \sum_{r=0}^{\theta_m^{-1}-1} \tau^r \right) \tau^{\frac{d+1}{2}} T_\sigma \right) \\ 
&= \sigma^j \left( -\sum_{k=0}^{i-1} t^k \right) \sum_{m=0}^{\frac{s}{2}-1} c_m \sigma^m \left( \sum_{r=0}^{\theta_{s-m}-1} \tau^r \right) \tau^{\frac{d+1}{2}} T_\sigma T_\tau \\ 
&= -i \sum_{m=0}^{\frac{s}{2}-1} c_m \left( \theta_{s-m} - 1 \right) T_\sigma T_\tau \\ 
\end{split} 
\]
This concludes the proof of Corollary \ref{Connecting Map Semi-direct Corollary}. \hfill $_\square$ \vspace{10pt}

In view of Theorem \ref{Semi-direct Checks} 
and Corollary \ref{Connecting Map Semi-direct Corollary},
the machinery in Theorem \ref{Positselski's Theorem} gives the following result.

\begin{theorem} \label{Connecting Map Semi-direct Theorem} 
In the semi-direct case we have the following 6-term exact sequence  
\[
\begin{tikzpicture}[descr/.style={fill=white,inner sep=1.5pt}]
        \matrix (m) [
            matrix of math nodes,
            row sep=3em,
            column sep=3em,
            text height=1.5ex, text depth=0.25ex
        ]
        { H^n(E,\mu_d)\oplus H^n(\cG, \overline{M_4}) & H^{n+1}(E,\mu_d) & H^n(\cG, \overline{M_4})   \\ 
        H^{n+1}(F,\mu_d) & H^{n+1}(E,\mu_d) & H^n(F,\mu_d)\oplus H^n(E',\mu_d) \\ 
        };

        \path[overlay, font=\scriptsize,>=latex]
        (m-1-1) edge [->]  node[descr,yshift=1.5ex] {$d_2 + h_3 $} (m-1-2) 
        (m-1-2) edge [->]  node[descr,yshift=1.5ex] {$d_3 $} (m-1-3) 
        (m-2-1) edge [->]  node[descr,yshift=1.5ex] {$d_1 $} (m-2-2) 
        (m-2-2) edge [->]  node[descr,yshift=1.5ex] {$h_1 \oplus d_2$} (m-2-3) 
        (m-1-3) edge[out=355,in=175,->] node[descr,yshift=0.3ex] {$\eta$} (m-2-1)

;
\end{tikzpicture} 
\] 
where $\eta$ is as described in Corollary \ref{Connecting Map Semi-direct Corollary}
and $d_1$ is scalar extension.

\end{theorem}

\begin{corollary} \label{Semi-direct Cohomological Kernel Result}
$\eta$ induces an isomorphism 
\[
\frac{H^n(\cG,\overline{M_4})}{d_3 H^n(\cG,\overline{M_3})} \overset{\cong}{\longrightarrow} H^n(E/F).
\]
\end{corollary} 
At the end of this paper, it will be shown that $d_3$ can be viewed as the corestriction map from $\widetilde{E}$ and that 
\[
\frac{H^n(\cG,\overline{M_4})}{\text{cor}_{\widetilde{E}/\widetilde{F}} H^n(\widetilde{E},\mu_d)} \cong H^n(E/F).
\]




\subsection{Some Examples} 
We conclude this section with a computation of $\eta: H^n(\cG,\overline{M_4}) \longrightarrow H^n(F,\mu_d)$ for $n=0$ and $n=1$.
\begin{enumerate} 
\item \underline{$\eta$ for $n=0$} \\
In the semi-direct case, $s$ is even and the non-trivial action of $\sigma^{\frac{s}{2}}$ on $\overline{M_4}$ is multiplication by $(-1)$. Hence $\overline{M_4}$ has no fixed points with $d$ being odd, and therefore $\eta = 0$ because $H^0(\cG,\overline{M_4})$ is trivial.  

Note that this is not true for the cyclic case, where $s=1$ is not even. In the cyclic case, $H^0(\cG,\overline{M_4}) = H^0(\cG,\mu_d) = \mu_d$ and $\eta$ is the well-known cup product with the character $-\chi_\alpha$ defined by the extension $E/F = F(\alpha)/F$. 
\item \underline{$\eta$ for $n=1$} \\ 
Let $\chi \in Z^1(\cG, \overline{M_4})$ be a crossed homomorphism. 
Then the identity $\sigma \cdot \chi(\tau) = \theta_1 \chi(\tau)$ may be deduced as from the cocycle condition as follows.  
\[ 
\begin{split} 
\sigma \cdot \chi(\tau) &= \sigma \cdot \chi(\tau) + \chi(\sigma) - \chi(\sigma) = \chi(\sigma \tau) - \chi(\sigma) = \chi(\tau^{\theta_1}\sigma) - \chi(\sigma) \\ 
&= \tau^{\theta_1} \cdot \chi(\sigma) + \chi(\tau^{\theta_1}) - \chi(\sigma) = \tau^{\theta_1} \cdot \chi(\sigma) + (\tau^{\theta_1 - 1} + \cdots + 1) \cdot \chi(\tau) - \chi(\sigma) \\ 
&\overset{(*)}{=} 1 \cdot \chi(\sigma) + (\theta_1) \cdot \chi(\tau) - \chi(\sigma) = \theta_1 \chi(\tau) \\ 
\end{split}
\] 
where the equality $(*)$ comes from $\tau$ acting trivially on $\overline{M_4}$. So for a unique $t \in \Z/d\Z$, 
\[
\chi(\tau) = t \sum_{m=0}^{\frac{s}{2}-1} \theta_1^{-m} \sigma^m \cB T_\tau
\]
where $t$ may be thought of as the $\sigma^0 \cB T_\tau$-coefficient of $\chi(\tau)$. 

In this case $\chi(\sigma)$ can be expressed as follows:
Let $s:\Z \longrightarrow \Z/d\Z$ such that $s_{m+\frac{s}{2}} = -s_m$ for every $m \in \Z$ so that 
\[
\chi(\sigma) = \sum_{m=0}^{\frac{s}{2}-1} s_m \sigma^m \cB T_\tau
\]
and more generally, 
\[
chi(\sigma) = \sum s_{\gamma_m} \sigma^{\gamma_m} \cB T_\tau
\]
for any choice of representatives $\gamma_1, \ldots, \gamma_{\frac{s}{2}-1} \in \Z$ of the cosets $\Z/(\frac{s}{2})\Z$. 
Suppose $g_1 \cN = \sigma^j \tau^i$ and $g_2 \cN = \tau^k \sigma^\ell$. Then 
\[
\chi(g_1,g_2) = k \cdot \chi(\tau) + \left( \sum_{\beta=0}^{\ell-1} \sigma^\beta \right)\cdot \chi(\sigma) = \sum_{m=0}^{\frac{s}{2}-1} \left( k t \theta_{-1}^m + \sum_{\beta=0}^{\ell-1} s_{m-\beta} \right) \sigma^m
\]
and hence 
\[
c_m = k t \theta_{-1}^m + \sum_{\beta=0}^{\ell-1}s_{m-\beta}.
\]
The connecting map formula in Corollary \ref{Connecting Map Semi-direct Corollary} may be applied to yield 
\[
\widetilde{\eta}(g_1,g_2) = -i \sum_{m=0}^{\frac{s}{2}-1} \left( \theta_1^{-m} - 1 \right) c_m = -i \sum_{m=0}^{\frac{s}{2}-1} \left( \theta_1^{-m} - 1 \right) \left( k t \theta_1^{-m} + \sum_{\beta=0}^{\ell-1}s_{m-\beta} \right).
\]

\end{enumerate} 


\section{Interpreting the Sequences}

\bs 

In this section we record some consequences of the sequences 
in Theorems \ref{Connecting Map Cyclic Theorem}, 
\ref{Connecting Map Dihedral Theorem}, and 
\ref{Connecting Map Semi-direct Theorem}.  As noted in the introduction, when 
$E/F$ is cyclic  one has the classical description of the relative
Brauer group,
\[
\frac{F^*}{N_{E/F}(E^*)}\cong
{\rm ker} ({\rm Br}_dF
\rightarrow {\rm Br}_dE).
\]
As an immediate consequence of Theorem 
\ref{Connecting Map Cyclic Theorem},  this classical
result generalizes as follows.

\begin{theorem}\label{Cyclic Application Theorem} If $E/F$ is cyclic
with $\mu_d\subset F$ we have
a description of the cohomological kernel for all $n\geq0$,
\[
\frac{{H}^n(F,\mu_d)}{{\rm cor}_{E/F}H^n(E,\mu_d)}\cong
{\rm ker} (H^{n+1}(F,\mu_d)
\rightarrow H^{n+1}(E,\mu_d)).
\]
\end{theorem}

The goal in this section is to find what generalizations 
of Theorem \ref{Cyclic Application Theorem} are possible in the
dihedral and semi-direct product cases considered in
the previous two sections.
We continue to assume the notation of the last section.
To better understand $M_3$ and $M_4$ we need to introduce more 
induced modules and subgroups.

\begin{definition}\label{MoreModules}
We set $T_2:=J+\sigma^{\frac{s}{2}}J\in {\Ind}^\cG_\cJ{\Z}$
and then set $\cT_2:=\sum_{j=0}^{\frac{s}{2}-1}{\Z}\cdot \sigma^jT_2\subset{\rm Ind}^\cG_\cJ{\Z}$. 
We note that $\cT_2$ is a $\cG$-submodule of ${\rm Ind}^\cG_\cJ{\Z}$ with $\sigma^{\frac{s}{2}}T_2=T_2$. We set
$\cJ'$ to be the unique group in $\cG$ of index $\frac{s}{2}$ containing $\cJ$ and denote by
$J'$ the corresponding subgroup of $G$. This means that 
$J'=\langle\tau,\sigma^{\frac{s}{2}}\rangle$.  
\end{definition}

\bs

By Galois theory, since 
$J'$ has index $\frac{s}{2}$ in $G$, $\cJ'={\rm Gal}(F_{sep}/F')$ where $[F':F]=\frac{s}{2}$ and $F\subset F'\subset \tilde{F}=F(\alpha)$ (since $\tilde{F}/F$ is cyclic, $F'$ is the unique such intermediate extension.)
We 
have the following.

\bs

\begin{lemma}\label{M4Structure} Given the above definitions.\\ (i) As $\cG$-modules,
$\cT_2\cong {\rm Ind}^\cG_{\cJ'}{\Z}$.\\
(ii)  As $\cG$-modules, $M_4:=M_3/M_3'\cong {\rm Ind}^\cG_\cJ{\Z}/\cT_2$.\\
In particular we have an exact sequence of $\cG$-modules.
\[
0\rightarrow{\rm Ind}^\cG_{\cJ'}{\Z}
\rightarrow{\rm Ind}^\cG_{\cJ}{\Z}
\rightarrow M_4\rightarrow 0.
\]
\end{lemma}

\noi {\bf Proof.} We note that ${\rm Ind}^\cG_{\cJ'}{\Z}\cong
\oplus_{j=0}^{\frac{s}{2}-1}{\Z}\cdot\sigma^j J'$ with $\sigma$
acting cyclicly and $\tau$ acting trivially on the summands.  
This is exactly how $G$ acts on the summands of $\cT_2$
and the map $\sigma^jT_2\mapsto \sigma^j\cJ'$, giving
 the isomorphism required
for (i).

For (ii), by definition $M_4$ is the free ${\Z}$-module with basis
$\ocB$, $\sigma\ocB$,\ldots,$\sigma^{\frac{s}{2}-1}\ocB$, with
trivial $\tau$-action and with $\sigma$ acting cyclically every summand except the last, 
which multiplies by -1 as it goes to the first summand: $\sigma \cdot \sigma^{\frac{s}{2}-1}\ocB = \sigma^{\frac{s}{2}}\ocB=-\ocB$. From this it follows
that the map $\sigma^jJ\mapsto\sigma^j\ocB$ defines a $G$-map
${\rm Ind}^\cG_\cJ{\Z}\rightarrow M_4$ with
kernel $\cT_2$ (the latter as $J+\sigma^{\frac{s}{2}}J\mapsto
\ocB+\sigma^{\frac{s}{2}}\ocB=\ocB-\ocB=0\in M_4$.)
This concludes the proof of Lemma \ref{M4Structure}. \hfill $\Box$

\bs

\noi {\bf Remarks.} (i) The ${\Z}$-ranks of the modules 
${\rm Ind}^\cG_{\cJ'}{\Z}, {\rm Ind}^\cG_{\cJ'}{\Z},
M_4$ are, respectively, $\frac{s}{2}$, $s$, and $\frac{s}{2}$.\\
(ii) Of course, all of these modules can be taken (mod $d$) and the same
results apply.

\bs

The cohomology of $M_4$ can be interpreted using the 
sequence of Lemma \ref{M4Structure}.  By definition we have $F\subseteq F'\subseteq \tF$ where
$[\tF:F]=s$ and $[F':F]=\frac{s}{2}$. Computing cohomology in $\mu_d$ and using
the fact that $(d,s)=1$  we know that $H^n(F',\mu_d)\rightarrow H^n(F,\mu_d)$
must be injective for all $n$.  In particular the long exact sequence in cohomology gives
exact sequences
\[
0\rightarrow H^n(F',\mu_d)\rightarrow H^n(\tF,\mu_d)\rightarrow
H^n(\cG,M_4)\rightarrow0.
\]
This means if we let
$\overline{H}^n(F,\mu_d):={\rm cok}(H^n(F',\mu_d)\rightarrow H^n(\tF,\mu_d))$
then
the Positselski connecting map in Theorem \ref{Connecting Map Semi-direct Theorem}
gives a map
$\overline{\eta}:\overline{H}^n(\tF,\mu_d){\rightarrow} H^{n+1}(F,\mu_d)$ that computes the cohomological kernels as noted next.

\bs

\begin{theorem}\label{First Interpretation Theorem}  In the above notation with $E/F$ being a semi-direct extension, we have
an exact sequence,
\[
\overline{H}^n(\tF,\mu_d)\stackrel{\overline{\eta}}
{\rightarrow}H^{n+1}(F,\mu_d)
\rightarrow H^{n+1}(E,\mu_d).
\]
\end{theorem}

\bs

To understand the kernel of $\overline{\eta}$ one needs to further understand the 
cohomology of $M_3$ and how it maps into the cohomology of $M_4$. For this purpose, if $\pi^*:H^n(\tF,\mu_d)\rightarrow H^n(\cG,M_4)$ 
is the induced map, we shall denote by
\[
N_3^n(E/F):=\pi^{*-1}({\rm im}(H^n(\cG,M_3)\rightarrow 
H^n(\cG,M_4)))\subseteq H^n(\tF,\mu_d)
\]
and then Theorem \ref{Connecting Map Semi-direct Theorem} gives the following result.

\begin{theorem}\label{Second Interpretation Theorem}  In the above notation with $E/F$ being a semi-direct extension, we have the following
characterization of the cohomological kernel $H^{n+1}(E/F,\mu_d)$,
\[
\frac{{H}^n(\tF,\mu_d)}{i_{\tF/F'}H^n(F',\mu_d)+N_3^n(E/F)}\cong
{\rm ker} (H^{n+1}(F,\mu_d)
\rightarrow H^{n+1}(E,\mu_d)).
\]
\end{theorem}

\bs

When interpreted loosely, this result can be understood as the analogue
of Theorem \ref{Cyclic Application Theorem} in the 
more general case case. (When $E/F$ is cyclic of degree $d$
and  $s=1$ we would have $\tF=F$ and $N_3^n(E/F)
=
{\rm cor}_{E/F}H^n(E,\mu_d)$. Also, the subfield $F'$ 
doesn't exist in  the cyclic case.)

\bs

Next we turn to $M_3\subset {\Z}[G]$.  We set $\cH':=\langle\cJ,
\sigma\tau\rangle={\Gal}(F_{sep}/E')$ where $E'$ is discussed
above. We know by Lemma \ref{Semi-Direct Bockstein}
the modules $M_3={\Z}
[G]\cdot\cB$, $M_3'={\Z}
[G]\cdot(1-\tau)\cB$ and $M_4={\Z}
[G]\cdot T_{\tau}\cB$.
By Lemma \ref{M4Structure} we have the following.

\begin{lemma}\label{OldTheorem12} The following diagram of $G$-modules and
$G$-maps is commutative
with exact rows and columns.
\[
\begin{tikzpicture}[descr/.style={fill=white,inner sep=1.5pt}]
        \matrix (m) [
            matrix of math nodes,
            row sep=2em,
            column sep=2em,
            text height=1.5ex, text depth=0.25ex
        ]
        { \ & 0 & 0 & 0 & \ \\
        0 & \cK' & \cK & \Z[G] T_\tau(1 + \sigma^{\frac{s}{2}}) & 0  \\ 
        0 & \Z[G](1-\tau) & \Z[G] & \Z[G] T_\tau & 0 \\ 
        0 & M_3' & M_3 & M_4 & 0 \\
        \ & 0 & 0 & 0 & \ \\
        };

        \path[overlay, font=\scriptsize,>=latex]
        (m-2-2) edge [->]  
        (m-2-3) 
        (m-2-3) edge [->]  node[descr,yshift=1.5ex, xshift=3ex] {$\cdot T_\tau $} (m-2-4) 
        (m-3-2) edge [->]  node[descr,yshift=1.5ex] {$\ $} (m-3-3) 
        (m-3-3) edge [->]  node[descr,yshift=1.5ex, xshift=.3ex] {$\cdot T_\tau$} (m-3-4) 
        (m-4-2) edge [->]  node[descr,yshift=1.5ex, xshift=-.2ex] {$\ $} (m-4-3) 
        (m-4-3) edge [->]  node[descr,yshift=1.5ex] {$\cdot T_\tau$} (m-4-4) 

        (m-2-2) edge [->]  node[descr,xshift=1ex] {$\ $} (m-3-2) 
        (m-3-2) edge [->]  node[descr,xshift=1.5ex] {$\cdot \cB $} (m-4-2) 
        (m-2-3) edge [->]  node[descr,xshift=1ex] {$\ $} (m-3-3) 
        (m-3-3) edge [->]  node[descr,xshift=1.5ex] {$\cdot \cB $} (m-4-3) 
        (m-2-4) edge [->]  node[descr,xshift=1ex] {$\ $} (m-3-4) 
        (m-3-4) edge [->]  node[descr,xshift=1.5ex] {$\cdot \cB $} (m-4-4) 

        (m-1-2) edge [->] (m-2-2) 
        (m-1-3) edge [->] (m-2-3) 
        (m-1-4) edge [->] (m-2-4) 
        (m-4-2) edge [->] (m-5-2) 
        (m-4-3) edge [->] (m-5-3) 
        (m-4-4) edge [->] (m-5-4) 
        (m-2-1) edge [->] (m-2-2) 
        (m-3-1) edge [->] (m-3-2) 
        (m-4-1) edge [->] (m-4-2) 
        (m-2-4) edge [->] (m-2-5) 
        (m-3-4) edge [->] (m-3-5) 
        (m-4-4) edge [->] (m-4-5) 
;
\end{tikzpicture} 
\] \\ 

The right column is that of Lemma \ref{M4Structure}, where
${\rm Ind}^\cG_{\cJ'}{\Z}={\Z}[G] T_{\tau}(1+\sigma^{\frac{s}{2}})$, and ${\rm Ind}^\cG_{\cJ}{\Z}={\Z}[G] T_{\tau}$.
Here $\cK:={\rm ker}(\cdot\cB:{\Z}[G]\rightarrow M_3)$ and
$\cK':={\rm ker}(\cdot\cB:{\Z}[G]\cdot(1-\tau)\rightarrow M'_3)$.
Moreover, the ${\Z}$-ranks of $\cK'$, ${\Z}[G]\cdot(1-\tau)$, 
$M'_3$ are $(s-1)(d-1)$, $s(d-1)$, $(d-1)$, resp., 
the ${\Z}$-ranks of $\cK$, ${\Z}[G]$, 
$M_3$ are $sd-\frac{s}{2}-(d-1)$, $sd$, $\frac{s}{2}+(d-1)$, resp., and the ${\Z}$-ranks of $
\Z[G] \cdot T_\tau(1 + \sigma^{\frac{s}{2}})$, ${\Z}[G]\cdot T_{\tau}$, 
$M_4$ are $\frac{s}{2}$, $s$, $\frac{s}{2}$, resp.
\end{lemma}

\bs

\noi {\bf Proof.} For commutativity, as the first set of downarrows are
inclusions as are the first rightarrows, the only question is the lower right
square. Since $T_{\tau}$ is central in ${\Z}[G]$, $\cB\cdot T_{\tau}=
T_{\tau}\cdot \cB$ and so the lower right square commutes. 

For exactness, the three downward arrows given by $\cdot \cB$ are surjective from the  decomposition of $M_3$ in Lemma \ref{M_3 Decomposition Semidirect}. The first two columns are then exact from the choices of $\cK$ and $\cK'$ as kernels of their respective $\cdot \cB$ maps. The last column is exact from Lemma \ref{M4Structure}, as mentioned above in this Theorem statement. 

The bottom row is second short exact sequence in the 4-term Positselski sequence in the semi-direct case, and above it is the Hilbert-90 short exact sequence from which it was constructed. 
With all other rows and columns being exact, it then follows 
by the usual diagram chase that that the top row is also exact. 





Finally the ${\Z}$-ranks of the bottom row are given in the remark
following Lemma \ref{M4Structure}.  The ${\Z}$-ranks of the two right columns are
clear by previous work, so the ranks of $\cK$ and $\cK'$ follow
by arithmetic. \hfill $\Box$

\bs

As an application of Theorem \ref{OldTheorem12} we can characterize $N_3^n(E/F)$
via the corestriction.

\begin{theorem}\label{OldTheorem13} In the above notation with $E/F$ being a semi-direct extension, we have 
the following
exact sequence calculating the cohomological kernel $H^{n+1}(E/F,\mu_d)$,
\[
\frac{{H}^n(\tF,\mu_d)}{{\rm cor}_{\tE/\tF}H^n(\tE,\mu_d)}\rightarrow H^{n+1}(F,\mu_d)
\rightarrow H^{n+1}(E,\mu_d).
\]
The first map is injective provided
$H^{n+1}(\cG,\cK)\rightarrow
H^{n+1}(\tE,\mu_d)$ is injective.
\end{theorem}

\noi {\bf Proof.} Consider the following diagram. 
\[\xymatrix{
 H^n(\tE,\mu_d)\ar[d]\ar[r]^{{\rm cor}}&
 H^n(\tF,\mu_d)\ar[d]\\
 H^n(\cG,{ M}_3)\ar[d]\ar[r]&
 H^n(\cG,{ M}_4)\ar[d]^{^0}\ar[r]&
 H^{n+1}(F,\mu_d)\ar[r]&H^{n+1}(E,\mu_d)\\
H^{n+1}(\cG,\cK)\ar[d]\ar[r]&
 H^{n+1}(F_1,\mu_d)\ar[d]\\
 H^{n+1}(\tE,\mu_d)\ar[r]^{{\rm cor}}&
 H^{n+1}(\tF,\mu_d)
}
\]
The middle row is exact by Theorem \ref{Connecting Map Semi-direct Theorem}.
The two columns are exact by the long exact sequence of cohomology applied to the middle and right columns of 
the diagram in Theorem \ref{OldTheorem12}. The diagram commutes since all 
maps are those induced by the diagram in 
Theorem \ref{OldTheorem12}. The map $H^{n+1}(F_1,\mu_d)\rightarrow
H^{n+1}(\tF,\mu_d)$ is injective since $[\tF:F_1]$ is prime 
to $d$.  Therefore the map $H^{n}(\tF,\mu_d)\rightarrow
H^n(\cG,\cM_4)$ is surjective.
This gives a surjective map
${H}^n(\tF,\mu_d)\rightarrow {\rm ker}(H^{n+1}(F,\mu_d)
\rightarrow H^{n+1}(E,\mu_d))).$
The exactness of the sequence follows by
noting the diagram shows ${\rm cor}_{\tE/\tF}
(H^n(\tE,\mu_d))$ has trivial image in 
$H^{n+1}(F,\mu_d)$. For the second statement, if
$H^{n+1}(\cG,\cK)\rightarrow
H^{n+1}(\tE,\mu_d)$ is injective then 
$H^n(\tE,\mu_d)\rightarrow  H^n(\cG,M_3)$ is 
surjective and the result follows by the exactness
of the second row. $\Box$

\bs

\noi {\bf Remark.} 
It is reasonable to conjecture that
$H^{n+1}(\cG,\cK)\rightarrow
H^{n+1}(\tE,\mu_d)$ is injective. But we need to understand
$\cK$ better.  This question will be studied in future work.

\bs

The section closes by looking at the case where $s=2$.
We
have $(1+\sigma)(1-\tau)=
(1+\sigma)(1-\sigma\tau)=(1-\tau)(1-\sigma\tau)$ and we find
\begin{eqnarray*}
C_i&=&\tau^i(1+\sigma)(1-\tau) = \tau^i(1-\sigma)(1-\sigma\tau)
= \tau^i(1-\tau)(1-\sigma\tau)\\
\cB&=& (1-\sigma)\tau^{\frac{d+1}{2}} = 
\tau^{\frac{d+1}{2}}-\tau^{-\frac{d+1}{2}}\sigma = \tau^{\frac{d+1}{2}} 
(1-\tau^{-1}\sigma)=\tau^{\frac{d+1}{2}} 
(1-\sigma\tau).
\end{eqnarray*}
From  this we find that $M_3={\Z}[J]\cdot(1-\sigma\tau)=
{\Z}[J]\cdot\cB$.
Looking at $M_3$ in this way may make what is going on when 
$s=2$ more transparent (in particular the relationship to $E'$.)
Even more, we have noted earlier that  both 
${\Z}[J]\cdot(1\pm\sigma\tau)$ have ${\Z}$-rank $d$,
and therefore as $ (1-\sigma\tau)(1+\sigma\tau)=0$
we know
the kernel of
the map $\cdot(1\pm\sigma\tau)$ is ${\Z}[J]\cdot(1\mp\sigma\tau)$.
 This leads to the following result.
 
\begin{lemma}\label{OldLemma14} When $s=2$ we have two exact sequences
\[
0\rightarrow M_3
\rightarrow{\Z}[G]
\rightarrow  {\rm Ind}^\cG_{\cH'}\rightarrow 0
\]
and
\[
0\rightarrow {\rm Ind}^\cG_{\cH'}
\rightarrow{\Z}[G]
\rightarrow M_3 \rightarrow 0.
\]
The second sequence coincides with the middle column of the diagram
of Lemma \ref{OldTheorem12} up to an automorphism of $M_3$
and therefore $\cK\cong {\rm Ind}^\cG_{\cH'}
$ in this case. 
\end{lemma}

\bs

\noi {\bf Proof.} We know that 
${\Z}[J]\cdot(1+\sigma\tau)\cong {\rm Ind}^\cG_{\cH'}$ and 
$M_3={\Z}[J]\cdot(1-\sigma\tau)$.  The exact sequences follow as
the kernel of
the map $\cdot(1\pm\sigma\tau)$ is ${\Z}[J]\cdot(1\mp\sigma\tau)$.
For the second statement, in Lemma 
\ref{OldTheorem12} the map 
${\Z}[G]\rightarrow M_3$ is multiplication $\cdot\cB$ where 
$\cB=\tau^{\frac{d+1}{2}} 
(1-\sigma\tau)$, whereas it is multiplication by $ (1-\sigma\tau)$ in the lemma.
However, multiplication by $\tau^{\frac{d+1}{2}}$ is an autormorphism of
$M_3$ so the result follows.
   $\Box$

\bs

In the dihedral case  ($s=2$), the exact sequence of Lemma \ref
{OldLemma14} and the long exact cohomology sequence
give the first column of the diagram in the proof of Theorem 
\ref{OldTheorem13},
\[
\cdots \rightarrow H^{n}(E',\mu_d) \rightarrow H^{n}(\tE,\mu_d)\rightarrow H^{n}(\cG,M_3)
\]
\[
\rightarrow H^{n+1}(E',\mu_d) \rightarrow H^{n+1}(\tE,\mu_d)\rightarrow H^{n+1}(\cG,M_3)\cdots.
\]
However, $[\tE:E']=2$ and $d$ is odd, so we know
$H^{n+1}(\cG,\cK)=H^{n+1}(E',\mu_d) \rightarrow H^{n+1}(\tE,\mu_d)$ is injective.
This gives the following application of Theorem 
\ref{OldTheorem13}.

\bs

\begin{theorem}\label{Dihedral Corestriction} In the dihedral case  ($s=2$)
the cohomological kernel $H^{n+1}(E/F,\mu_d)$ is given by
\[
(*_8) \ \ \ \ \ \frac{{H}^n(\tF,\mu_d)}{{\rm cor}_{\tE/\tF}H^n(\tE,\mu_d)}
\stackrel{\cong}{\rightarrow} {\rm ker}(H^{n+1}(F,\mu_d)
\rightarrow H^{n+1}(E,\mu_d)).
\]
\end{theorem}

\noi {\bf Proof.} Since
$H^{n+1}(\cG,\cK)=H^{n+1}(E',\mu_d) \rightarrow H^{n+1}(\tE,\mu_d)$ is injective the result is immediate by Theorem 
\ref{OldTheorem13}.
 $\Box$

\bs

\bs

\Large
\noi {\bf Acknowledgements.}
\normalsize

\bs

\noi The author was supported by the James B. Axe Foundation. The author was also supported by the University of California at Santa Barbara as a graduate student, while many of these ideas were developed. These results appeared in the author's PhD dissertation. The author is very grateful to Bill Jacob for guidance and support during every step of the development of this paper.

\bibliographystyle{amsalpha}

\end{document}